\documentclass[journal]{IEEEtran}

\usepackage{commands}
\usepackage{color}
\usepackage{hyperref}
\usepackage{algorithmicx}
\usepackage[ruled]{algorithm2e}

\title{Nearest Neighbors for Matrix Estimation Interpreted as Blind Regression for Latent Variable Model}

\author{Yihua~Li,
        Devavrat~Shah,
        Dogyoon~Song,
				Christina~Lee~Yu
\thanks{This work was performed while all authors were affiliated with the Laboratory for Information and Decision Systems and the department of EECS at Massachusetts Institute of Technology. Shah is a member and director of the Statistics and Data Science Center at MIT.}
\thanks{A preliminary version of this work was presented at the Neural Information Processing Systems
Conference in December 2016 under the name ``Blind Regression: Nonparametric Regression for
Latent Variable Models via Collaborative Filtering''. The results have been significantly improved,
strengthened, and expanded with new extensions since the preliminary version, and thus 
this manuscript has limited overlap with the preliminary version of this work.}}

\begin{document}

\maketitle

\begin{abstract}
We consider the setup of nonparametric {\em blind regression} for estimating the entries of a large $m \times n$ matrix, when provided with a small, 
random fraction of noisy measurements. We assume that all rows $u \in [m]$ and columns $i \in [n]$ of the matrix are associated to latent features 
$\featrow{u}$ and $\featcol{i}$ respectively, and the $(u,i)$-th entry of the matrix, $A(u, i)$ is equal to $f(\featrow{u}, \featcol{i})$ for a latent function $f$. Given noisy 
observations of a small, random subset of the matrix entries, our goal is to estimate the unobserved entries of the matrix as well as to ``de-noise'' 
the observed entries. As the main result of this work, we introduce a nearest-neighbor-based estimation algorithm, and 
establish its consistency when the underlying latent function $f$ is Lipschitz, the underlying latent space is a bounded diameter Polish space,
%
and the random fraction of observed entries in the matrix is at least $\max \big( m^{-1 + \delta}, n^{-1/2 + \delta} \big)$, for any $\delta > 0$. 
As an important byproduct, our analysis sheds light into the performance of the classical collaborative filtering algorithm for matrix completion, 
which has been widely utilized in practice.  Experiments with the MovieLens and Netflix datasets suggest that our algorithm provides a principled 
improvement over basic collaborative filtering and is competitive with matrix factorization methods. Our algorithm has a natural extension to 
the setting of tensor completion via flattening the tensor to matrix. 
When applied to the setting of image 
in-painting, which is a $3$-order tensor, we find that our approach is competitive with respect to state-of-art tensor completion algorithms across 
benchmark images.
\end{abstract}

\begin{IEEEkeywords}
Blind regression, matrix estimation, matrix completion, tensor estimation, tensor completion, latent variable model, collaborative filtering, nearest neighbor methods
\end{IEEEkeywords}

%


\section{Introduction} 

The problem of matrix completion has received enormous attention in the past decade: 
consider an $m \times n$ matrix $A$ of interest. Suppose we observe a subset of the entries of an $m \times n$ matrix $Z$,  which is a noisy version of $A$, such that each $(u,i)$-th entry $Z(u,i)$ is a random variable with $\E[Z(u,i)] = A(u,i)$ for $u \in [m], ~i \in [n]$\footnote{We shall utilize notation $[m] = \{1,\dots m\}$.}. 
The goal of matrix completion is to recover matrix $A$ given partial observations from $Z$. 


\subsection{Our Contributions}

We provide a similarity-based nearest neighbor algorithm akin to popular collaborative filtering with theoretical performance guarantees under the latent variable model. 
In addition, we extend the analysis of our algorithm to the setting of tensor completion through flattering tensor to matrix. To our knowledge, this is the 
first theoretical analysis for similarity-based collaborative filtering algorithms, shedding insight into the widespread success of this popular heuristic for 
the past two decades. The algorithm we introduce is a simple variant of classical collaborative filtering, in which we compute similarities between pairs of 
rows and pairs of columns by comparing their common overlapped entries. Our model assumes that each row and column is associated to a latent variable, 
i.e. a vector of hidden features, and that the data entry is in expectation equal to some unknown function of those latent variables. 
We assume the latent space is a complete, separable metric space aka Polish space equipped with a Borel probability measure. 
The key regularity condition that we require is that the image of the latent space by the latent function has a small effective covering 
number with respect to the push-forward measure\footnote{Let $(\cX, d, \mu)$ denote the latent space. Let 
$B(x, r) = \{ z \in \cX: d(x,z) \leq r \}$. Given $\eps, r > 0$, the effective covering number of $\cX$ with respect to $\mu$
\[	N_{\textrm{eff}}(\cX, r, \eps) \triangleq \inf_{I, S \subset \cX}\left\{ |I| \text{ s.t. } S \subset \cup_{x \in I} B(x, r) 
	\text{ and }\mu(\cX \setminus S) \leq \eps \right\}.	\]}
which Borel measure over Polish space naturally satisfies. 

Given this latent variable model, we prove that the estimate produced by this algorithm is consistent as long as the fraction of entries that are observed is at least $\max(m^{-1 + \delta}, n^{-1/2 + \delta})$ for some $\delta > 0$ for an $m \times n$ matrix (for precise statement, See Corollary \ref{coro:polish} and its implication). We provide experiments using our method to predict ratings in the MovieLens and Netflix datasets. The results suggest that our algorithm improves over basic collaborative filtering and is competitive with factorization-based methods.

We also discuss that the algorithm and analysis can be extended to tensor completion by flattening the tensor to a matrix.
We implemented our method for predicting missing pixels in image in-painting, which showed that our method is competitive with existing spectral methods used for tensor completion.

The algorithm that we propose 
has similarities to classical non-parametric nearest neighbor method, cf. \cite{chenshah} and
kernel regression, which also relies on approximations by local smoothing, cf. \cite{mack1982weak, WandJones94}. However, since kernel regression and other similar methods use explicit knowledge of the input features, their analysis and proof techniques do not extend to our context. Instead of using distance in the unknown latent space, the algorithm weighs datapoints according to similarities that are computed from the data itself. Our analysis shows that although the similarities between the data points may not reflect the distance between the latent features, they essentially reflect the functional distances (in the $L^2$ sense) between the latent function restricted to the pair of rows (or columns) associated with the data points, which is sufficient to guarantee that the datapoints with high similarities are indeed similar in value.


\subsection{Related Literature}

The primary methods used to solve the problem in the literature include neighbor-based approaches, such as collaborative filtering, and spectral approaches, which include low-rank matrix factorization or minimization of a loss function with respect to spectral constraints. 

\subsubsection*{Spectral Methods}

In the recent years, there have been exciting intellectual developments in the context of spectral approaches such as matrix factorization. All matrices admit a singular-value decomposition, such that they can be uniquely factorized. The goal of the factorization-based method is to recover row and column singular vectors accurately from the partially observed, noisy matrix $Z$ and subsequently estimate the matrix $A$. \cite{srebro2004generalization} was one of the earliest works to suggest the use of low-rank matrix approximation in this context. Subsequently, statistically efficient approaches were suggested using optimization-based estimators, proving that matrix factorization can fill in the missing entries with sample complexity as low as $rm \log m$ for an $m \times m$ matrix, where $r$ is the rank of the matrix \cite{candes2009exact, rohde2011estimation, keshavan56matrix, negahban2012restricted, jain2013low}. There has been an exciting line of ongoing work to make the resulting algorithms faster and scalable \cite{Fazel03, LiuVandenberghe10, Cai08, Lin09, Shen09, Mazumder10}.

\cite{XuMassoulieLelarge14} proposed a spectral clustering method for inferring the edge label distribution for a network sampled from a generalized stochastic block model. The model is similar to the proposed latent variable model introduced in Section \ref{sec:setup}, except that the edges are labeled by one of finitely many labels in a symmetric setup with $m = n$, 
and the goal is to estimate the label distribution in addition to the expected label. When the expected function has a finite spectrum decomposition, i.e. low rank, then they provide a consistent estimator for the sparse data regime, with $\Omega(m \log m)$ samples. When the function is only approximately low rank (e.g. the class of general Lipschitz functions), for a fixed rank $r$ approximation, the mean squared error bounds converge to a positive constant which captures the low rank approximation gap.
{That is, $\Omega(m \log m)$ samples are not sufficient to guarantee consistent estimation for the entire class of Lipschitz functions. }  
    
Many of these approaches are based on the structural assumption that the underlying matrix is {\em low-rank} and the matrix entries are reasonably ``incoherent''. Unfortunately, the low-rank assumption may not hold in practice. The recent work \cite{GantiBalzanoWillett2015} makes precisely this observation, showing that a simple non-linear, monotonic transformation of a low-rank matrix could easily produce an effectively high-rank matrix, despite few free model parameters. They provide an algorithm and analysis specific to the form of their model, which achieves sample complexity of $O((mn)^{2/3})$ for an $m \times n$ matrix. However, their algorithm only applies to functions $f$ which are a nonlinear monotonic transformation of the inner product of the latent features. 
\cite{LeeKimLebanonSingerBengio16} propose an algorithm for estimating locally low rank matrices, however their algorithm assumes prior knowledge of the ``correct'' kernel function between pairs of rows and columns which is not known a priori.

\cite{Chatterjee15} proposes the universal singular value thresholding estimator (USVT) inspired by low-rank matrix approximation. Somewhat interestingly, it argues that under the latent variable model considered in this work (see Section \ref{sec:setup}), the USVT algorithm provides an accurate estimate for any Lipschitz function. However, to guarantee consistency of the USVT estimator for an $m \times m$ (i.e. $m = n$) matrix, it requires observing $\Omega\left(m^{\frac{2(d+1)}{(d+2)}} \right)$ many entries out of the $m^2$ total entries, where $d$ is the dimension of the latent space in which the row and column latent features belong. In recent work, \cite{Xu2017} extends the analysis of USVT for graphon estimation, assuming a generative latent variable model for binary observation matrices representing networks. If the latent function is $\alpha$-Holder smooth, author establishes that the spectrum decays polynomially, and thus the MSE of the USVT estimator is bounded above by $O \big((mp)^{-\frac{2\alpha}{2\alpha+d}} \big)$, which converges to zero as long as $p = \omega(m^{-1})$ (for an $m \times m$ symmetric matrix).




\subsubsection*{Collaborative Filtering}

The term collaborative filtering was coined by \cite{goldberg92}, and this technique is widely used in practice due to its simplicity and ability to scale.
There are two main paradigms in neighborhood-based collaborative filtering: the user-user paradigm and the item-item paradigm. To recommend items to a user in the user-user paradigm, one first looks for similar users, and then recommends items liked by those similar users. In the item-item paradigm, in contrast, items similar to those liked by the user are found and subsequently recommended. Much empirical evidence exists that the item-item paradigm performs well in many cases \cite{Linden03,korenHandbook,Ning2015}. There have also been many heuristic improvements upon the basic algorithm, such as normalizing the data, combining neighbor methods with spectral methods, combining both user and item neighbors, and additionally optimizing over interpolation weights given to each datapoint within the neighborhood when computing the final prediction \cite{BellKoren07,Koren08,WangDeVriesReinders06}.

Despite the widespread success of similarity-based collaborative filtering heuristics, the theoretical understanding of these method is very limited. In recent works, latent mixture models have been introduced to explain the collaborative filtering algorithm as well as the empirically observed superior performance of item-item paradigms, c.f. \cite{BreslerChenShah14, BreslerShahVoloch15}.
However, these results assume binary ratings and a specific parametric model, such as a mixture distribution model for preferences across users and movies. We hope that by providing an analysis for collaborative filtering within a nonparametric model, we can provide a better understanding of collaborative filtering.

Within the context of dense graphon estimation, when the entries in the data matrix are binary and the sample probability $p$ is constant $\Theta(1)$, there have been a few theoretical results that prove convergence of the mean squared error for similarity-based methods \cite{AiroldiCostaChan13, ZhangLevinaZhu15}. They hinge upon computing similarities between rows or columns  by comparing commonly observed entries, similar to collaborative filtering. Similar to our result, they are able to prove convergence for the class of Lipschitz functions. However, \cite{AiroldiCostaChan13} assumes that the algorithm is given multiple instances of the sampled dataset, which is not available in our formulation. \cite{ZhangLevinaZhu15} is weaker than our result in that it assumes $p = \Theta(1)$, however they are able to handle a more general noise model, when the entries are binary. The similarity between a pair of vertices is computed from the maximum difference between entries in the associated rows of the second power of the data matrix, which is computationally more expensive than directly comparing rows in the original data matrix.

\subsubsection*{Tensor Completion}
Recently there have been efforts to extend decomposition methods or neighborhood-based approaches to the context of tensor completion, however this has proven to be significantly more challenging than matrix completion due to the complication that tensors do not have a canonical decomposition such as the singular value decomposition (SVD) for a matrix. This property makes obtaining a decomposition for a tensor challenging. The survey \cite{kolda2009tensor} elaborates on these challenges. There have been recent developments in obtaining efficient tensor decompositions in form of rank-1 tensors (tensors obtained from one vector), presented in \cite{anandkumar2014tensor}. This has been especially effective in learning latent variable models and estimating missing data as shown in, for example \cite{jain2014provable, oh2014learning}.

Many results in tensor estimation take the approach of flattening the tensor to a matrix and subsequently apply matrix estimation algorithms \cite{liu2013tensor,gandy2011tensor,tomioka2010estimation,tomioka2011statistical}. A $d$-order tensor where each dimension is length $n$ would be flattened to a $n^{\lfloor d/2 \rfloor} \times n^{\lceil d/2 \rceil}$ matrix, resulting in a sample complexity of $\Omega(n^{\lceil d/2 \rceil} \text{polylog}(n) )$. 
Subsequently there has been a line of work extending spectral methods, local iterative methods, or the sum of squares method to the specific tensor structure to obtain improved sample complexities of $\Omega(n^{d/2} \text{polylog}(n) )$ \cite{jain2014provable, bhojanapalli2015new, BarakMoitra16, PotechinSteurer17, xia2017statistically, xia2017polynomial, MontanariSun18}. 

Beyond tensor decomposition, there have been recent developments in the context of learning latent variable models or mixture distributions also called non-negative matrix factorization, c.f. \cite{arora2012learning, arora2012computing}.



\subsection{Organization of the Paper}

In Section \ref{sec:setup} we setup the formal model and problem statement and discuss the assumptions needed for our analysis. In Section \ref{sec:estimator} we introduce the basic form of our algorithm, which is similar to the user-user variant of collaborative filtering. We present heuristic variants of the algorithm that perform well in practice. In Section \ref{sec:main.result} we present the main theoretical results of our paper as they pertain to matrix completion, showing provable convergence of the user-user (and by symmetry item-item) variant of our algorithm. In Section \ref{sec:discussion} we provide a discussion of our results and comparison with other models. In Section \ref{sec:tensor_results} we discuss how to extend the algorithm and analysis to tensor completion. In Section \ref{sec:exp} we present experimental results from applying our methods to both matrix completion in the context of predicting movie ratings, and tensor completion in the context of image inpainting. The detailed proofs are presented in the Appendix.
\section{Setup}\label{sec:setup}


\subsection{Our Model}\label{sec:model}
Suppose that there is an unknown $m \times n$ matrix $A$ which we would like to estimate. We observe only a fraction 
of the total $mn$ entries of $A$ with some noise added. Let $\obs \subset [m] \times[n]$ denote the index set of observed entries. 
Specifically, we observe entries of data matrix $Z$ that is generated as follows. Let $M \in \{0,1\}^{m \times n}$ be a binary matrix, which we call the masking matrix. We let $A \in \Reals^{m \times n}$ 
denote the signal matrix and $N \in \Reals^{m \times n}$ denote the noise matrix. For each $(u,i) \in [m] \times [n]$, 
\begin{equation}\label{eqn:gen_model}
	Z(u,i) = \begin{cases}
			A(u,i) + N(u,i)	&	\text{when }M(u,i) = 1,\\
			\text{unknown}	&	\text{when }M(u,i) = 0.
		\end{cases}
\end{equation}
For later use, we let $\Omega \triangleq \{ (u,i) \in [m] \times [n]: M(u,i) = 1 \}$ denote the set of index pairs of the observed entries. 

\subsubsection{Latent Variable Model}\label{sec:lv_model}
We assume $A(u,i)$ is generated by the following model.
\begin{itemize}
	\item
	Nonparametric model: there exists a latent function $f$ such that
	\begin{equation}\label{eqn:latent_model}
		A(u,i) = f( \featrow{u}, \featcol{i})	
	\end{equation}
	for all $(u,i) \in [m] \times [n]$. Here, $\featrow{u}, \featcol{i}$ denote latent variables associated with row $u$ 
	and column $i$, respectively.
	\item
	Regularity Assumptions:
	\begin{itemize}
		\item
		For all $u \in [m]$, $\featrow{u} \in \latsprow$, where $(\latsprow, \drow)$ is a complete, separable metric space a.k.a. Polish space equipped with a Borel probability measure $\murow$ and $\featrow{u}$ 
		is drawn i.i.d. according to $\murow$.
		\item
		For all $i \in [n]$, $\featcol{i} \in \latspcol$, where $(\latspcol, \dcol)$ is a complete, separable metric space aka Polish space equipped with a Borel probability measure $\mucol$ and $\featcol{i}$ 
		is drawn i.i.d. according to $\mucol$.
		\item
		Latent function $f$ is bounded, i.e. for all $\alpha \in \latsprow$ and $\beta \in \latspcol$,
$	\big| f( \alpha, \beta) \big| \leq \fbound.$ 	
\item Without loss of generality, we shall assume that there exists $\alpha \in \latsprow$ and $\beta \in \latspcol$
such that $f(\alpha, \beta) = 0$.
\end{itemize}

Our results will additionally require that the local measure is well-behaved such that there is a sufficient mass of nearest neighbors. In particular we introduce two concrete models that have good local neighborhood properties:
\begin{enumerate}
\item[(a)] {\em Finite Types:}
Let $\latsprow$ be equipped with the discrete metric\footnote{$\drow(x_1, x_2) = 1$ if and only if $x_1 \neq x_2$.} topology and suppose $\murow$ has finite support in $\latsprow$ with $\supp(\murow)$ denoting the support of $\murow$.
\item[(b)] {\em Lipschitz Latent Function:}
Assume that the function $f$ is $\lip$-Lipschitz in the sense that 
\begin{align*}
	&\big| f( \alpha_1, \alpha_1) - f( \beta_1, \beta_2) \big| \\
		&\qquad\leq \lip \max \big( \drow(\alpha_1, \beta_1), \dcol(\alpha_2, \beta_2) \big).
\end{align*}
\end{enumerate}

\end{itemize}

\subsubsection{Noise}
We assume $N(u,i)$ is a centered, sub-gaussian random variable for all $(u,i) \in [m] \times [n]$ such that
\begin{itemize}
	\item
	Centered: $\bbE N(u,i) = 0$ for all $(u,i) \in [m] \times [n]$.
	\item
	Sub-gaussian\footnote{The Orlicz $\psi_2$-norm of a random variable $X$ is defined as $\| X \|_{\psi_2} \triangleq \inf_{t > 0} 
		\{ \Exp{\exp\big(\frac{|X|^2}{t^2}\big)} \leq 2 \}$.}: there exists $\sigma > 0$ such that $\| N(u,i) \|_{\psi_2} \leq \sigma$ 
		for all $(u,i) \in [m] \times [n]$.
	\item
	Independent and identically distributed (i.i.d.): $N(u,i)$'s are independent and identically distributed.
\end{itemize}

\subsubsection{Masking}
We assume $M$ is a random matrix with each entry drawn as per $\textrm{Bernoulli}(p)$ for some $p \in (0,1]$, i.i.d. That is, for each 
$(u,i) \in [m] \times [n]$, 
\[	
	M(u,i) = \begin{cases}
			1	&	\text{with probability }p,\\
			0	&	\text{with probability }1-p.
			\end{cases}	
\]

Our algorithm will involve one threshold parameter which will need to be tuned. As a result of our analysis we can specify the optimal choice of the threshold as a function of $n,m,p,\sigma, L, \fbound$ which trades off between the bias and variance. In practice, all of these parameters are assumed to be unknown and the threshold parameter can be chosen via cross-validation.


\subsection{Problem Statement: Blind Regression}
Note that the latent variable model representation in \eqref{eqn:latent_model} is not the unique representation, as there exists multiple
equivalent representations. Suppose that one applies a measure-preserving transformation $T$ on the latent feature space 
$\latsprow \times \latspcol$, and take the push-forward of $f$ with respect to $T$ ($f \circ T^{-1}$) as the new latent function. 
This new representation -- the pair of the latent space and the latent function -- yields the same data generation process. Therefore, 
the question of estimating the function $f$ itself is not well posed, and we focus our energy on predicting the values $A(u,i)$. 

\begin{problem}[Blind Regression]\label{prob:blind_regression}
	Given $Z \in \Reals^{m \times n}$ that is partially observed on $\Omega \subset [m] \times [n]$ as described in Section \ref{sec:model}, 
	we want to estimate the underlying matrix $A$.
\end{problem}
We call the problem of interest {\em Blind Regression} for the following reason. In the setting of {\em Regression}, one observes data containing features and associated labels; the goal is to learn the functional relationship (or model) between features and labels assuming that labels are noisy observations. In our setting, tuples $(\featrow{u}, \featcol{i})$ are the relevant (but unobserved) features and $Z(u,i)$ are noisy observations of associated labels $A(u,i) = f(\featrow{u},\featcol{i})$. We want to predict the value of $f(\featrow{u},\featcol{i})$ for all pairs $(u,i)$ for $u \in [m], i \in [n]$. In the sense that we want to predict the function $f$ evaluated on new points given previous data, the task has the feel of Regression. However, the features $(\featrow{u}, \featcol{i})$ are {\em latent}; and thus we use the term {\em Blind Regression}. 

Given an estimator $\hat{A}$ for the unknown matrix $A \in \Reals^{m \times n}$ of interest, we use the mean-squared error (MSE) to evaluate the performance of the estimator, defined as 
\begin{equation}\label{eqn:MSE}
    \MSE(\hat{A}) = \Exp{ \frac{1}{mn}\sum_{u=1}^m \sum_{i=1}^n \left( \hat{A}(u,i) - A(u,i)\right)^2 }.
\end{equation}
The expectation here is taken over all sources of randomness in the data generation process: (i) realization of the latent variables; 
(ii) realization of the noise variables; and (iii) masking. 
An estimator $\hat{A}$ is called consistent if $\lim_{m, n \to \infty} \MSE(\hat{A}) = 0$. For a consistent estimator, we also want to 
establish an upper bound on the rate of convergence rate for the mean-squared error. Now we pose a follow-up question.
\begin{problem}
	Can we achieve a consistent estimator $\hat{A}$ for $A$ in the setup of Problem \ref{prob:blind_regression}? 
	If so, how fast does $\MSE(\hat{A})$ decay to $0$ for given $p = p(m,n)$ as the problem size $m, n \to \infty$? 
\end{problem}


\subsection{Exchangeability and Latent Variable Model}
The latent variable model is well motivated and arises as a canonical representation for row and column exchangeable data, cf.  \cite{Aldous81} and  \cite{Hoover82}. 
Suppose that our data matrix $Z$ is a particular realization of the first $m \times n$ entries of a random array $\bZ = \left\{\bZ(u,i) \right\}_{(u, i) \in \Nats \times \Nats}$, which satisfies
\begin{align}\label{eq:exchangeable}
 \bZ(u,i) \overset{d}{=} \bZ\left(\sigma(u), \tau(i)\right) ~\text{ for all } (u,i),
\end{align}
for every pair of permutations\footnote{The permutations over $\Nats$ are defined in the usual manner where only finitely many indices are {\em permuted}.} $\sigma, \tau$ 
of $\Nats$. We use $\overset{d}{=}$ to denote that the joint distribution over $\{Z(u,i)\}_{u,i}$ is equivalently distributed as the joint distribution over $\{\bZ\left(\sigma(u), \tau(i)\right)\}_{u,i}$, i.e. the random variables on both sides have the same distribution. Random array $\bZ$ satisfying 
\eqref{eq:exchangeable} is called exchangeable\footnote{To be precise, {\em separately} row and column exchangeable.}. 
For an interested reader, \cite{Austin12} and \cite{OrbRoy2015} present overviews of exchangeable arrays.  

In practice, the use of exchangeable arrays as a model is appropriate for variety of reasons. For example, in the setting of a recommendation system with anonymized data, 
this property may be reasonable if the order of the users in the system does not intrinsically carry information about the type of user; or in other words, if a user in the system 
could equally likely have been located in any row of the dataset. 

In addition to exchangeability being quite a reasonable property for a wide variety of applications, it also leads to a convenient latent variable representation. The Aldous-Hoover representation theorem provides a succinct characterization for such exchangeable arrays. According to the theorem (see Corollary 3.3 in \cite{OrbRoy2015} for example), a random 
data array $\bZ$ is exchangeable if and only if it can also be represented as
\begin{align} \label{eq:exch_form}
 \bZ\left( u,i \right) \overset{d}{=} f_{\theta}\big( \frow(u), \fcol(i), \theta_{\text{entry}}(u,i) \big) ~\text{ for all }(u,i)
\end{align}
where $\theta$, $\big\{\frow(u)\big\}_{u \in \Nats}$, $\big\{\fcol(i)\big\}_{i \in \Nats}$, $\big\{\theta_{\text{entry}}(u,i)\big\}_{(u,i) \in \Nats \times \Nats}$ are independent random variables 
drawn uniformly at random from the unit interval $[0,1]$, and $f_{\theta}$ is a measurable function indexed by the realization of $\theta$. 
As described in \cite{OrbRoy2015}, this suggests the following generative model: 
\begin{enumerate}
	\item	Sample an instance of $\theta \sim U[0,1]$ determining the governing function $f_{\theta}$.
	\item For every row $u \in [m]$ and every column $i \in [n]$, independently sample uniform random variables 
		$\frow(u) \sim U[0,1]$, $\fcol(i) \sim U[0,1]$, $\theta_{\text{entry}}(u,i) \sim U[0,1]$.
	\item Compute the realized data matrix $Z$ according to 
	\begin{align*}
	Z(u,i) = f_{\theta}\big( \frow(u), \fcol(i), \theta_{\text{entry}}(u,i) \big).
	\end{align*}
\end{enumerate}

By comparing the model from \eqref{eq:exch_form} with the latent variable model described in Section \ref{sec:lv_model} (cf. \eqref{eqn:gen_model}, \eqref{eqn:latent_model}), 
we can see that the latent variable model considered in this work is a restricted subclass of exchangeable models that additionally impose an additive noise model and 
regularity conditions on the function $f_{\theta}$ in exchange for a more general latent space and associated probability measure.

\medskip
\paragraph*{Equivalence of Models}
In our model we have conditioned on the universal index\footnote{Equivalently, we may think that $f_{\theta} = f$ for all $\theta \in [0,1]$. Note that we do not know what $f$ is a priori.} 
$\theta$, such that given partial observations from matrix $Z$ for a particular $f_{\theta}$, our goal is to learn predicted outcomes of the realized $f_{\theta}$. Our model takes 
the form of
\begin{align*}
 	Z(u,i) =  f(\featrow{u}, \featcol{i}) + \eta(u,i)
\end{align*}
where $\{\featrow{u}\}_{u \in [m]}, \{\featcol{i}\}_{i \in [n]}, \{\eta(u,i)\}_{(u,i) \in [m] \times [n]}$ are sampled independently. We can transform this to the form of \eqref{eq:exch_form} by 
considering $f$ to be equal to the realized function $f_{\theta}$, considering the latent variables $\featrow{u} \sim \prxa$ and $\featcol{i} \sim \prxb$ to be higher dimensional 
representations of $\frow(u)$ and $\fcol(i)$ in spaces $\latspA$ and $\latspB$, and considering the noise term $\eta_{\text{entry}}(u,i)$ to be generated by applying some transformation 
to the variable $\theta(u,i)$. Given these transformations, it becomes equivalent that 
\begin{align*}
 	Z(u,i) &= f_{\theta}\big( \frow(u), \fcol(i), \theta_{\text{entry}}(u,i)  \big)	\\
		 &=  f(\featrow{u}, \featcol{i}) + \eta(u,i).	
\end{align*}
	Instead of allowing $f$ to be any arbitrary measurable function over $[0,1] \times [0,1]$,
our model additionally imposes regularity conditions on $f$ by requiring it to be bounded and either finite types or Lipschitz continuous with respect to a higher dimensional representation 
$\latspA \times \latspB$. From a modeling perspective, we are effectively transferring 
the {\em model complexity} from a potentially complex measurable latent function over $[0,1] \times [0,1]$ to a simpler (e.g. Lipschitz) latent function over a potentially 
more complex latent variable space $\latspA \times \latspB$. The simple functional form provides analytic tractability for establishing theoretical results.


\subsection{Comparison with Other Models}

In this section we discuss how our model relates to other models considered in the literature. In particular we would like to clarify which models are captured in the latent variable model and which are not. The key assumptions of our model are (1) bounded or Lipschitz latent function, (2) additive i.i.d. sub-gaussian noise, and (3) latent space being finite or more generally bounded Polish space. 

\subsubsection{Low-rank Assumption}
Suppose that the latent spaces $\latsprow$ and $\latspcol$ are equal dimensional vector spaces and the latent function $f$ is a bilinear form. This is a low-rank model with its rank being 
equal to $\dim \latsprow = \dim \latspcol$. 
If we assume $\latsprow$ and $\latspcol$ have bounded diameters and are equipped with some metric topology, then this model is a specific case of our model with Lipschitz constant being $\big(\diam \latsprow + \diam \latspcol \big)$ (up to a multiplicative constant). Our latent variable model additionally assumes that the latent feature vectors are randomly sampled i.i.d. according to some probability measure, whereas typical low rank results allow for arbitrarily chosen latent feature vectors that satisfy incoherence.

Note that we utilize 
a different measure of `model complexity' than the rank of the parameter matrix. 
When the underlying function is truly bilinear and $\min \big\{ \dim \latsprow, \dim \latsprow \big\} \ll \diam \latsprow + \diam \latsprow$, then rank of the data matrix could 
be a better measure to capture the essential complexity of the model; note that the rank of the data matrix does not scale as the diameter of the latent space increases. 
However, if the underlying model is nonlinear, the data matrix is likely to have full rank, e.g., when $f(\featrow{u}, \featcol{i}) = g(\left< \featrow{u}, \featcol{i}\right>)$ for 
some monotone increasing nonlinear function $g$. In particular, when the size of the matrix far exceeds the `intrinsic model complexity' (and hence so does the rank 
of the data matrix), it would make more sense to use our model. 

Suppose that there is an underlying low rank matrix, yet the observation is an entrywise nonlinear monotone transformation of the low rank matrix; this is also known as single index models. \cite{GantiBalzanoWillett2015} pointed out that adding nonlinearity could easily result in the underlying matrix no longer being low rank, requiring more complex approaches to estimate. In contrast, our model very easily handles nonlinearities as long as they satisfy local smoothness conditions.

\subsubsection{Biclustering, Submatrix Detection, Stochastic Block Model}
The models used in biclustering, submatrix detection, and planted clique assume that there is some submatrix for which the data is significantly shifted in expectation compared to the rest of the matrix which is assumed to be uniform plus noise \cite{balakrishnan2011statistical,ma2015computational}. This would correspond to ``finite types'' in our model, which can be modeled as a piecewise constant function. The desired goal is to detect or identify the deviant submatrix corresponding to a subset of the rows and columns. In contrast our task focuses on estimating the expected matrix. Given an estimation algorithm that could guarantee max error bounds entrywise or for each row/column, we could simply threshold to obtain an estimate for the deviant submatrix. Another distinction is that the literature in submatrix detection focuses on understanding thresholds of the minimal size submatrix that is detectable. However, in our setting as the latent feature vectors are sampled from a fixed distribution, the size of the deviant submatrix that our model studies would always be proportional to $O(mn)$ as long as the probability of sampling a deviant row or column is bounded below by a fixed constant. The stochastic block model also assumes finite types modeled through a piecewise constant function, which would fit within our latent variable model. It is commonly used to study the task of clustering, where the goal is to recover a true underlying partitioning of the rows and columns rather than only estimate the expected matrix \cite{abbe2016exact}. Their setting specifically assumes a binary observation model, while our model assumes identically distributed additive sub-gaussian noise. 

\subsubsection{Strong Stochastic Transitivity, Statistical Seriation, Low Permutation-Rank Matrices}
There have been many models introduced to study matrices that arise from pairwise comparisons for the purposes of ranking or estimation. Many parametric models such as Bradley-Terry-Luce or Thurstone, can be described by a smooth function of scalar latent row and column variables. Our model adds the additional assumption that the latent variables must be sampled i.i.d. from a bounded diameter subspace. There have also been a series of nonparametric models introduced, in particular the strong stochastic transitivity assumption, which only requires that there exists some permutation of the rows and columns such that the entries are monotonically increasing along the permuted rows and columns \cite{shah2017stochastically, chatterjee2019estimation}. 
Statistical seriation assumes that there is some permutation of the rows and columns such that after the permutation, all rows and columns have the same shape, which could include monotonicity or other shape constraints \cite{flammarion2019optimal}. 
Low permutation-rank matrices assume that a mixture model over strong stochastic transitive matrices, i.e. a mixture over permutation matrices \cite{shah2018low}.

Our latent variable model can encompass a subset of such matrices by imposing monotonicity assumptions on our latent function. However our model additionally assumes latent
space being Polish space endowed with Borel measure, which are not required for permutation matrices. Similarly, a subset of mixtures of permutation matrices can be modeled 
by latent functions which are a mixture of different monotonic functions; however this may not be able to fully encompass all low permutation-rank models. If one were to impose 
boundedness on the entries, then it is possible that as the size of the matrix grows, strong stochastic transitivity might imply the existence of good nearest neighbor rows, 
which would then allow our algorithm to also perform well, but that remains an important direction to explore.

%
%


\section{Algorithm} \label{sec:estimator} 
	
Our algorithm builds on intuition from local approximation methods such as kernel regression. Therefore it takes the form of a similarity-based method, 
which first defines a kernel, i.e. similarity between pairs of rows or columns, and then computes the estimate for each matrix entry by averging over 
datapoints that are determined to be close to the entry of interest according to the `similarity.' 

We present the general form of our algorithm in Section \ref{sec:general_alg}. In Section \ref{sec:user_user}, we describe the basic user-user 
fixed radius nearest neighbor algorithm\footnote{This is equivalent to a variant of the classical similarity based collaborative filtering methods.} as 
an example of the algorithm with a concrete choice of similarity and the averaging scheme to define the estimate. In Section \ref{sec:other_algs}, 
we describe other variations of our algorithm, e.g., a variation of the algorithm that combines both row and column similarities to compute 
the kernel between datapoints. 

The user-user nearest neighbor algorithm serves as our prime example and we provide theoretical guarantees for a vanishing upper bound on 
its mean squared error later in Section \ref{sec:main.result}. Also, we show experimental results that suggest combining row and column similarities improve 
the quality of estimates; see Section \ref{sec:exp}.


\subsection{General Form}\label{sec:general_alg}
Our algorithm takes (i) the data matrix $Z$, (ii) a rule for sample splitting, and (iii) a thresholding parameter $\eta \geq 0$ as its inputs. 
The algorithm outputs $\hA \in \Reals^{m \times n}$. 

\begin{algorithm} [h] 
    \SetKwInOut{Input}{Input}
    \SetKwInOut{Output}{Output}
    \Input{$Z \in \Reals^{m \times n}; \eta \geq 0$}	
    \Output{$\hA \in \Reals^{m \times n}$}
    
    \begin{algorithmic}[1]
	\State Split the observations in $\Omega$. 
	For each $(u,i) \in \Omega$, put it in $\Omega_1$ or $\Omega_2$. 
	\State For each $(u,i) \in [m] \times [n]$, determine the set of reliable neighbors $ \baseest(u,i) \subset \Omega_2$ with 
	similarities/dissimilarities that are computed based on $\Omega_1$.
	\State For each $(u,i) \in [m] \times [n]$, compute $\hA(u,i)$ using $\baseest(u,i)$.
   \end{algorithmic}

    \caption{Generic Description of the Algorithm}\label{alg:main_algorithm}
\end{algorithm}

Specifically, our algorithm determines the set of reliable neighbors by considering the `behavioral' similarity in the function values. 
We need to define a similarity/dissimilarity statistic that can capture such similarities as well as can be computed from the data matrix. 
In Section \ref{sec:user_user}, we describe a version of our algorithm that utilizes the squared $\ell_2$ distance as a measure of 
dissimilarity, cf. Algorithm \ref{alg:user_algorithm}. This version of algorithm is the prime example we consider in this work and we provide 
a theoretical guarantee on its performance in Section \ref{sec:main.result}.

We remark that some similarity functions can be preferred over others, depending on the model assumptions. For example, the squared 
$\ell_2$ distance is a natural choice when the latent space is assumed to be Euclidean, while the cosine similarity (= the angle between 
the latent variables) would be a more faithful measure of similarity when the latent space is the projective space (= spherical). We briefly 
discuss this matter in Section \ref{sec:other_algs} with additional examples of variations of our algorithm.

\begin{remark}
There are many variations in (1) how we determine the set of reliable neighbors (e.g., by defining similarity/dissimilarity functions), 
and (2) how we compute the final estimate based on the observations at those neighbors.
\end{remark}

\begin{remark}
We note that sample splitting is done for the ease of analysis, and is not essential in executing the algorithm.
\end{remark}

\subsection{Prime Example: User-user Fixed Radius Nearest Neighbor (of the $0$-th order)}\label{sec:user_user}
In this section, we describe our prime example of the generic algorithm described in Algorithm \ref{alg:main_algorithm}, namely, the 
user-user fixed radius neighbor algorithm. In this version of algorithm,
\begin{itemize}
	\item we measure the dissimilarity between two rows with the squared $\ell_2$ distance between the rows at the overlapping column indices; 
		specifically, we define the dissimilarity between two rows $u, v \in [m]$ as 
		\begin{align}
			&\dissimrow(u,v)	= \frac{1}{ | \basesim(u) \cap \basesim(v) | }	\nonumber\\
				& \qquad\qquad
				\times  \sum_{j \in \basesim(u) \cap \basesim(v)} \big( Z(u,j) - Z(v,j) \big)^2,	\label{eqn:user_dissim}
		\end{align}
		with the convention $0/0 = \infty$, and 
	\item for each $(u,i) \in [m] \times [n]$, we estimate $\hA(u,i)$ by averaging $Z(v,i)$ for $v \in [m]\setminus \{u\}$ such that $M(v,1) = 1$ 
		and $v$ is similar to the row $u$. We define that $v$ is similar to $u$ if the squared $\ell_2$ distance between their rows is no greater than some threshold $\eta  \geq 0$, 
		which is a tunable parameter input to the algorithm.
\end{itemize}
See Algorithm \ref{alg:user_algorithm} for the full description. We state the algorithm for estimating a single entry $(u,i)$, which affects the sample splitting rule. 
The sample splitting is used for a cleaner analysis, but we believe that the results should extend without sample splitting as well.

\begin{algorithm} [h!] 
    \SetKwInOut{Input}{Input}
    \SetKwInOut{Output}{Output}
    \Input{$Z \in \Reals^{m \times n}; (u,i) \in [m] \times [n]; \eta \geq 0$}	
    \Output{$\hA(u,i)$}
    
    \begin{algorithmic}[1]
	\State Split the observations in $\Omega$: for each $(v,j) \in \Omega$, 
	
		\noindent put it in $\Omega_1$ if $j \neq i$, and put it in $\Omega_2$ otherwise.
	\State Determine the set of reliable neighbors using $\Omega_1$.
		\begin{itemize}
		\item
		For each $v \in [m]$, define 
		\[	\basesim(v) \triangleq \{ j \in [n]: (v,j) \in \Omega_1  \},	\]
		which is the adjacency set of $v$ denoting the observed columns for row $v$. These entries are 
		
		\noindent used to determine its similarity with other rows.
		\item
		For each $v \in [m]$, estimate the dissimilarity\\
		
		\noindent between two rows $u, v \in [m]$ as \eqref{eqn:user_dissim}; we let $\dissimrow(u,v) = \infty$ if $| \basesim(u) \cap \basesim(v) | = 0$.
		\item
		We consider $\{ v \in [m]: \dissimrow(u,v) \leq \eta \}$ as 
		
		\noindent the set of reliable neighbor rows of $u$.
		\item
		Define 
		\begin{align*}
			\baseest(u,i) &\triangleq \big\{ (v,j) \in [m] \times [n] \textrm{ such that}\\
				&\quad \dissimrow(u,v) \leq \eta ~\text{and}~ (v,j) \in \Omega_2  \big\}.
		\end{align*}
		\end{itemize}
	\State Compute the estimate for $(u,i)$ as
		\begin{equation*}
			\hA(u,i) =	\frac{1}{| \baseest(u,i) |} \sum_{(v,j) \in \baseest(u,i)} Z(v,j)
		\end{equation*}
		when $| \baseest(u,i) | \neq \emptyset$; we let $ \hA(u,i) = 0$ when $| \baseest(u,i) | = \emptyset$.
		
		\noindent{\scriptsize* We may replace the trivial estimate $0$ with any value in $f(\latsprow, \latspcol)$ to handle the exception $| \baseest(u,i) | = \emptyset$.}
   \end{algorithmic}

    \caption{User-user Fixed Radius Nearest Neighbor}\label{alg:user_algorithm}
\end{algorithm}

\begin{remark}
Note that Algorithm \ref{alg:user_algorithm} is equivalent to the classical user-user fixed radius nearest neighbor collaborative filtering algorithm. 
The algorithm analyzed in the preliminary version of this work, cf. \cite{lee2016blind}, is similar but not identical to Algorithm \ref{alg:user_algorithm}. 
The algorithm in \cite{lee2016blind} is motivated by the kernel regression of the first order and it is asymptotically equivalent to a mean-adjusted variant 
of the user-user $k$-nearest neighbor collaborative filtering algorithm. On the other hand, Algorithm \ref{alg:user_algorithm} does not adjust the estimate 
the empirical means and it is kernel regression of the zeroth order in that sense.
\end{remark}

\subsubsection*{Intuition behind the Algorithm}

The following equations provide  intuition for the reason why Algorithm \ref{alg:user_algorithm} works. 
Assuming the row latent variables are instantiated as $\featrow{1}, \ldots, \featrow{m}$ and $|\baseest(u,i)| = \nbase(\eta)$, the mean squared error of 
the resulting estimate is given as\footnote{In fact, $\baseest(u,i)$ is not deterministic but random. We provide 
a complete analysis with formal proofs in Section \ref{sec:main.result} and Appendix \ref{sec:proof_main_thm}. }
\begin{align}
	&\Expc{ \big( \hA(u,i) - A(u,i) \big)^2 } \nonumber\\
		&= \Expc{ \bigg( \frac{1}{ \nbase } \sum_{(v,i) \in \baseest(u,i)} \big[  A(v,i) - A(u,i) + N(v,i) \big] \bigg)^2  }\nonumber\\
		&\leq \max_{ (v,i) \in \baseest(u,i)} \big\| f(\featrow{u}, \cdot) - f(\featrow{v}, \cdot) \big\|_{L^2}^2 + \frac{2 \sigma^2}{\nbase(\eta)}.	\label{eqn:est_upper_bound}
\end{align}
Here, $\Expc{ ~\cdot~ } = \bbE\big[ ~\cdot ~\big| ~\featrow{1}, \ldots, \featrow{m}, |\baseest(u,i)|$$ = \nbase(\eta)\big]$ denotes the conditional expectation.
This expression shows that the expected squared error conditioned on row latent variables is directly related to the squared $L^2$ distance between 
the slices of the latent function $f$ associated to rows $u$ and $v$. The good news is that the $L^2$ distance can in fact be estimated from the data itself. 
For any pair of row indices $u, v \in [m]$,
\begin{align*}
	&\Exp{ \big( Z(u,j) - Z(v,j) \big)^2 ~\big|~ \featrow{u}, \featrow{v} }\\
		&\qquad= \Exp{ \big( A(u,j) - A(v,j) \big)^2 ~\big|~ \featrow{u}, \featrow{v} }\\
			&\qquad\qquad + \Exp{ \big( N(u,j) - N(v,j) \big)^2 ~\big|~ \featrow{u}, \featrow{v} }\\
		&\qquad= \big\| f(\featrow{u}, \cdot) - f(\featrow{v}, \cdot) \big\|_{L^2}^2 + 2 \sigma^2.
\end{align*}
Because of the concentration of measure, we expect the dissimilarity between rows $u$ and $v$ defined in \eqref{eqn:user_dissim} to be close to 
$ \big\| f(\featrow{u}, \cdot) - f(\featrow{v}, \cdot) \big\|_{L^2}^2 + 2 \sigma^2$ as the size of the overlap $ | \basesim(u) \cap \basesim(v) | $ increases. 

\begin{remark}
	For fixed $m$ rows and $n$ columns, choosing a large $\eta$ leads to the increase in the size of $\baseest(u,i)$, which results in the increase in the first term in \eqref{eqn:est_upper_bound} 
	and the decrease in the second term in \eqref{eqn:est_upper_bound}. This demonstrates a bias-variance tradeoff associated to the choice of algorithmic 
	parameter $\eta \geq 0$.
\end{remark}

\begin{remark}
	Our algorithm does not require any prior knowledge on the model parameters such as the regularity parameters $\lip$, diameter of ${\latsprow}$ and 
	the noise variance $\sigma$. The algorithmic parameter $\eta \geq 0$ can be chosen arbitrarily. However, we need $\eta$ to be in a certain range 
	to achieve good theoretical guarantee for the upper bound on its MSE. See Theorem \ref{thm:main_mse_user}.
\end{remark}


\subsection{Additional Examples of Our Algorithm: Other Variations}\label{sec:other_algs}
In this section, we exhibit other variations of our algorithm and provide intuition for them. However, we do not formally prove error bounds for 
these variations in this paper.

\subsubsection{Item-item and User-item Variants}
In Section \ref{sec:user_user}, we described a user-user fixed radius nearest neighbor algorithm that utilizes the (dis-)similarity between rows, $\dissimrow(u,v)$ for 
$u, v \in [m]$. We can apply the same idea to use the similarity between columns, $\dissimcol(i,j)$ for $i, j \in [n]$. Moreover, once the similarity between 
rows and columns are identified, we can define similarity between any pair of index tuples $(u,i), (v, j) \in [m] \times [n]$ based on the $\dissimrow(u,v)$ and 
$\dissimcol(i,j)$. For example, we may define $\dissim((u,i), (v,j)) = \max \{ \dissimrow(u,v), \dissimcol(i,j) \}$.

\subsubsection{Kernels for Weighting}
Once the (dis-)similarities between rows and columns are identified, we need a weighting scheme to define the estimates. For example, in Algorithm \ref{alg:user_algorithm}, we use the the hard-thresholding kernel 
in which the weight is $1$ if the row dissimilarity is less than $\eta$ and the weight is $0$ otherwise. However, this choice is not necessary 
and our algorithm can be implemented with any choice of kernels (i.e. weighting scheme).

\begin{example}[Gaussian Kernel Weights]
Given $(u,i) \in [m] \times [n]$, suppose that we have computed dissimilarity between $(u,i)$ and $(v,j)$ for all $(v,j) \in \baseest(u,i)$, e.g., by defining
$\dissim((u,i), (v,j)) = \max \{ \dissimrow(u,v), \dissimcol(i,j) \}$. Then, given a bandwith parameter $\lambda \in \RealsP$, we define the weights according 
to a Gaussian kernel\footnote{One may doubt why we call this Gaussian kernel instead of Laplacian kernel. The reason is that we treat $\dissim((u,i), (v,j))\big)$ 
as a proxy of the {\em squared} $L^2$ distance.} as: for each $(v,j) \in [m] \times [n]$,
\begin{equation*}
	w_{ui}(v, j) = \begin{cases}
				e^{-\lambda ~\dissim((u,i), (v,j)) } & \text{if }(v,j) \in \baseest(u,i),	\\
				0		&	\text{else},
		\end{cases}
\end{equation*}
and then define
\begin{equation*}
	\hA(u,i) = \frac{\sum_{(v,j)} w_{ui}(v,j) Z(v,j) }{ \sum_{(v,j)} w_{ui}(v,j)}.
\end{equation*}
When $\lambda \to \infty$, the estimate $\hA(u,i)$ only depends on the `nearest' neighbor $(v,j)$. On the other hand, when $\lambda = 0$, the algorithm 
equally averages all the $Z(v,j)$ for $(v,j) \in \baseest(u,i)$. Empirically, this variant of the algorithm seems to perform very well with an appropriate selection 
of the bandwidth parameter $\lambda$, which can be tuned using cross validation.
\end{example}

\subsubsection{Other Ways to Define Dissimilarity} We suggest other alternatives to define similarities.
\paragraph{Higher Order Information} 
Algorithm \ref{alg:user_algorithm} detects reliable neighbors $(v,j) \in \baseest(u,i)$ and estimates $\hA(u,i)$ by averaging $Z(v,j)$. This procedure is 
equivalent to traditional local smoothing, or kernel regression of the $0$-th order, except that we needed to estimate the `proximity' of the 
latent features from the data matrix. If we consider $Z(v,j)$ as a constant function, then $\hA(u,i)$ can be viewed as the evaluation 
of the average of the constant function $Z(v,j)$ at $(u,i)$.

The next question naturally arises: ``Can we come up with an algorithm that is equivalent to higher-order kernel regression algorithm?'' 
For example, one can build linear regression estimators centered at $(v,j)$ for each $(v,j) \in \baseest(u,i)$, instead of considering the constant function 
at the level of $Z(v,j)$. Evaluating the average of such linear estimators at $(u,i)$ will yield a different, probably more refined, estimate $\hA(u,i)$. 
This turns out to be equivalent to kernel regression of the $1$-st order with latent features. This version of algorithm is analyzed in the preliminary 
version of this paper, cf. \cite{lee2016blind}.

\paragraph{Beyond the Euclidean Latent Space}
In our proposed algorithm, we defined the dissimilarities between the rows and the columns using empirical $L^2$ distance of the slices of the latent function. 
Depending on the geometry of the latent space, other ways of defining dissimilarities can be favorable. For example, suppose that only the `direction' of the rows in the data matrix matters, 
while the amplitude does not carry much information. Then it would make sense to measure the dissimilarity between rows by estimating their `angles' instead 
of their Euclidean distance\footnote{To be precise, we want to define the dissimilarity between the equivalent classes of the rows. In this example, the equivalence 
relation is defined by positive scaling}. In fact, this is what cosine similarity measures, which is commonly used in classical 
collaborative filtering.

More importantly, the canonical latent space in the data generation process can be non-Euclidean. For example, the latent variables could be drawn from 
a sphere\footnote{more generally, from a Riemannian manifold with non-zero curvature}. Cosine similarity makes use of this knowledge, while the $L^2$ 
dissimilarity does not -- it essentially isometrically embeds the sphere to a higher-dimensional Euclidean space and then uses 
the metric of the ambient space. This could be problematic because one may need a much larger dimension for Euclidean embedding than the intrinsic dimension. 
Additionally in some cases it may be impossible to find an isometric embedding to the Euclidean space of any finite dimension; in particular, when the latent space is negatively curved. 
One prominent example of such latent space is the tree equipped with the geodesic distance; it is impossible to find an isometric embedding of an infinitely growing 
$3$-regular tree to $\Reals^d$ for any $d < \infty$.

Both the question of finding a good similarity/dissimilarity function that utilizes the structure of the latent space and the question of systematically exploiting 
higher-order information of the latent function remain important direction of future research.

\subsection{Computational Complexity}\label{sec:computational_complexity}

In order to compute $\dissimrow(u,v)$, we need to average differences over the set $| \basesim(u) \cap \basesim(v) |$. As the sparsity of each row is in expectation $np$, the cost of the similarity computation is in expectation bounded by $O(np)$. We need to compute $\dissimrow(u,v)$ for all $\binom{m}{2}$ pairs. For each row $u$, we need to compute the set $\{v: \dissimrow(u,v) \leq \eta\}$, which takes $O(m)$. To compute the final estimate for each entry $(u,i)$ we need to average over the set $\baseest(u,i)$, which in expectation is bounded above by $mp$. Therefore, the computational complexity is bounded by $O(m^2 (np) + m^2 + m n (mp)) = O(m^2 np)$.

The algorithm can be parallelized into two phases. In phase 1, we can in parallel compute $\dissimrow(u,v)$ for all $\binom{m}{2}$ pairs, each taking computation time $O(np)$. In phase 2, we can compute the estimates $\hat{A}(u,i)$ in parallel for all $mn$ indices of the matrix, each taking computation time $O(m + mp)$. Therefore assuming we have $O(m^2 + mn)$ processors that can compute simultaneously, the parallelized computational complexity is $O(np + m)$.
		 
\section{Main Results} \label{sec:main.result}


\subsection{Theorem Statement}

We present a theorem which upper bounds the MSE of the estimate produced by the user-user fixed radius nearest neighbor algorithm 
presented in Section \ref{sec:user_user} (cf. Algorithm \ref{alg:user_algorithm}). Recall from our problem statement in Section \ref{sec:setup} 
that $p$ is the probability each entry is observed, 
$\fbound$ is an upper bound on magnitude of the latent function,
and $\sigma^2$ is an upper bound on the noise variance\footnote{In fact, $\| N(u,i) \|_{\psi_2} 
\leq \sigma$. We may identify $\sigma^2$ with an upper bound on the variance, up to an absolute constant.}.

\medskip
\subsubsection{Informally Stated Asymptotic Upper Bounds} 

We present an informal version of our main theorem under a simplified setting for ease of
exposition. The general form of main result is stated in Theorem \ref{thm:main_mse_user} which is
rather involved. To that end, consider a simplified setting where there are only finite types of 
row latent variables. In that case, we obtain the following result.
\begin{theorem*}[Informal Statement of Corollary \ref{coro:finite}]
Suppose that the measure for the latent row space has finite support. 
If $p \gg \max\Big(\frac{1}{m}, \sqrt{\frac{\log m}{n}}\Big)$, then
\[	\MSE(\hat{A}) = O \Bigg( (\fbound + \sigma)^2 \bigg(\sqrt{\frac{\log (m-1)}{(n-1)p^2}} \vee \frac{|\supp(\murow)|}{(m-1)p} \bigg)\Bigg).	\]
\end{theorem*}
\noindent Note that the prefix constant $\fbound + \sigma$ is the measure of variability (complexity) for the latent function (up to noise).

There are two main sources of the estimation error: (i) the error in estimation of similarities (step 2 of Algorithm \ref{alg:user_algorithm}) 
and (ii) the error in approximation by smoothing (step 3 of Algorithm \ref{alg:user_algorithm}). Intuitively, empirical estimates of similarity/dissimilarity 
between two rows becomes more accurate as $(n-1)p^2$, which is equal to the expected size of the overlap, increases. This is captured by the term 
$\sqrt{\frac{\log (m-1)}{(n-1)p^2}}$. 

Given the similarities are sufficiently accurate, the estimation error is dominated by the approximation error. There are $|\supp(\murow)|$ types 
of rows and $(m-1)p$ number of available samples, hence, $| \baseest(u,i) | \approx \frac{(m-1)p}{ | \supp(\murow)|}$ in expectation. This captures the second term in the error. This also reflects the impact of local geometry of the measure $\murow$ on the estimation error.

\medskip
\subsubsection{Formal Statement of the Main Theorem} 
We define the 
function $\localprob$ for $x \in \latsprow$ and $r \geq 0$ as
\begin{align*}
\localprob(x, r) = \bP_{\featrow{v} \sim \murow}\left(\big\| f(x, \cdot) - f(\featrow{v}, \cdot) \big\|_{L^2}^2 \leq r\right).
\end{align*}
%
%
%
%
%
%
Theorem \ref{thm:main_mse_user} presents an upper bound on the mean-squared error (cf. \eqref{eqn:MSE}) of 
the user-user Fixed radius nearest neighbor algorithm described in Algorithm \ref{alg:user_algorithm}. 
Before presenting the theorem statement, we remark that 
\begin{align*}
	\MSE(\hA)&= \Exp{\big( \hA(u,i) - A(u,i) \big)^2 }\\
			&= \Exp{\big( \hA(1,1) - A(1,1) \big)^2 }.
\end{align*}
due to the exchangeability of the model and the linearity of the expectation. In Theorem \ref{thm:main_mse_user}, 
we provide an upper bound on $\Exp{\big( \hA(1,1) - A(1,1) \big)^2 ~\big|~ \featrow{1} }$, conditioned on $\featrow{1}$, 
which is the latent feature of the first row.

\begin{theorem}[Main Theorem]\label{thm:main_mse_user}
Let $\hA$ be the estimator returned by the Algorithm \ref{alg:user_algorithm}. Let constant $K \triangleq \Big(  \frac{2\fbound}{\sqrt{\ln 2}} + 2 \sigma \Big)^2$. Suppose that our algorithm uses threshold parameter $\eta \geq \eta' 
+ K \max\bigg( \sqrt{\frac{4 \log (m-1)}{c(n-1)p^2}}, \frac{4 \log (m-1)}{c(n-1)p^2} \bigg)$ for some $\eta' \geq 2 \sigma^2$. Then
\begin{align*}
	&\Exp{\big( \hA(1,1) - A(1,1) \big)^2 ~\big|~ \featrow{1} }\\
	&\leq \big( \eta - 2\sigma^2 \big) + CK \bigg[  \frac{1}{\sqrt{ (n-1) p^2}} + \exp \big(-c (n-1)p^2\big) \bigg]\\
	&	+ 2 \bbE_{ \featrow{1}} \Bigg[ \bigg[  (m-1)p  \cdot  
		\localprob\Big(\featrow{1}, \eta' - 2\sigma^2\Big) \bigg]^{-1}\Bigg] \\
	&	+ \big( \fbound^2 + 1 \big)\\
		&	\qquad \times \bbE_{ \featrow{1}}\Bigg[ \exp\bigg( - \frac{(m-1)p}{8} \cdot 
			\localprob\Big(\featrow{1}, \eta' - 2\sigma^2\Big) \bigg) \Bigg]\\
	&	+ \big[ (m-1) + \fbound^2 \big] \bigg[  \frac{1}{(m-1)^2} + \exp \Big( - \frac{(n-1)p^2}{8} \Big) \bigg].
\end{align*}
In the above expression, $C, c > 0$ are absolute constants, and $\localprob(x, r) = \bP_{\featrow{v} \sim \murow}\left(\big\| f(x, \cdot) - f(\featrow{v}, \cdot) \big\|_{L^2}^2 \leq r\right)$ for any $x \in  \latsprow$ and $r > 0$. 
\end{theorem}

The threshold $\eta \geq 0$ is a tunable parameter that our algorithm uses. Our analysis captures the bias-variance tradeoff associated with 
the parameter $\eta$. When $\eta$ is too large, our upper bound becomes loose because the algorithm utilizes information from 
less reliable neighbors; on the other hand, when $\eta$ is chosen too small, the resulting estimate count only on a few neighbors 
and suffers from large variance. 

The first two terms of the MSE bound come from bounding the bias of the estimator. The first term $\big( \eta - 2\sigma^2 \big)$ reflects the bias that is unavoidable due to the selection of the threshold $\eta$ and the error in estimating $\dissimrow(u,v)$, and the second term bounds the tail of the bias along with the bad event that $| \basesim(u) \cap \basesim(v) |$ is too small to produce a good estimate of $\dissimrow(u,v)$. The remaining three terms bound the variance of the estimator; the error bound for the trivial estimates (when $\baseest(1,1) = \emptyset$), given in the form of $\fbound^2 \Prob{\baseest(1,1) = \emptyset} $ is also subsumed in these three terms.

Taking expectation of the bound in Theorem \ref{thm:main_mse_user} with respect to $\featrow{1} \sim \murow$, we obtain the desired 
bound on MSE, i.e.
\[	\MSE(\hA) = \bbE_{\featrow{1} \sim \murow} \Big[ \Exp{\big( \hA(1,1) - A(1,1) \big)^2 ~\big|~ \featrow{1} } \Big].	\]

\subsection{Implications} 

Theorem \ref{thm:main_mse_user} provides an upper bound on MSE for Algorithm \ref{alg:user_algorithm}, which depends on 
the underlying probability measure $\murow$ and the latent function $f$ through $\localprob$. Here we evaluate this implicit bound for
three special examples:
\begin{itemize}
\item The latent space has finitely many elements, or equivalently $\murow$ has finite support.
\item The latent space is unit hypercube in a finite dimensional space, latent function is Lipschitz. 
\item The latent space is a complete, separate metric space aka Polish space with bounded diameter, latent function is Lipschitz.
\end{itemize}  

\medskip 
\subsubsection{Finite Types} 
We state the following Corollary of Theorem \ref{thm:main_mse_user} when $\murow$ 
has finite support. 
\begin{coro}[Finite support]\label{coro:finite}
Suppose that $\latsprow$ is equipped with the discrete metric\footnote{$\drow(x_1, x_2) = 1$ if and only if $x_1 \neq x_2$.} topology 
and $\murow$ has finite support in $\latsprow$ with $\supp(\murow)$ denoting the support of $\murow$. Let 
$K = \Big(  \frac{2\fbound}{\sqrt{\ln 2}} + 2 \sigma \Big)^2$, $\eta' = 2\sigma^2$, and 
\[
\eta = \eta' + K \max\bigg( \sqrt{\frac{4 \log (m-1)}{c(n-1)p^2}}, ~\frac{4 \log (m-1)}{c(n-1)p^2} \bigg).
\]
If $p \geq \max\bigg(\frac{8}{m-1}, ~\Big(4 \vee \sqrt{\frac{2}{c}} \Big)\cdot \sqrt{\frac{\log(m-1)}{n-1}}\bigg)$, then 
\begin{align}\label{cor.MSE_finite}
	\MSE(\hat{A}) 
		&\leq C K \Bigg[\max\bigg( \sqrt{\frac{4 \log (m-1)}{c(n-1)p^2}}, ~\frac{4 \log (m-1)}{c(n-1)p^2} \bigg) \nonumber\\
		&	\qquad+ \frac{|\supp(\murow)|}{(m-1)p} \Bigg]
\end{align}
where $C, c > 0$ are absolute constants. 
\end{coro}
That is, as long as $p \gg \max\Big(\frac{1}{m}, \sqrt{\frac{\log m}{n}}\Big)$, under discrete measure with finite support 
\[	\MSE(\hat{A}) = O \Bigg( (\fbound + \sigma)^2 \bigg(\sqrt{\frac{\log (m-1)}{(n-1)p^2}} \vee \frac{|\supp(\murow)|}{(m-1)p} \bigg)\Bigg).	\]

\begin{proof}[Proof of Corollary \ref{coro:finite}]
Our interest is in bounding 
\begin{align}\label{eq.term1}
\bbE_{ \featrow{1}} \Bigg[ \bigg[  (m-1)p  \cdot  \localprob\Big(\featrow{1}, \eta' - 2\sigma^2\Big) \bigg]^{-1}\Bigg] 
\end{align}
and 
\begin{align}\label{eq.term2}
\bbE_{ \featrow{1}}\Bigg[ \exp\bigg( - \frac{(m-1)p}{8} \cdot \localprob\Big(\featrow{1}, \eta' - 2\sigma^2\Big) \bigg) \Bigg]
\end{align}
from the statement of Theorem \ref{thm:main_mse_user} to obtain the desired result. Since we are considering discrete metric with measure having finite support, 
for any $r \geq 0$ and $x \in \supp(\murow)$,  $\localprob(x, r) \geq \murow(x)$. Given choice of $\eta'$, $\eta' - 2\sigma^2 \geq 0$. Therefore, \eqref{eq.term1} can be written as
\begin{align}\label{eq.term1a}
&\bbE_{ \featrow{1}} \Bigg[ \bigg[  (m-1)p  \cdot  \localprob\Big(\featrow{1}, \eta' - 2\sigma^2\Big) \bigg]^{-1}\Bigg] \nonumber\\
& \qquad \leq \sum_{x \in \supp(\murow)} \murow(x) \times \frac{1}{(m-1)p  \murow(x)} \nonumber  \\
& \qquad= \frac{| \supp(\murow)| }{(m-1)p}.
\end{align}
Similarly, \eqref{eq.term2} reduces to 
\begin{align}\label{eq.term2a}
&\bbE_{ \featrow{1}}\Bigg[ \exp\bigg( - \frac{(m-1)p}{8} \cdot \localprob\Big(\featrow{1}, \eta' - 2\sigma^2\Big) \bigg) \Bigg] \nonumber\\
& \qquad\leq \sum_{x \in \supp(\murow)} \murow(x) \times \exp\bigg( - \frac{(m-1)p}{8} \murow(x)  \bigg) \nonumber \\
& \qquad\leq | \supp(\murow)| ~\Bigg\{ \sup_{\theta \in [0,1]} \theta \exp\bigg( - \frac{(m-1)p}{8} \theta  \bigg) \Bigg\} \nonumber \\
& \qquad= | \supp(\murow)| ~\frac{8\exp(-1)}{(m-1)p}.
\end{align}
In above, we have used an easy to verify fact that $ \sup_{\theta \in [0,1]}  \theta \exp\bigg( - \frac{(m-1)p}{8} \theta  \bigg)$ is
achieved for $\theta = \frac{8}{(m-1)p} \in [0,1]$. Replacing \eqref{eq.term1a} and \eqref{eq.term2a} in the statement of Theorem \ref{thm:main_mse_user}, and realizing that
other terms in the bound of  Theorem \ref{thm:main_mse_user} are asymptotically smaller order, we obtain \eqref{cor.MSE_finite}.
\end{proof}

\medskip 
\subsubsection{Uniform Measure over $[0,1]^d$, Lipschitz Latent Function}
Let $\lip$ denote the Lipschitz constant of the latent function $f$. We consider the setting where $\murow$ is uniform Lebesgue measure on $[0,1]^d$, the unit cube in $d$ dimension.
We define
$B(x,r) = \{x' \in \latsprow: \drow(x, x') \leq r \}$ for any $r > 0$. 
By Lipschitzness, $\drow(\featrow{1}, \featrow{v}) \leq \frac{\sqrt{\eta' - 2\sigma^2}}{L}$ implies $\big\| f(\featrow{1}, \cdot) - f(\featrow{v}, \cdot) \big\|_{L^2}^2 \leq \eta' - 2\sigma^2$. Therefore, 
\[\localprob(\featrow{1}, \eta'-2\sigma^2) \geq 
\murow\Bigg( B\Big(\featrow{1}, \frac{\sqrt{\eta' - 2\sigma^2}}{L} \Big)\Bigg).\]
There exists universal constants $\alpha, \beta > 0$ such that for any $d \geq 1$, and $x \in [0,1]^d$ and $r > 0$, 
\begin{align*}
\mbox{Vol}(B(x,r)) & \geq \min(1, ~\alpha \beta^d  r^d). 
\end{align*} 
We shall assume that $\latspcol = [0,1]^d$ as well. Therefore, 
$$ \sup_{\alpha, \beta, \alpha', \beta' \in [0,1]^d} | f(\alpha, \beta) - f(\alpha', \beta')| \leq L \times \sqrt{d}, $$
and hence $\fbound \leq L \sqrt{d}$ since there exists $\alpha \in \latsprow, \beta \in \latspcol$ such that
$f(\alpha, \beta) = 0$.  Now, we state the following implication of Theorem \ref{thm:main_mse_user} in this setting. 
\begin{coro}[Unit cube and the Lebesgue measure]\label{coro:uniform}
Let $\murow$ be the uniform measure on $\latsprow = [0,1]^d$. Let $K = \big(  \frac{2L \sqrt{d}}{\sqrt{\ln 2}} + 2 \sigma \big)^2$, 
$ \eta' = 2 \sigma^2 +  \alpha^{2/d} \beta^2 L^2 \big((m-1)p\big)^{-\frac{2}{d+2}}$, and 
\[
	\eta = \eta' +  K \max\bigg( \sqrt{\frac{4 \log (m-1)}{c(n-1)p^2}}, ~\frac{4 \log (m-1)}{c(n-1)p^2} \bigg).
\]
If $p \geq \max\bigg(\frac{8}{m-1}, ~\Big(4 \vee \sqrt{\frac{2}{c}} \Big)\cdot \sqrt{\frac{\log(m-1)}{n-1}}\bigg)$, then 
\begin{align}\label{cor.MSE_uniform}
	\MSE(\hat{A}) 
		&\leq C K \Bigg[ \max\bigg( \sqrt{\frac{4 \log (m-1)}{c(n-1)p^2}}, ~\frac{4 \log (m-1)}{c(n-1)p^2} \bigg)	\nonumber\\
		& \qquad+  \alpha^{2/d} \beta^2 L^2 \big[ (m-1)p \big]^{-\frac{2}{d+2}}	\nonumber\\
		& \qquad+  \exp\Big( - \frac{1}{8} \big[ (m-1)p \big]^{\frac{2}{d+2}} \Big) \Bigg],
\end{align}
where $C, C', c > 0$ are absolute constants.
\end{coro}
That is, as long as $p \gg \max\Big(\frac{1}{m}, \sqrt{\frac{\log m}{n}}\Big)$, under the uniform measure, 
\begin{align*}
	&\MSE(\hat{A})\\
	& = O \Bigg( (L\sqrt{d} + \sigma)^2 \bigg( \sqrt{\frac{\log (m-1)}{(n-1)p^2}} \vee L^2 \Big(\frac{1}{(m-1)p}\Big)^{\frac{2}{d+2}} \bigg) \Bigg).
\end{align*}
\begin{proof}[Proof of Corollary \ref{coro:uniform}]
Similar to proof of Corollary \ref{coro:finite}, our interest is in bounding \eqref{eq.term1} and \eqref{eq.term2}.
To that end, under
uniform Lebesgue measure $\murow$, on $d$ dimensional unit cube $\latsprow$, it follows that for any $x \in \latsprow$ and $r > 0$
\begin{align*}
\localprob(\featrow{1}, \eta'-2\sigma^2) 
&\geq  \murow\bigg( B(\featrow{1}, \frac{\sqrt{\eta' - 2\sigma^2}}{L} \bigg) \\
& \geq \min\bigg(1, \alpha \beta^d \Big(\frac{\sqrt{\eta' - 2\sigma^2}}{L}\Big)^d\bigg).
\end{align*}
Therefore, by choice of $\eta'$
\begin{align}\label{eq.term1b}
&\bbE_{ \featrow{1}} \Bigg[ \bigg[  (m-1)p  \cdot  \localprob\Big(\featrow{1}, \eta' - 2\sigma^2\Big) \bigg]^{-1}\Bigg] \nonumber\\
& \qquad \leq \frac{1}{(m-1)p} \times \frac{\alpha \beta^d L^d}{(\eta' - 2\sigma^2)^{d/2}} \nonumber \\
& \qquad = \Big((m-1)p\Big)^{-\frac{2}{d+2}}.
\end{align}
Similarly, 
\begin{align}\label{eq.term2b}
&\bbE_{ \featrow{1}}\Bigg[ \exp\bigg( - \frac{(m-1)p}{8} \cdot \localprob\Big(\featrow{1}, \eta' - 2\sigma^2\Big) \bigg) \Bigg] \nonumber\\
& \qquad \leq \exp\bigg( - \frac{(m-1)p}{8}  \times \frac{(\eta' - 2\sigma^2)^{d/2}}{\alpha \beta^d L^d}\bigg) \nonumber \\
& \qquad \leq \exp\bigg( - \frac18 \Big((m-1)p\Big)^{\frac{2}{d+2}}\bigg).
\end{align}
From \eqref{eq.term1b} and \eqref{eq.term2b}, and statement of Theorem \ref{thm:main_mse_user}, 
and realizing that
other terms in the bound of  Theorem \ref{thm:main_mse_user} are asymptotically smaller order, 
we conclude \eqref{cor.MSE_uniform}.
\end{proof}

\medskip 
\subsubsection{Bounded Polish Space, Lipschitz Latent Function} 
Let $\lip$ denote the Lipschitz constant of the latent function $f$. 
We define
$B(x,r) = \{x' \in \latsprow: \drow(x, x') \leq r \}$ for any $r > 0$. 
By Lipschitzness, $\drow(\featrow{1}, \featrow{v}) \leq \frac{\sqrt{\eta' - 2\sigma^2}}{L}$ implies $\big\| f(\featrow{1}, \cdot) - f(\featrow{v}, \cdot) \big\|_{L^2}^2 \leq \eta' - 2\sigma^2$. Therefore, 
\begin{align}\label{eq.c3.0}
\localprob(\featrow{1}, \eta'-2\sigma^2) & \geq \murow\bigg( B(\featrow{1}, \frac{\sqrt{\eta' - 2\sigma^2}}{L} \bigg).
\end{align}
We consider the row latent space $\latsprow$ to be a complete, separable metric space, i.e. a Polish space\footnote{Recall that a metric space is separable if it has a countable dense subset, and 
complete if every Cauchy sequence converges within the space.}. Let the diameter of the space be bounded, i.e. there exists finite $D > 0$ such 
$\sup_{x, x' \in \latsprow} \drow(x, x') \leq D$. We shall also assume that the diameter of $\latspcol$ is bounded by $D$ as well. This implies 
\begin{align*}
\fbound & \leq L D, 
\end{align*}
since there exists $\alpha \in \latsprow, \beta \in \latspcol$ such that $f(\alpha, \beta) = 0$. 

An important consequence of $\latsprow$ being Polish is that $\murow$ is {\em tight},
i.e. for any $\delta > 0$, there exists a compact set $S_\delta \subseteq \latsprow$ such that $\murow(S_\delta) \geq 1-\delta$. 
Due to compactness of $S_\delta$ and the space being Polish, there exists a finite number of balls of any given radius $\varepsilon > 0$ of choice
such that they cover $S_\delta$. That is, for any $\varepsilon, \delta > 0$, the effective covering number $N_{\text{eff}}(\latsprow, \varepsilon, \delta)$ is always finite. 
Let $B(x_i, \varepsilon)$ for $i \in [N_{\text{eff}}(\latsprow, \varepsilon, \delta)]$ denote the collection of balls so that
$S_\delta \subseteq \cup_{i \in [N_{\text{eff}}(\latsprow, \varepsilon, \delta)]} B(x_i, \varepsilon)$. By construction,
$\murow\Big(\cup_{i \in [N_{\text{eff}}(\latsprow, \varepsilon, \delta)]} B(x_i, \varepsilon)\Big) \geq \murow(S_\delta) \geq 1-\delta$. 
Define $\bad = \{ i \in [N_{\text{eff}}(\latsprow, \varepsilon, \delta)]~:~ \murow(B(x_i, \varepsilon)) < \delta/N_{\text{eff}}(\latsprow, \varepsilon, \delta)\}$.
Let $\good = \{ i \in [N_{\text{eff}}(\latsprow, \varepsilon, \delta)]~:~ i \notin \bad\}$. Therefore, 
\begin{align*}
\murow\Big(\cup_{i \in \bad} B(x_i, \varepsilon)\Big) & < |\bad| \frac{\delta}{N_{\text{eff}}(\latsprow, \varepsilon, \delta)} \nonumber \\
	& \leq~ \delta, 
\end{align*}
since $|\bad| \leq N_{\text{eff}}(\latsprow, \varepsilon, \delta)$. Therefore, 
\begin{align*}
\murow\Big(\cup_{i \in \good} B(x_i, \varepsilon)\Big) & \geq \murow(S_\delta) - \murow\Big(\cup_{i \in \bad} B(x_i, \varepsilon)\Big)	\nonumber \\
	&\geq~ 1-2\delta.
\end{align*}
For any $x \in \cup_{i \in \good} B(x_i, \varepsilon)$, there exists $i \in \good$ such that $x \in B(x_i, \varepsilon)$. Hence, 
$B(x_i, \varepsilon) \subseteq B(x, 2\varepsilon)$. Therefore, for any $x \in S^\good(\delta, \varepsilon) \equiv \cup_{i \in \good} B(x_i, \varepsilon)$,
\begin{align*}
\murow\Big(B(x, 2\varepsilon)\Big) & \geq \frac{\delta}{N_{\text{eff}}(\latsprow, \varepsilon, \delta)}
\end{align*}
where $\murow\Big(S^\good(\delta, \varepsilon)\Big) ~\geq~ 1-2\delta$. 
By \eqref{eq.c3.0} it follows that for $\featrow{1} \in S^\good(\delta, \varepsilon)$ with $\murow\Big(S^\good(\delta, \varepsilon)\Big) \geq 1-2\delta$,
\begin{align}\label{eq.c3.4}
\localprob(\featrow{1}, 4 \varepsilon^2 L^2) & \geq  \frac{\delta}{N_{\text{eff}}(\latsprow, \varepsilon, \delta)}
\end{align}
for any choice of $\varepsilon > 0,~\delta \in (0,1/2)$. 
Let us choose $\varepsilon = \varepsilon(\delta)$ where  $4 \varepsilon^2(\delta) L^2 = \delta$. Replacing this choice of $\varepsilon$ in \eqref{eq.c3.4}, we obtain 
that for any $\delta > 0$, there exists a set $S^\prime(\delta) \equiv S^\good(\delta, \varepsilon(\delta))$ such that (i) $\murow\big(S^\prime(\delta)\big) \geq 1-2\delta$ and (ii) for all $\featrow{1} \in S^\prime(\delta)$,
\begin{align}\label{eq.c3.5}
	 \localprob(\featrow{1}, \delta) & \geq  \frac{\delta}{N_{\text{eff}}(\latsprow, \frac{\sqrt{\delta}}{2L}, \delta)}.
\end{align}
For all $t \geq N_{\text{eff}}(\latsprow, \frac{1}{2L}, 1)$, we define function $\delta^{\star}(t)$ to be
\begin{align*}
\delta^\star(t) & = \inf \Big\{ \delta : \delta^2 \geq \frac{N_{\text{eff}}(\latsprow, \frac{\sqrt{\delta}}{2L}, \delta)}{t} \Big\}.
\end{align*}
For any $\delta > 0$, $N_{\text{eff}}(\latsprow, \frac{\sqrt{\delta}}{2L}, \delta)$ is finite, and thus for 
$t \geq \delta^{-2}N_{\text{eff}}(\latsprow, \frac{\sqrt{\delta}}{2L}, \delta)$, $\delta^\star(t) \leq \delta$. We can verify that $\delta^\star(t)$ is also
monotonically non-increasing. Therefore, it follows that 
\begin{align*}
\lim_{t\to\infty} \delta^\star(t) & = 0.
\end{align*}
Using the above developed machinery, now we are ready to bound $\MSE(\hat{A})$ as summarized in Corollary \ref{coro:polish}.
\begin{coro}[Bounded Polish Space]\label{coro:polish}
Let $\murow$ be any measure on $\latsprow$, which is assumed to be a bounded diameter Polish space. Let the diameter of 
$\latsprow$ and $\latspcol$ be bounded above by $D$. Let $K = \big(  \frac{2L \sqrt{d}}{\sqrt{\ln 2}} + 2 \sigma \big)^2$, 
$\eta' = 2 \sigma^2 +  \delta^\star((m-1)p)$, and 
\[
	\eta = \eta' +  K \max\Big( \sqrt{\frac{4 \log (m-1)}{c(n-1)p^2}}, ~\frac{4 \log (m-1)}{c(n-1)p^2} \Big).
\]
If $p \geq \max\bigg(\frac{8  N_{\text{\normalfont eff}}(\latsprow, \frac{1}{2L}, 1)}{m-1}, ~ \Big(4 \vee \sqrt{\frac{2}{c}} \Big)\cdot \sqrt{\frac{\log(m-1)}{n-1}}\bigg)$, then
\begin{align}\label{cor.MSE_polish}
	\MSE(\hat{A}) 
		&\leq C K \Bigg[ \max\bigg( \sqrt{\frac{4 \log (m-1)}{c(n-1)p^2}}, ~\frac{4 \log (m-1)}{c(n-1)p^2} \bigg) \nonumber\\
		&	\qquad +  \delta^\star((m-1)p) \Bigg],
\end{align}
where $C, C', c > 0$ are absolute constants, and $\delta^\star(t) \to 0$ as $t \to \infty$.
\end{coro}
As an immediate consequence of Corollary \ref{coro:polish}, it follows that as long as $p \gg \max\Big(\frac{1}{m}, \sqrt{\frac{\log m}{n}}\Big)$, 
for any measure on a bounded Polish space, 
\[	\MSE(\hat{A}) \to 0~\quad\text{as}\quad ~m,n \to \infty.\]
This implies that our estimator is consistent for any latent variable model with bounded Polish space as long as $p \geq \max(m^{-1+\delta}, n^{-\frac12 + \delta})$ for
any $\delta > 0$ and $\log m = o(n^\delta)$.  

If $\latsprow = [0,1]^d$, for any $\delta \in (0,1)$, $N_{\text{eff}}(\latsprow, \varepsilon, \delta) = O(\varepsilon^{-d})$. 
Therefore, $\delta^*(t) = O(t^{-\frac{2}{d+4}} L^{\frac{2d}{d+4}})$. We can conclude that 
as long as $p \gg \max\Big(\frac{1}{m}, \sqrt{\frac{\log m}{n}}\Big)$, and $\latsprow = [0,1]^d$, 
\begin{align}\label{eq.mse.uniform.polish}
	\MSE(\hat{A}) &\leq C K \Bigg[ \max\bigg( \sqrt{\frac{4 \log (m-1)}{c(n-1)p^2}}, ~\frac{4 \log (m-1)}{c(n-1)p^2} \bigg) \nonumber\\
		&	\qquad+  L^2 \big[ (m-1)p \big]^{-\frac{2}{d+4}} \Bigg],
\end{align}
with universal constant $C > 0$. It is worth noting the similarity between \eqref{cor.MSE_uniform} and \eqref{eq.mse.uniform.polish} -- the only
difference is the $d+4$ instead of $d+2$ in the denominator of the exponent for term $(m-1)p$. This is precisely the minimal cost of generalizing from the specific uniform
distribution to any arbitrary distribution. 
\begin{proof}
Since $p \geq \frac{8  N_{\text{eff}}(\latsprow, \frac{1}{2L}, 1)}{m-1}$, as argued before, $\delta^\star((m-1)p)$ is well defined and we shall
choose $\delta = \delta^\star((m-1)p)$. From the choice of $\eta'$, for any $\featrow{1} \in S^\prime(\delta)$, it follows from \eqref{eq.c3.5} that 
$\localprob(\featrow{1}, \delta)  \geq  \frac{\delta}{N_{\text{eff}}(\latsprow, \frac{\sqrt{\delta}}{2L}, \delta)}$. Therefore, for any 
$\featrow{1} \in S^\prime(\delta)$,
\begin{align}\label{eq.c3.1.0}
&\bbE_{ \featrow{1}} \Bigg[ \bigg[  (m-1)p  \cdot  \localprob\Big(\featrow{1}, \eta' - 2\sigma^2\Big) \bigg]^{-1}\Bigg] \nonumber\\
& \qquad \leq \frac{1}{(m-1)p} \times \frac{N_{\text{eff}}(\latsprow, \frac{\sqrt{\delta}}{2L}, \delta)}{\delta} \nonumber \\
& \qquad \leq \delta, 
\end{align}
where the last inequality follows from the fact that $\delta = \delta^\star((m-1)p)$.  Similarly, 
\begin{align}\label{eq.c3.1.1}
& \bbE_{ \featrow{1}}\Bigg[ \exp\bigg( - \frac{(m-1)p}{8} \cdot \localprob\Big(\featrow{1}, \eta' - 2\sigma^2\Big) \bigg) \Bigg] \nonumber\\
& \qquad \leq \exp\bigg( - \frac{(m-1)p}{8}  \times \frac{\delta}{N_{\text{eff}}(\latsprow, \frac{\sqrt{\delta}}{2L}, \delta)} \bigg) \nonumber \\
& \qquad \leq \exp\bigg( - \frac{1}{8 \delta} \bigg),
\end{align}
where again the last inequality follows from the fact that $\delta = \delta^\star((m-1)p)$. Replacing 
\eqref{eq.c3.1.0} and \eqref{eq.c3.1.1} in the statement of Theorem \ref{thm:main_mse_user}, and realizing that the
other terms in the bound of  Theorem \ref{thm:main_mse_user} are asymptotically smaller order and 
$\exp\big( - \frac{1}{8 \delta} \big) \leq \delta$ for any $\delta \in (0,1)$, we obtain
\begin{align}\label{eq.c3.1.2}
&\Exp{\big( \hA(1,1) - A(1,1) \big)^2 ~\Big|~ \featrow{1}  \in S'(\delta^\star((m-1)p))} \nonumber\\
& \qquad \leq C_1 K 
	\Bigg[ \max\bigg( \sqrt{\frac{4 \log (m-1)}{c(n-1)p^2}}, ~\frac{4 \log (m-1)}{c(n-1)p^2} \bigg)	\nonumber\\
& \qquad\qquad + \delta^\star((m-1)p)\Bigg]
\end{align}
By Cauchy-Schwarz inequality and the fact that $\murow(S'(\delta^\star((m-1)p))) \geq 1-2\delta^\star((m-1)p))$, 
we have
\begin{align}\label{eq.c3.1.3}
&\Exp{\big( \hA(1,1) - A(1,1) \big)^2 \Ind \big( { \featrow{1}  \notin S'(\delta^\star((m-1)p)) } \big)} \nonumber\\
& \leq \sqrt{\Exp{\big( \hA(1,1) - A(1,1) \big)^4} } \Prob{\featrow{1}  \notin S'(\delta^\star((m-1)p))} \nonumber \\
& \leq 2\delta^\star((m-1)p) \sqrt{\Exp{\big( \hA(1,1) - A(1,1) \big)^4}}. 
\end{align}
In what follows, we shall argue that 
\begin{align}\label{eq.c3.1.4}
\Exp{\big( \hA(1,1) - A(1,1) \big)^4} & = O\Big( \fbound^4 + \sigma^4 \Big). 
\end{align} 
By \eqref{eq.c3.1.2}, \eqref{eq.c3.1.3} and \eqref{eq.c3.1.4}, the main claim \eqref{cor.MSE_polish} follows. With that in mind, we shall now
establish \eqref{eq.c3.1.4} to conclude the proof. 
Recall from Algorithm \ref{alg:user_algorithm} that
\[	\hA(1,1) = 	 \frac{1}{| \baseest(1,1) |} \sum_{(v,1) \in \baseest(1,1)} Z(v,1) \] 
in our user-user fixed radius nearest neighbor algorithm when $\Nbase \geq 1$ and $\hA(1,1) = 0$ when $\Nbase = 0$. Therefore,
\begin{align*}
	&\Exp{\big( \hA(1,1) - A(1,1) \big)^4}\\
		&\qquad = \Exp{(A(1,1))^4 ~\Ind\big( { \Nbase = 0 }\big)} \\
		&\qquad\qquad+ \Exp{\big( \hA(1,1) - A(1,1) \big)^4 ~\Ind\big( { \Nbase \geq 1 }\big)}.
\end{align*}
By the assumption that the magnitude of the latent function $f$ is bounded by $\fbound$,
\begin{align*}
	\Exp{(A(1,1))^4 ~\Ind\big( { \Nbase = 0 }\big)  }
		&\leq \fbound^4.
\end{align*}
By introducing the notations 
\begin{align*}
	X &= \frac{\Ind\big( { \Nbase \geq 1 }\big)}{| \baseest(1,1) |} \sum_{(v,1) \in \baseest(1,1)} \big( A(v,1) - A(1,1) \big), \\
	Y &=  \frac{\Ind\big( { \Nbase \geq 1 }\big)}{| \baseest(1,1) |} \sum_{(v,1) \in \baseest(1,1)} N(v,1),
\end{align*}
we can rewrite the expected squared error as 
\begin{align}\label{eq.c3.1.5}
	&\Exp{\big( \hA(1,1) - A(1,1) \big)^4 ~\Ind\big( { \Nbase \geq 1 }\big)}	\nonumber\\
	&\qquad = \bbE\Big[ \big(X + Y\big)^4 \Big] ~=~ \sum_{k=0}^4 {4 \choose k} \bbE[X^k] \bbE[Y^{4-k}]. 
\end{align}
where we have used the independence of terms in $X$ and $Y$. It immediately follows that $X$ is a bounded random
variable with $|X| \leq 2\fbound$. Therefore, 
\begin{align}\label{eq.c3.1.6}
\bbE[X^k] & \leq \fbound^k, ~~\text{for}~~k \geq 1.
\end{align} 
The randomness influencing the selection of $\baseest(1,1)$ is independent of the random variables in the summation of term $Y$. Conditioned 
on $|\baseest(1,1)| = \ell$ for any $\ell \geq 1$, $Y$ is simply a summation of $\ell$ i.i.d. random variables, each
with ${\psi_2}$ norm bounded by $\sigma$ and zero mean. Therefore, it follows that conditioned on $|\baseest(1,1)| = \ell$
for any $\ell \geq 1$, 
\begin{align}\label{eq.c3.1.7}
\bbE{Y^k} & \leq \begin{cases} 0 & k = 1 \\
\sigma^2 & k = 2 \\
\sqrt{6} \sigma^3 & k = 3 \\
2 \sigma^4 & k = 4
\end{cases}
\end{align}
From \eqref{eq.c3.1.5}-\eqref{eq.c3.1.7}, the inequality \eqref{eq.c3.1.4} follows. This completes the proof of Corollary \ref{cor.MSE_polish}.
\end{proof}

\section{Discussion}\label{sec:discussion}

\subsection{Intuition}
\subsubsection{Structural Assumptions and Sample Complexity}
Without any structure assumed on the latent function $f$ and the latent spaces $\latsprow, \latspcol$, one would need $mn$ number of 
samples to recover the matrix. Our framework relies on two key assumptions to reduce the sample complexity: (i) (most of) the latent spaces 
can be covered by balls centered at a few representative points; and (ii) the latent function $f$ is regular (Lipschitz) and hence the proximity 
of two points in the latent space results in the similarity of the function values, which is observable from data.

Given a latent space $(\latsprow, \drow, \murow)$, we want to cover the entire space $\latsprow$ minus a small fraction that has negligible measure 
$\eps > 0$, with a collection of balls of radius $\delta > 0$. In other words, our model considers the metric entropy\footnote{The smallest number of bits 
that suffices to specify every point $x$ in the set with accuracy $\delta$ in the metric $\drow$.} of $\latsprow$ minus $\eps$-fraction 
as a measure of complexity. However, we don't have direct access to the latent features, but the distance between two points must be estimated from 
the pattern of the associated function values. If the latent function were an isometry, the distance in the image would faithfully reflect the distance in the domain. 
Lipschitzness assumption on $f$ is a robust analogue of the isometry assumption; if the latent function $f$ is $\lip$-Lipschitz, then the distance in the 
domain (= latent space) cannot be inflated more than $L$ times in the image; therefore, the image of the latent function for close neighbors in the latent space 
will behave in a similar fashion. With the Lipschitz assumption, the cost of indirect measurement is no more than $\lip$, and it suffices to consider a 
$\frac{\delta}{\lip}$-covering of the latent space instead of a $\delta$-covering.

\subsubsection{MSE Upper Bounds and the Rate of Convergence}
The minimax optimal rate for nonparametric regression is $O(N^{-2/(2+d)})$ where $N$ is the number of observations uniformly distributed in a $d$-dimensional hypercube. In our algorithm we use the data across columns to learn similarities between the rows, yet when we actually compute the final estimates, we only use the datapoints in each column separately. The number of observations per column is $mp$ in expectation. In the case when row latent features are sampled from a uniform measure over a $d$-dimensional cube, the second term of our MSE bound from Corollary \ref{coro:uniform} is $O((mp)^{-2/(2+d)})$. This indicates that the second term of our MSE bound is optimal for any estimator that uses entries in each column separately. In particular, even if the algorithm were given oracle knowledge of the row latent features it could not do better as long as it constrains itself to estimating each column separately. An algorithm which would average both similar rows and columns could be able to improve beyond this bound. The first term of our MSE bound comes from the step of estimating $\dissimrow(u,v)$ for all pairs of rows; however this first term is quickly dominated by the second term as $d$ increases, suggesting that the second term dominates the MSE. This suggests that our analysis is tight for our estimator. 

Our result can be compared with the upper bound on the MSE for the UVST estimator as presented in Theorem 2.7 of \cite{Chatterjee15}. For simplicity, consider the setting of a square matrix, i.e. $m = n$. For a matrix sampled from the latent variable model with latent variable dimension $d$, their theorem guarantees that
        \begin{equation}\label{eqn:USVT_bound}
	   \MSE(\hA^{\USVT}) \leq C \frac{m^{-\frac{1}{d+2}}}{\sqrt{p}}	
	\end{equation}
for some constant $C$ as long as $p \geq m^{-1 + \delta}$. This upper bound is meaningful only when $p > m^{-\frac{2}{d+2}}$, because the MSE bound in \eqref{eqn:USVT_bound} is bounded below by $C$ when $p \leq m^{-\frac{2}{d+2}}$. However, requiring $p > m^{-\frac{2}{d+2}}$ can be too restrictive when the latent dimension $d$ is large since it means that we need to sample almost every entry to achieve a nontrivial bound. Chatterjee's result stems from showing that a Lipschitz function can be approximated by a piecewise constant function, which upper bounds the rank of the (approximate) target matrix. This global discretization results in a large penalty with regards to the dimension of the latent space.  


In recent work, \cite{Xu2017} extends the analysis of USVT for graphon estimation, when the observation matrix is binary. He shows if the latent function is $\alpha$-H\"older smooth, the spectrum decays polynomially, and thus the MSE of the USVT estimator is bounded above by 
\[	\MSE(\hA^{\USVT}) = O((np)^{-\frac{2\alpha}{2\alpha+d}}).\]
This refined bound shows that at least for the binary observation case, the USVT estimator is consistent as long as $p = \omega(n^{-1})$, independent from the latent dimension $d$, relying on the quick decay of the spectrum.

Our algorithm and analysis provides a vanishing upper bound on the MSE whenever $p \geq \max \left\{ m^{-1 + \delta}, n^{-1/2 + \delta} \right\}$, also independent of the latent dimension. In fact even as $d$ grows with $m$, as long as $d = o(\log m)$, our analysis guarantees that our algorithm achieves a vanishing MSE. Our analysis relies on the ``local'' structure of the latent space; even if the latent dimension increases, we only need to ensure that there are sufficiently many close neighbor rows so that nearest neighbor averaging produces a good estimate.


\subsection{Future Work}

Our analysis seems tight due to the comparison with minimax rates for noparametric regression, however it is likely that a model which could average entries across both rows and columns would be able to improve the MSE bounds. In particular, the minimax optimal rates for nonparametric regression would imply a lower bound of 
\[	\MSE = \Omega((nmp)^{-\frac{2}{2+2d}}),\]
in the setting where both row and column latent features are sampled uniformly from a $d$-dimensional hypercube. This is achievable with an oracle estimator that were given knowledge of the latent features, but we do not know of information theoretic lower bounds for the MSE that are specific to the latent variable model where features are unobserved. For specific settings such as when the function $f$ when considered as an integral operator has finite spectrum, it is equivalent to low-rank models, for which lower bounds have been characterized. For specific noise models such as the binary observation model which corresponds to the graphon generative model for random graphs, \cite{GaoLuZhou15,KloppTsybakovVerzelen15} show that variants of the least squares estimator achieve optimal rates, but unfortunately they are not polynomial time computable.

From an implementation perspective, the similarity based algorithm proposed in this work, similar to classical collaborative filtering methods, is easy to implement and scales extremely well to large datasets, as it naturally enjoys a parallelizable implementation. Furthermore, the operation of finding $k$ nearest neighbors can benefit from computational advances in building scalable approximate nearest neighbor indices, cf. \cite{indyk1, indyk2}. 

Next we discuss some natural extensions and directions for future work. In our model, the latent function $f$ is assumed to be Lipschitz. However, the proof only truly utilizes the fact that ``locally'' the function value does not oscillate too wildly. Intuitively, this suggests that the result may extend to a broader class of functions, beyond Lipschitz functions. For example, a function with bounded Fourier coefficients does not oscillate too wildly, and thus it may behave well for the purposes of analyzing our algorithm.

Another possible direction for extension is related to the measurement of similarity and the sample complexity. Our current algorithm measures the similarity of rows $u$ and $v$ from their overlapping observed entries, which critically determines the sample complexity requirement of $np^2 \gg1$. However, for sparser datasets without overlaps, we may be able to reveal the similarity by instead comparing distribution signatures such as moments or comparing them through their ``extended" neighborhoods \cite{borgs2017thy}.

As a concluding remark, we would like to mention that the latent variable model is a fairly general model and there is a large body of related applications. Some of the popular recent examples, which are special cases of latent variable model, include Stochastic blockmodels for community detection, the Bradley-Terry model for ranking from pair-wise comparison data and the Dawid-Skene model for low-cost crowd sourcing. Another prominent example of latent variable model is the generative model for random graphs referred to as a Graphon, which has been shown to be the limit of a sequence of graphs. We refer interested readers to \cite[Section 2.4]{Chatterjee15} for an excellent overview on the broad applicability of the latent variable model.

\section{Extending Beyond Matrices to Tensors} \label{sec:tensor_results}
A natural extension beyond matrix completion is the problem of completing a tensor of higher ($>2$) order. Given an unknown tensor $T$ of order $\tau$ with dimensions 
$n_1\times \dots \times n_{\tau}$, suppose that we observe a fraction of its entries corrupted by noise. Similar to matrix completion, the goal in tensor completion is 
to estimate the missing entries in the tensor from the noisy partial observations, as well as to ``de-noise'' the observed entries. 

\subsection{Short Background}
The tensor completion problem is important within a wide variety of applications, including recommendation systems, multi-aspect data mining \cite{Kolda2008, Sun2009}, 
and machine vision \cite{LMWY2013, Zhang2014, inpainting-survey}.

Although tensor completion has been widely studied, there is still a wide gap in understanding, unlike matrix completion. This gap 
partially stems from the hardness of tensor decomposition, as most recovery methods rely on retrieving hidden algebraic structure 
through the framework of low-rank factorization. Tensors do not have a canonical decomposition such as the singular value 
decomposition (SVD) for a matrix. 

There is a factorization scheme, namely the CANDECOMP/PARAFAC (CP) decomposition, which factorizes the tensor as a sum of 
rank-1 tensors (outer product of vectors). However, it is known that finding the rank of a tensor is NP-Complete, which makes it 
computationally intractable. Also, there are known ill-posedness \cite{Silva2008} issues with CP-based low-rank approximation. 

There are other kinds of decompositions such as the Tucker decomposition. Approaches based on Tucker decomposition essentially 
unfold (matricize or flatten) the tensor, and make use of matrix completion theory and methods \cite{Gandy2011, Signoretto2011, 
Tomioka2011, LMWY2013, MHWG2014}.

\subsection{Setup for Tensor Completion}
\subsubsection{Revisiting Our Model}
The nonparametric model presented in Section \ref{sec:model} naturally extends beyond the bivariate case that corresponds to matrices, 
to multivariate setups encompassing higher-order tensors. We recap the latent variable model for a tensor following similar assumptions.

Suppose that there is an unknown $n_1\times \dots \times n_{\tau}$ tensor $T_A$ that we would like to estimate. We observe only a fraction 
of the total $\prod_{i=1}^{\tau} n_i$ entries of $T_A$ with some noise added. Let $\obs \subset [n_1] \times \cdots \times [n_{\tau}]$ denote 
the index set of observed entries. 

Precisely, we consider the following data generation model. Our data tensor $T_Z$ is an $n_1\times \dots \times n_{\tau}$ tensor that is 
generated as follows. Let $T_M \in \{0,1\}^{n_1\times \dots \times n_{\tau}}$ be a binary matrix, which we call the masking tensor. 
We let $T_A \in \Reals^{n_1\times \dots \times n_{\tau}}$ denote the signal tensor and $T_N \in \Reals^{n_1\times \dots \times n_{\tau}}$ denote 
the noise tensor. For each $\valpha = (\alpha_1, \ldots, \alpha_{\tau}) \in [n_1]\times \cdots \times [n_{\tau}]$, 
\begin{equation*}
	T_Z(\valpha) = \begin{cases}
			T_A(\valpha) + T_N(\valpha)	&	\text{when }T_M(\valpha) = 1,\\
			\text{unknown}				&	\text{when }T_M(\valpha) = 0.
		\end{cases}
\end{equation*}

Similar to the model assumptions for the matrix case in Section \ref{sec:lv_model}, we assume $T_A(\valpha)$ is generated 
by the following latent variable model equipped with certain regularity assumptions. 
\begin{itemize}
	\item
	Nonparametric model: there exists a latent function $f$ such that
	\begin{equation*}
		T_A(\valpha) = f( x_1(\alpha_1), \dots, x_{\tau}(\alpha_{\tau}) )	
	\end{equation*}
	for all $\valpha \in [n_1] \times \dots \times [n_{\tau}]$. Here, $x_1(\alpha_1), \dots, x_{\tau}(\alpha_{\tau})$ denote latent variables associated 
	with index $\alpha_i$ in the $i$-th coordinate of the tensor $T_A$, respectively.
	\item
	Regularity Assumptions
	\begin{itemize}
		\item
		For each $i \in [\tau]$, the latent variables $x_i(\alpha_i) \in \cX_i$ for all $\alpha_i \in [n_i]$, where $(\cX_i, d_i)$ is 
		a metric space such that $\diam \cX_i = \sup_{\alpha, \beta \in \cX_i} d_i(\alpha, \beta) \leq D_i$. Moreover, $\cX_i$ 
		is equipped with a Borel probability measure $\mu_i$ and $x_i(\alpha_i)$ is drawn i.i.d. according to $\mu_i$. 
		\item
		Latent function $f$ is bounded, specifically that there exists a constant $\fbound$ such that
		for all $\valpha$, $|f( \valpha)| \leq \fbound$.
		\item Latent function $f$ is $\lip$-Lipschitz in the sense that 
		\[	\big| f( \valpha) - f( \vbeta) \big|	\leq \lip \max_{i \in [\tau]} \big( d_i(\alpha_i, \beta_i) \big).	\]
	\end{itemize}
\end{itemize}
Note that our model assumptions for the noise matrix and the masking matrix are stated in an entrywise fashion and 
readily extends to their tensor analogues, cf. Section \ref{sec:model}.

\subsubsection{Flattening a Tensor to a Matrix}
A $\tau$-order tensor $T \in \Reals^{n_1 \times \cdots \times n_{\tau}}$ can be viewed as a $\tau$-dimensional array of numbers. 
It is possible to `flatten' the tensor $T$ to a matrix (i.e. 2-dimensional array) by rearranging the numbers in the $\tau$-dimensional array.
Formally, given a set $U \subset [\tau]$, we define $\flatten_U(T)$ to be a $\prod_{i \in U} n_i$ by $\prod_{j \in [\tau] \setminus U} n_j$ matrix obtained by flattening the original tensor. 
Without loss of generality\footnote{by taking transpose of $T$}, we may assume $U = [\upsilon]$ for some $0 \leq \upsilon \leq \tau$. We index 
the rows and the columns of $\flatten_U(T)$ using a $\upsilon$-tuple and a $(\tau-\upsilon)$-tuple, respectively. That is to say, for any 
$(\alpha_1, \ldots, \alpha_{\upsilon}) \in [n_1] \times \cdots \times [n_{\upsilon}]$ and $(\bar{\alpha}_1, \ldots, \bar{\alpha}_{\tau-\upsilon}) \in [n_{\upsilon+1}] \times \cdots \times [n_{\tau}]$, 
we have
\begin{align*}
	&\flatten_U(T)\big((\alpha_1, \cdots, \alpha_{\upsilon});(\bar{\alpha}_1, \cdots, \bar{\alpha}_{\tau-\upsilon}) \big)\\
	&\qquad=	T( \alpha_1, \cdots, \alpha_{\upsilon}, \bar{\alpha}_1, \cdots, \bar{\alpha}_{\tau-\upsilon}).
\end{align*}

\subsection{Tensor Completion Algorithm}\label{sec:user_user_tensor}

We describe our generic recipe for tensor completion in Algorithm \ref{alg:algorithm_tensor}. We remark here that 
any matrix estimation algorithm can be used as the matrix estimation subroutine in Step 2.

\begin{algorithm} [h!] 
    \SetKwInOut{Input}{Input}
    \SetKwInOut{Output}{Output}
    \Input{$T_Z \in \Reals^{n_1 \times \dots \times n_{\tau}}; U \subset [\tau]; \eta \geq 0$}	
    \Output{$\hTA \in \Reals^{n_1 \times \dots \times n_{\tau}}$}
    
    \begin{algorithmic}[1]
	\item Flatten $T_Z$ to $\flatten_U(T_Z)$
	\item Run a matrix estimation subroutine (e.g., Algorithm \ref{alg:main_algorithm}) with $\flatten_U(T_Z)$ and $\eta$ 
		to obtain $\widehat{\flatten_U(T_A)}$
	\item Reshape $\widehat{\flatten_U(T_A)}$ to obtain $\widehat{T_A}$
   \end{algorithmic}
    \caption{Generic Description of `Tensor Completion with Flattening' Algorithm}\label{alg:algorithm_tensor}
\end{algorithm}

\subsection{Analysis}
Suppose that Algorithm \ref{alg:user_algorithm} is used as the matrix estimation subroutine in Algorithm \ref{alg:algorithm_tensor}. 
We can obtain an MSE upper bound for the tensor completion algorithm, which is similar to that stated in Theorem \ref{thm:main_mse_user}. We do not include a formal theorem statement and its proof here, but we briefly point out what remains 
unchanged and what needs to be modified in the proof of Theorem \ref{thm:main_mse_user} (cf. Appendix \ref{sec:proof_main_thm}) 
to obtain its counterpart for tensor completion.

Given a $\tau$-order tensor $T_Z \in \Reals^{n_1 \times \dots \times n_{\tau}}$ and an index set $U \subset [\tau]$, we consider the flattened 
matrix, $\flatten_U(T_Z) \in \Reals^{ \prod_{i \in U}n_i \times \prod_{j \in [\tau] \setminus U} n_j }$. Since this matrix is obtained by 
flattening $T_Z$, the `rows' and the `columns' of $\flatten_U(T_Z)$ are not fully exchangeable -- they satisfy only `partial' 
exchangeability induced by the exchangeability in $T_Z$. As a result, we cannot assume the latent variables associated with 
$\prod_{i \in U}n_i$ rows of $\flatten_U(T_Z)$ are drawn i.i.d. from a latent space in the current setup. In fact, there are only $\sum_{i \in U}n_i$ number of independent latent random variables ($n_i$ from $\cX_i$ for each $i \in U$) associated to the $\prod{i \in U}n_i$ rows of the flattened matrix. 
Similarly, the latent variables associated to the columns of $\flatten_U(T_Z)$ are also not sampled i.i.d. as they involve shared coordinates in the original tensor. This difference adds a complication to the analysis of 
Algorithm \ref{alg:algorithm_tensor}, but its effect is limited to the Step 1 in the proof of Lemma \ref{lem:signal_technical}. 

Following the discussion in Appendix \ref{sec:proof_outline}, we observe that there are $\sum_{i=1}^{\tau} n_i + 2 \prod_{i=1}^{\tau} n_i$ 
independent sources of randomness in our model for $T_Z$: $ \big\{ x_i(\alpha_i)\big\}_{i \in [\tau]\atop \alpha_i \in [n_i]}$, 
$\big\{ N(\valpha) \big\}_{\valpha \in [n_1] \times \dots \times [n_{\tau}]}$, $\big\{ M(\valpha) \big\}_{\valpha \in [n_1] \times \dots \times [n_{\tau}]}$.
For the sake of simplicity, we may assume $U = [\upsilon]$ for some $0 \leq \upsilon \leq \tau$. For $\valpha^{(1)} \in [n_1] \times \dots \times [n_{\upsilon}]$ 
and $\valpha^{(2)} \in [n_{\upsilon+1}] \times \dots \times [n_{\tau}]$, we let $\vfeatrow{\valpha^{(1)}} = \big( x_1(\alpha^{(1)}_1), \dots, 
x_{\upsilon}(\alpha^{(1)}_{\upsilon}) \big)$ and $\vfeatcol{\valpha^{(2)}} = \big( x_{\upsilon+1}(\alpha^{(2)}_1), \dots, x_{\tau}(\alpha^{(2)}_{\tau-\upsilon}) \big)$, respectively. 
Also, we let $\vvec{1}$ denote a sequence of an appropriate length with all coordinates being $1$.

By the same exchangeability argument, we can upper bound the MSE by the sum of signal component and the noise component as in the proof of 
Theorem \ref{thm:main_mse_user} -- see Eq. \eqref{eqn:mse_upper.0}, Lemmas \ref{lem:signal_MSE_upper} and \ref{lem:noise_MSE_upper} 
in Appendix \ref{sec:proof_main_thm}):
\begin{align}
	&\MSE(\hTA) =\Exp{\big( \hTA(\vvec{1},\vvec{1}) - T_A(\vvec{1},\vvec{1}) \big)^2 }	\nonumber\\
		&\leq
			\bbE_{\Frandom \setminus \vfeatcol{\vvec{1}} } 
			\bigg[ \max_{(\vvec{v},\vvec{i}) \in \baseest(\vvec{1},\vvec{1}) } \Big\|  f(\vfeatrow{\vvec{v}}, \cdot ) - f(\vfeatrow{\vvec{1}}, \cdot) \Big\|_{L^2}^2	\bigg] \nonumber\\
		&+ C \sigma^2 \bbE \bigg[ \frac{1}{\vNbase} ~\Ind\big( { \vNbase \geq 1 }\big)  \bigg].	\label{eqn:MSE_upper_tensor}
\end{align}
It remains to bound each of the two terms in \eqref{eqn:MSE_upper_tensor}, which can be accomplished by concentration inequalities. 

The second term in \eqref{eqn:MSE_upper_tensor} can be bounded by a very similar argument as in Lemma \ref{lem:noise_technical}, but we need to adjust for the fact that $\vfeatrow{\cdot}$ is associated to a tuple of $U$ latent variables, that are shared across different rows. In Step 1 of the proof or Lemma \ref{lem:noise_technical}, we would use Chernoff bound separately for each dimension of the tensor $i \in U$ to argue that for each entry $x_1(\alpha^{(1)}_i)$, there are sufficiently many ``nearest neighbor'' coordinates with similar latent variables in the $i$-th dimension of the tensor. The number of nearest neighbor rows to $\vfeatrow{\valpha^{(1)}}$ would then be lower bounded by the product of the number of nearest neighbors for each coordinate $i \in U$.

We also need a small modification in bounding the first term in \eqref{eqn:MSE_upper_tensor}. In Step 1 of the proof of Lemma \ref{lem:signal_technical}, we use Bernstein's inequality to prove the concentration of
\begin{align*}
	&\dissimrow(\vvec{1},\vvec{v})\\
	& = \frac{1}{|\basesim(\vvec{1}) \cap \basesim(\vvec{v})|} 
		\sum_{\vvec{j} \in \basesim(\vvec{1}) \cap \basesim(\vvec{v})} \big( T_Z(\vvec{1},\vvec{j}) - T_Z(\vvec{v},\vvec{j}) \big)^2.
\end{align*}
We observe that we cannot use Bernstein's inequality in the current setup of tensor completion, because the summands, 
$ \big( T_Z(\vvec{1},\vvec{j}) - T_Z(\vvec{v},\vvec{j}) \big)^2$, are no longer independent, unlike the setup of matrix completion. 

Nevertheless, the dependence between the summands is still reasonably weak. We may consider $\dissimrow(\vvec{1},\vvec{v})$ 
as a function of the independent random variables, $\big\{ x_i(\alpha_i)\big\}_{i \in [\tau]\atop \alpha_i \in [n_i]}$, 
$\big\{ N(\valpha) \big\}_{\valpha \in [n_1] \times \dots \times [n_{\tau}]}$, $\big\{ M(\valpha) \big\}_{\valpha \in [n_1] \times \dots \times [n_{\tau}]}$.
Then we are able to prove concentration of $\dissimrow(\vvec{1},\vvec{v})$ to its expectation, using a different tool, e.g., 
by modified log-Sobolev inequality. Once we show $\dissimrow(\vvec{1},\vvec{v}) \approx \big\|  f(\featrow{\vvec{v}}, \cdot ) 
- f(\featrow{\vvec{1}}, \cdot) \big\|_{L^2}^2 + 2\sigma^2$ in a similar form as in \eqref{eqn:conc_dissim}, the rest of the proof of 
Lemma \ref{lem:signal_technical} can be reused.

\begin{remark}
In order to quickly get a sense of the appropriate concentration result, we consider the case where $n_i = n$ ($n \geq 2$) for all $i \in [\tau]$, $p = 1$, $|U| = \upsilon$ and the noise is bounded, i.e., $|N(\valpha)| \leq \gamma$ for all $\valpha \in [n]^{\tau}$. We also assume $\diam \cX_i \leq D$ for all $i = \upsilon + 1, \dots, \tau$. Given $\featrow{\vvec{1}}$ 
and $\featrow{\vvec{v}}$, we can observe that $\dissimrow(\vvec{1},\vvec{v})$ is a function of $\big\{ x_i(\alpha_i)\big\}_{i = \upsilon + 1, \ldots, \tau\atop \alpha_i \in [n]}$, 
$\big\{ N(\vvec{1}, \vbeta) \big\}_{\vbeta \in [n]^{\tau-\upsilon}\atop \vbeta \neq \vvec{1}}$, and $\big\{ N(\vvec{v}, \vbeta) \big\}_{\vbeta \in [n]^{\tau-\upsilon}\atop 
\vbeta \neq \vvec{1}}$. Note that $|\basesim(\vvec{1}) \cap \basesim(\vvec{v})| = n^{\tau-\upsilon} - 1$.

We consider the influence that each independent random variable exerts on the function $\dissimrow(\vvec{1},\vvec{v})$. It is easy to verify that 
$\dissimrow(\vvec{1},\vvec{v})$ satisfies the bounded difference property. To verify the property, let us introduce shorthand notations 
$g \big[ x_i(\alpha_i) \big] \triangleq  \dissimrow(\vvec{1},\vvec{v})\big[ x_i(\alpha_i); ~\textrm{the rest}  \big] $ and 
$h \big[ N(\vvec{1}, \vbeta) \big] \triangleq \dissimrow(\vvec{1},\vvec{v})\big[ N(\vvec{1}, \vbeta);~ \textrm{the rest}  \big] $. Then we observe that
\begin{align*}
	\sup_{x_i(\alpha_i), x'_i(\alpha_i)} \Big| g \big[ x_i(\alpha_i) \big] - g \big[ x'_i(\alpha_i) \big] \Big|
		&\leq \frac{n^{\tau-\upsilon-1}-1}{n^{\tau-\upsilon}-1} \lip^2 D^2\\
		&\leq \frac{\lip^2 D^2}{n},\\
	\sup_{N(\vvec{1},\vbeta), N'(\vvec{1},\vbeta)} \Big| h \big[ N(\vvec{1}, \vbeta) \big] - h \big[ N'(\vvec{1}, \vbeta) \big] \Big|
		&\leq \frac{4 \lip D \gamma + 5 \gamma^2}{n^{\tau-\upsilon}-1}.
\end{align*}
Applying the bounded difference inequality (e.g., McDiarmid's), it follows that for any $t > 0$,
\begin{align*}
	&\Prob{ \dissimrow(\vvec{1},\vvec{v}) - \bbE \dissimrow(\vvec{1},\vvec{v}) \geq t}\\
		&\leq \exp \left( - \frac{2t^2}{ (\tau-\upsilon) n \big( \frac{\lip^2 D^2}{n} \big)^2 + 2 (n^{\tau-\upsilon}-1) \big(  \frac{4 \lip D \gamma + 5 \gamma^2}{n^{\tau-\upsilon} -1} \big)^2 } \right)\\
		&\leq \exp \left( - \frac{2 n t^2}{(\tau-\upsilon) \lip^4 D^4} \bigg(1 + \frac{2n}{n^{\tau - \upsilon } -1} \Big( \frac{ 4 \lip D \gamma + 5\gamma^2}{\lip^2 D^2} \Big)^2 \bigg)^{-1} \right).
\end{align*}
Note that we may assume $\tau - \upsilon \geq 1$ and therefore $\frac{2n}{n^{\tau-\upsilon} -1} \leq 4$ (because $n \geq 2$). As a result, we can expect 
$\big| \dissimrow(\vvec{1},\vvec{v}) - \bbE \dissimrow(\vvec{1},\vvec{v}) \big| \lesssim \frac{( \lip^2 D^2 + \gamma^2 ) \sqrt{\tau - \upsilon}}{\sqrt{n}}$. 
This concentration inequality is $\sqrt{\tau - \upsilon}$ times weaker than the result (Bernstein's inequality) from Step 1 of the proof of Lemma \ref{lem:signal_technical}, but its dependence on $n$ remains the same.
\end{remark}

\begin{remark}
The sample complexity of our matrix estimation algorithm requires that $p \geq \max(m^{-1 + \delta}, n^{-1/2 + \delta})$, where $m$ is the number of rows and $n$ is the number of columns. Therefore, the optimal flattening of the tensor that would minimize sample complexity is the flattening that such that 
\[\prod_{i \in U} n_i \approx \bigg(\prod_{j \in [\tau] \setminus U} n_j \bigg)^2.\]
If $n_i = n$ for all $i \in [\tau]$, then the optimal flattening would result in matrix dimensions of $n^{\lfloor\tau/3\rfloor} \times n^{\lceil 2\tau/3\rceil}$. In this case the natural extension of our analysis should result in MSE convergence with sample complexity of $p \geq n^{-\lfloor\tau/3\rfloor + \delta}$.
\end{remark}

\section{Experiments} \label{sec:exp}

In this section we present experimental results from applying the User-Item Gaussian Kernel variant of our algorithm to real datasets. We state in Algorithm \ref{alg:user_gauss_algorithm} the specific algorithm variant that is used in the experiments. We did not sample split as it was primarily used for the analysis but is not necessary in practice. The implemented algorithm adds a few modifications:
\begin{itemize}
\item $\dissimrow(u,v)$ is computed according to the {\em sample variance} of the differences between two rows
\item we use Gausian kernel weights, and we combine both row and column dissimilarities by taking the maximum
\item we use an estimator $A(u,i) \approx (Z(v,i) + Z(u,j) - Z(v,j))$ motivated by first order Taylor series approximation
\end{itemize}
The resulting algorithm is similar to the mean-adjusted variant of collaborative filtering, except for the addition of combining both row and column dissimilarities.

\begin{algorithm} [h!] 
    \SetKwInOut{Input}{Input}
    \SetKwInOut{Output}{Output}
    \Input{$Z \in \Reals^{m \times n}; \lambda \geq 0, \beta \in \IntegersP$}	
    \Output{$\hA$}
    
    \begin{algorithmic}[1]
		\item
		For each $v \in [m]$ and for each $i \in [n]$, define 
		\begin{align*}
			\basesim^{\text{row}}(v) &\triangleq \{ j \in [n]: (v,j) \in \Omega  \} \quad\text{and }\\
			\basesim^{\text{col}}(i) &\triangleq \{ v \in [n]: (v,i) \in \Omega  \}.
		\end{align*}
		Let us denote $Q^{\text{row}}(u,v) = \basesim^{\text{row}}(u) \cap \basesim^{\text{row}}(v)$ and $Q^{\text{col}}(i,j) = \basesim^{\text{col}}(i) \cap \basesim^{\text{col}}(j)$.
		\item
		For each $(u,v) \in [m]^2$, estimate the dissimilarity between two rows with the sample variance between common entries 
		\begin{align*}
			&\dissimrow(u,v) = \frac{1}{2 | Q^{\text{row}}(u,v)| ( |Q^{\text{row}}(u,v)| - 1)}\\
			&\times \sum_{(i,j) \in (Q^{\text{row}}(u,v))^2} \left( (Z_{u,i} - Z_{v,i}) - (Z_{u,j} - Z{v,j}) \right)^2.
		\end{align*}
		
		\item
		For each $(i,j) \in [n]^2$, estimate the dissimilarity between two columns with the sample variance between common entries 
		\begin{align*}
			&\dissimcol(i,j) = \frac{1}{2 |Q^{\text{col}}(i,j)| ( |Q^{\text{col}}(i,j)| - 1)} \\
			&\times \sum_{(u,v) \in (Q^{\text{col}}(i,j))^2} \left((Z_{u,i} - Z_{u,j})- (Z_{v,i} - Z_{v,j}) \right)^2.
		\end{align*}
		
		\item 
		Define the weights
		\begin{align*}
			w_{ui}(v, j) &= e^{-\lambda ~\max \{ \dissimrow(u,v), ~\dissimcol(i,j) \} }\\
				&\qquad\times \Ind_{\{|Q^{\text{row}}(u,v)| \geq \beta\}} \Ind_{\{|Q^{\text{col}}(i,j)| \geq \beta\}}
		\end{align*}	
		
		\item Compute the estimate for $(u,i)$ according to
		\begin{equation*}
			\hA(u,i) =  \frac{\sum_{(v,j)} w_{ui}(v,j) (Z(v,i) + Z(u,j) - Z(v,j)) }{ \sum_{(v,j)} w_{ui}(v,j)}.
		\end{equation*}
   \end{algorithmic}

    \caption{User-user Gaussian Kernel Algorithm with First-Order Estimates}\label{alg:user_gauss_algorithm}
\end{algorithm}


\subsection{Matrix completion}

    We evaluated the performance of our algorithm on predicting user-movie ratings for the MovieLens 1M and Netflix datasets. 
    We chose the overlap parameter $\beta = 2$ to ensure the algorithm is able to compute an estimate for all missing entries. When $\beta$ is larger, the algorithm enforces rows (or columns) to have more commonly rated movies (or users). Although this increases the reliability of the estimates, it also reduces the fraction of entries for which the estimate is defined. 

    We compared our method with user-user collaborative filtering, item-item collaborative filtering, and SoftImpute from \cite{mazumder2010tensor}. We chose the classic mean-adjusted collaborative filtering method, in which the weights are proportional to the cosine similarity of pairs of users or items (i.e. movies). SoftImpute is a matrix-factorization-based method which iteratively replaces missing elements in the matrix with those obtained from a soft-thresholded SVD.

\begin{figure}[ht]
    \centering
        \includegraphics[width=\linewidth]{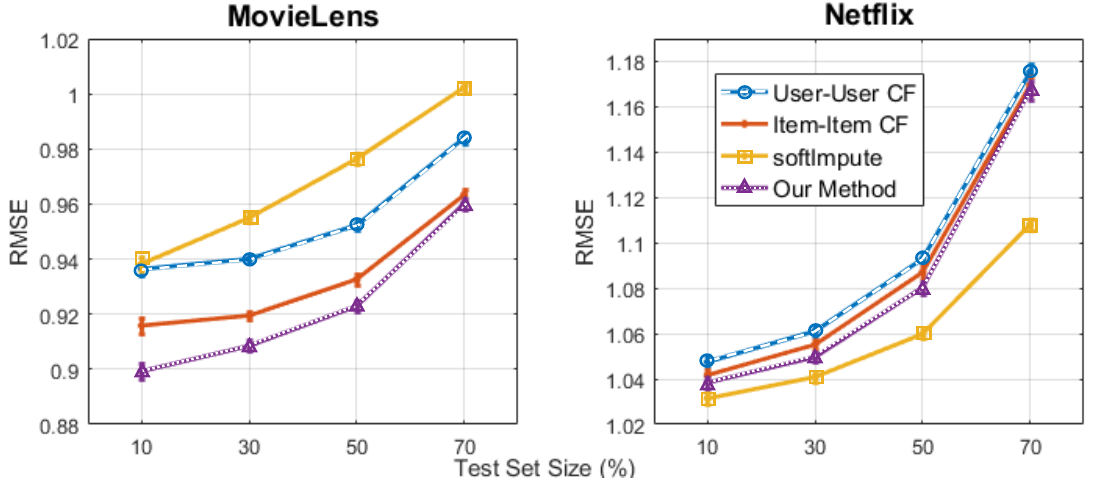}
        \caption{Performance of algorithms on Netflix and MovieLens datasets with 95$\%$ confidence interval.  $\lambda$ values used by our algorithm are 2.8 (10$\%$), 2.3 (30$\%$), 1.7 (50$\%$), 1 (70$\%$) for MovieLens, and  1.8 (10$\%$), 1.7 (30$\%$), 1.6 (50$\%$), 1.5 (70$\%$) for Netflix. }
    \label{fig:result}
    \end{figure}

    The MovieLens 1M data set contains about 1 million ratings by 6000 users of 4000 movies. The Netflix data set consists of about 100 million movie ratings by 480,189 users of about 17,770 movies. For both MovieLens and Netflix data sets, the ratings are integers from 1 to 5. From each dataset, we generated 100 smaller user-movie rating matrices, in which we randomly subsampled 2000 users and 2000 movies. 
    For each rating matrix, we randomly select and withhold a percentage of the known ratings for the test set, while the remaining portion of the data set is revealed to the algorithm for computing the estimates (or training). After the algorithm computes its predictions for all the missing user-movie pairs, we evaluate the Root Mean Squared Error (RMSE) of the predictions compared to the ratings from the withheld test set. Figure \ref{fig:result} plots the RMSE of our method along with classic collaborative filtering and SoftImpute evaluated against $10\%$, $30\%$, $50 \%$, and $70 \%$ withheld test sets. The RMSE is averaged over 100 subsampled rating matrices, and $95\%$ confidence intervals are provided.
    
Figure \ref{fig:result} suggests that our algorithm achieves a systematic improvement over classical user-user and item-item collaborative filtering. SoftImpute performs worse than all methods on the MovieLens dataset, but it performs better than all methods on the Netflix dataset. 


\subsection{Tensor completion}
\begin{figure*}[t]
\captionsetup[subfigure]{labelformat=empty}

    
         
            \begin{subfigure}[h]{0.15\textwidth}
            \caption{Original}
        \includegraphics[width=\textwidth]{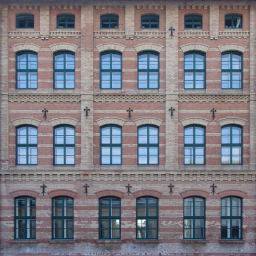}
        	\caption{}
        \label{fig:facade}
    \end{subfigure}
       \hfill 
    \begin{subfigure}[h]{0.15\textwidth}
            \caption{Degraded}
        \includegraphics[width=\textwidth]{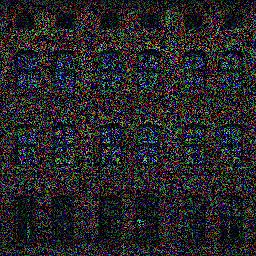}
        	   \caption{}
        \label{fig:facade_prediction70}
    \end{subfigure}
   \hfill 
    \begin{subfigure}[h]{0.15\textwidth}
            \caption{FaLRTC}
        \includegraphics[width=\textwidth]{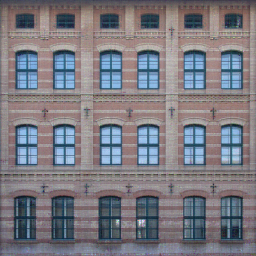}
                \caption{0.0924}
        \label{fig:facade_FaLRTC70}
    \end{subfigure}
       \hfill 
    \begin{subfigure}[h]{0.15\textwidth}
            \caption{TenAlt}
        \includegraphics[width=\textwidth]{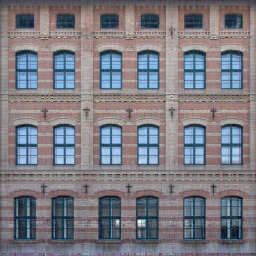}
                \caption{0.099}
        \label{fig:facade_alt70}
    \end{subfigure}
           \hfill 
    \begin{subfigure}[h]{0.15\textwidth}
            \caption{FBCP}
        \includegraphics[width=\textwidth]{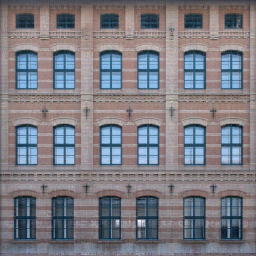}
                \caption{0.12}
        \label{fig:facade_FBCP70}
    \end{subfigure}
          \hfill 
    \begin{subfigure}[h]{0.15\textwidth}
        \caption{Our Method}
        \includegraphics[width=\textwidth]{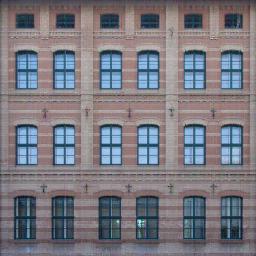}
                \caption{0.0869}
        \label{fig:facade_tensor70}
    \end{subfigure}

        \begin{subfigure}[h]{0.15\textwidth}
        \includegraphics[width=\textwidth]{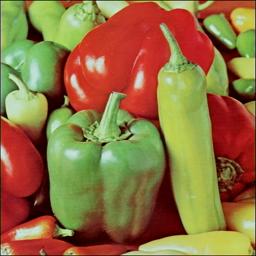}
        \caption{}
        \label{fig:pepper}
    \end{subfigure}
   \hfill 
        \begin{subfigure}[h]{0.15\textwidth}
        \includegraphics[width=\textwidth]{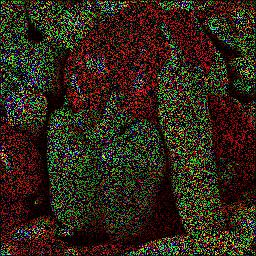}
        \caption{}
        \label{fig:pepper_predictio70}
    \end{subfigure}
   \hfill 
    \begin{subfigure}[h]{0.15\textwidth}
        \includegraphics[width=\textwidth]{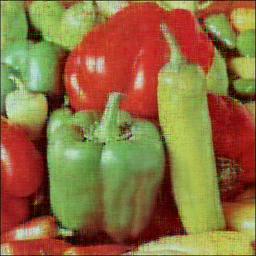}
        \caption{0.1101}
        \label{fig:pepper_FaLRTC70}
    \end{subfigure}
       \hfill 
    \begin{subfigure}[h]{0.15\textwidth}
        \includegraphics[width=\textwidth]{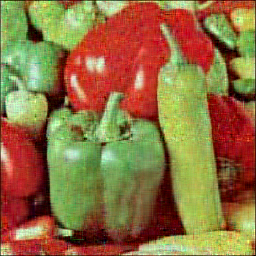}
        \caption{0.1182}
        \label{fig:pepper_alt70}
    \end{subfigure}
           \hfill 
    \begin{subfigure}[h]{0.15\textwidth}
        \includegraphics[width=\textwidth]{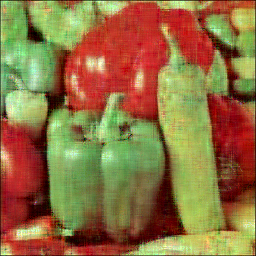}
        \caption{0.154}
        \label{fig:pepper_FBCP70}
    \end{subfigure}
          \hfill 
    \begin{subfigure}[h]{0.15\textwidth}
        \includegraphics[width=\textwidth]{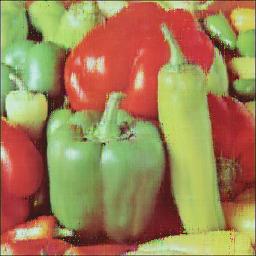}
        \caption{0.109}
        \label{fig:pepper_tensor70}
    \end{subfigure}    
        
        \caption{Recovery results for Facade and Pepper images with 70$\%$ of missing entries. RSE is reported under  the recovery images.}
      \label{fig:facade_pepper_visual}
\end{figure*}

   We consider the problem of image inpainting for evaluating the performance of tensor completion algorithm.  Inpainting is the process of reconstructing lost or deteriorated parts of image or videos. Such methods, in particular, have revitalized the process of recovery old artifacts in museum world which was historically done by conservators or art restorers. An interested reader is referred to a recent survey \cite{inpainting-survey} for summary of the state of art on methods and techniques. We compare performance of our algorithm against existing methods in the literature on the image inpainting problem. 
Figure \ref{fig:facade_pepper_visual} shows a sample of the image inpainting results for the facade and pepper images when $70\%$ of the pixels are removed.

 An image can be represented as a $3^{\text{rd}}$-order tensor where the dimensions are rows $\times$ columns $\times$ RGB. In particular we used three images (Lenna, Pepper, and Facade) of dimensions 256 $\times$ 256 $\times$ 3. For each image, a percentage of the pixels are randomly removed, and the missing entries are filled in by various tensor completion algorithms.

    For the implementation of our tensor completion method, we collapsed the last two dimensions of the tensor (columns and RGB) to reduce the image to a matrix, and applied our method. 
    We set the overlap parameter $\beta = 2$. 
        We compared our method against fast low rank tensor completion (FaLRTC) \cite{liu2013}, alternating minimization for tensor completion (TenAlt) \cite{jain2014provable}, and fully Bayesian CP factorization (FBCP) \cite{fbcp}, which extends the CANDECOMP/PARAFAC(CP) tensor factorization with automatic tensor rank determination.

\begin{figure*}[h]
    \centering
        \includegraphics[width=0.8\linewidth]{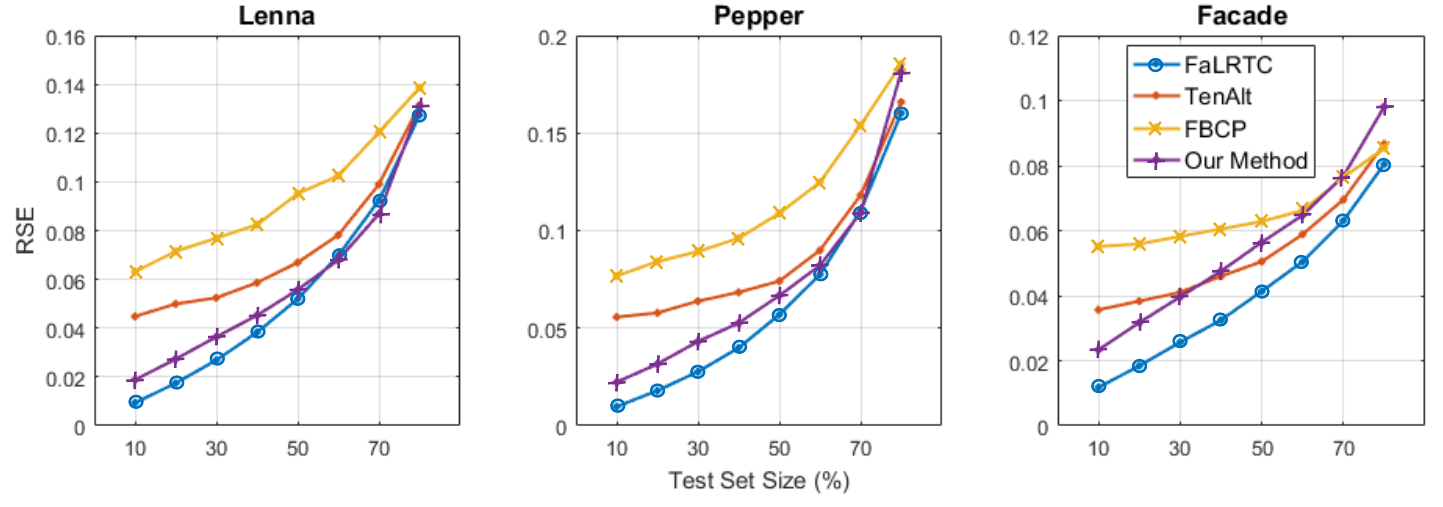}
        \caption{Performance comparison between different tensor completion algorithms based on RSE vs testing set size. For our method, we set overlap parameter $\beta$ to 2.
}
    \label{fig:tensor_result}
    \end{figure*}

    To evaluate the outputs produced by each method, we computed the relative squared error (RSE), defined as
        \[\text{RSE} = \frac{\sum_{i,j,k \in E} (\hat{Z}(i, j, k) - Z(i, j, k))^2}{\sum_{i,j,k \in E} (Z(i,j,k) - \bar{Z})^2}, \]
    where $\bar{Z}$ is the average value of the true entries. Figure $\ref{fig:tensor_result}$ plots the RSE achieved by each tensor completion method on the three images, as a function of the percentage of pixels removed. The results demonstrate that our tensor completion method is competitive with existing tensor factorization based approaches, while maintaining a naive simplicity. In short, a simple algorithm works nearly as good as the best algorithm for this problem!
    
\section*{Acknowledgment}
This work is supported in parts by ARO under MURI award 133668-5079809, by NSF under grants CMMI-1462158, CMMI-1634259 and TRIPODS, by DARPA under grant W911NF-16-1-0551, and additionally by a Samsung Scholarship, Siebel Scholarship, NSF Graduate Fellowship, and Claude E. Shannon Research Assistantship.

\bibliographystyle{IEEEtran}
\bibliography{BR_bibliography}

\begin{thebibliography}{10}
\providecommand{\url}[1]{#1}
\csname url@samestyle\endcsname
\providecommand{\newblock}{\relax}
\providecommand{\bibinfo}[2]{#2}
\providecommand{\BIBentrySTDinterwordspacing}{\spaceskip=0pt\relax}
\providecommand{\BIBentryALTinterwordstretchfactor}{4}
\providecommand{\BIBentryALTinterwordspacing}{\spaceskip=\fontdimen2\font plus
\BIBentryALTinterwordstretchfactor\fontdimen3\font minus
  \fontdimen4\font\relax}
\providecommand{\BIBforeignlanguage}[2]{{%
\expandafter\ifx\csname l@#1\endcsname\relax
\typeout{** WARNING: IEEEtran.bst: No hyphenation pattern has been}%
\typeout{** loaded for the language `#1'. Using the pattern for}%
\typeout{** the default language instead.}%
\else
\language=\csname l@#1\endcsname
\fi
#2}}
\providecommand{\BIBdecl}{\relax}
\BIBdecl

\bibitem{chenshah}
G.~Chen and D.~Shah, \emph{Explaining the success of nearest neighbor method in
  prediction}.\hskip 1em plus 0.5em minus 0.4em\relax Foundations and Trends in
  Machine Learning, 2018.

\bibitem{mack1982weak}
Y.~Mack and B.~W. Silverman, ``Weak and strong uniform consistency of kernel
  regression estimates,'' \emph{Zeitschrift f{\"u}r Wahrscheinlichkeitstheorie
  und verwandte Gebiete}, vol.~61, no.~3, pp. 405--415, 1982.

\bibitem{WandJones94}
M.~P. Wand and M.~C. Jones, \emph{Kernel smoothing}.\hskip 1em plus 0.5em minus
  0.4em\relax Crc Press, 1994.

\bibitem{srebro2004generalization}
N.~Srebro, N.~Alon, and T.~S. Jaakkola, ``Generalization error bounds for
  collaborative prediction with low-rank matrices,'' in \emph{Advances In
  Neural Information Processing Systems}, 2004, pp. 1321--1328.

\bibitem{candes2009exact}
E.~J. Cand{\`e}s and B.~Recht, ``Exact matrix completion via convex
  optimization,'' \emph{Foundations of Computational mathematics}, vol.~9,
  no.~6, pp. 717--772, 2009.

\bibitem{rohde2011estimation}
A.~Rohde, A.~B. Tsybakov \emph{et~al.}, ``Estimation of high-dimensional
  low-rank matrices,'' \emph{The Annals of Statistics}, vol.~39, no.~2, pp.
  887--930, 2011.

\bibitem{keshavan56matrix}
R.~Keshavan, A.~Montanari, and S.~Oh, ``Matrix completion from a few entries,''
  \emph{IEEE Trans. Inf. Theory}, vol.~56, no.~6, 2009.

\bibitem{negahban2012restricted}
S.~Negahban and M.~J. Wainwright, ``Restricted strong convexity and weighted
  matrix completion: Optimal bounds with noise,'' \emph{The Journal of Machine
  Learning Research}, vol.~13, no.~1, pp. 1665--1697, 2012.

\bibitem{jain2013low}
P.~Jain, P.~Netrapalli, and S.~Sanghavi, ``Low-rank matrix completion using
  alternating minimization,'' in \emph{Proceedings of the 45th annual ACM
  symposium on Theory of computing}.\hskip 1em plus 0.5em minus 0.4em\relax
  ACM, 2013, pp. 665--674.

\bibitem{Fazel03}
M.~Fazel, H.~Hindi, and S.~P. Boyd, ``Log-det heuristic for matrix rank
  minimization with applications to hankel and euclidean distance matrices,''
  in \emph{Proceedings of ACC}, vol.~3.\hskip 1em plus 0.5em minus 0.4em\relax
  IEEE, 2003, pp. 2156--2162.

\bibitem{LiuVandenberghe10}
Z.~Liu and L.~Vandenberghe, ``Interior-point method for nuclear norm
  approximation with application to system identification,'' \emph{SIAM Journal
  on Matrix Analysis and Applications}, vol.~31, no.~3, pp. 1235--1256, 2010.

\bibitem{Cai08}
D.~Cai, X.~He, X.~Wu, and J.~Han, ``Non-negative matrix factorization on
  manifold,'' in \emph{Data Mining, 2008. ICDM'08. Eighth IEEE International
  Conference on}.\hskip 1em plus 0.5em minus 0.4em\relax IEEE, 2008, pp.
  63--72.

\bibitem{Lin09}
Z.~Lin, A.~Ganesh, J.~Wright, L.~Wu, M.~Chen, and Y.~Ma, ``Fast convex
  optimization algorithms for exact recovery of a corrupted low-rank matrix,''
  \emph{CAMSAP}, vol.~61, 2009.

\bibitem{Shen09}
B.-H. Shen, S.~Ji, and J.~Ye, ``Mining discrete patterns via binary matrix
  factorization,'' in \emph{Proceedings of the 15th ACM SIGKDD international
  conference}.\hskip 1em plus 0.5em minus 0.4em\relax ACM, 2009, pp. 757--766.

\bibitem{Mazumder10}
R.~Mazumder, T.~Hastie, and R.~Tibshirani, ``Spectral regularization algorithms
  for learning large incomplete matrices,'' \emph{The Journal of Machine
  Learning Research}, vol.~11, pp. 2287--2322, 2010.

\bibitem{XuMassoulieLelarge14}
J.~Xu, L.~Massouli{\'e}, and M.~Lelarge, ``Edge label inference in generalized
  stochastic block models: from spectral theory to impossibility results.'' in
  \emph{COLT}, 2014, pp. 903--920.

\bibitem{GantiBalzanoWillett2015}
R.~S. Ganti, L.~Balzano, and R.~Willett, ``Matrix completion under monotonic
  single index models,'' in \emph{Advances in Neural Information Processing
  Systems}, 2015, pp. 1864--1872.

\bibitem{LeeKimLebanonSingerBengio16}
\BIBentryALTinterwordspacing
J.~Lee, S.~Kim, G.~Lebanon, Y.~Singer, and S.~Bengio, ``Llorma: Local low-rank
  matrix approximation,'' \emph{Journal of Machine Learning Research}, vol.~17,
  no.~15, pp. 1--24, 2016. [Online]. Available:
  \url{http://jmlr.org/papers/v17/14-301.html}
\BIBentrySTDinterwordspacing

\bibitem{Chatterjee15}
S.~Chatterjee, ``Matrix estimation by universal singular value thresholding,''
  \emph{The Annals of Statistics}, vol.~43, no.~1, pp. 177--214, 2015.

\bibitem{Xu2017}
J.~Xu, ``Rates of convergence of spectral methods for graphon estimation,'' in
  \emph{ICML}, 2017.

\bibitem{goldberg92}
D.~Goldberg, D.~Nichols, B.~M. Oki, and D.~Terry, ``Using collaborative
  filtering to weave an information tapestry,'' \emph{Commun. ACM}, 1992.

\bibitem{Linden03}
G.~Linden, B.~Smith, and J.~York, ``Amazon.com recommendations: Item-to-item
  collaborative filtering,'' \emph{IEEE Internet Computing}, vol.~7, no.~1, pp.
  76--80, 2003.

\bibitem{korenHandbook}
Y.~Koren and R.~Bell, ``Advances in collaborative filtering,'' in
  \emph{Recommender Systems Handbook}.\hskip 1em plus 0.5em minus 0.4em\relax
  Springer US, 2011, pp. 145--186.

\bibitem{Ning2015}
X.~Ning, C.~Desrosiers, and G.~Karypis, \emph{Recommender Systems
  Handbook}.\hskip 1em plus 0.5em minus 0.4em\relax Springer US, 2015, ch. A
  Comprehensive Survey of Neighborhood-Based Recommendation Methods, pp.
  37--76.

\bibitem{BellKoren07}
\BIBentryALTinterwordspacing
R.~M. Bell and Y.~Koren, ``Scalable collaborative filtering with jointly
  derived neighborhood interpolation weights,'' in \emph{Proceedings of the
  2007 Seventh IEEE International Conference on Data Mining}, ser. ICDM
  '07.\hskip 1em plus 0.5em minus 0.4em\relax Washington, DC, USA: IEEE
  Computer Society, 2007, pp. 43--52. [Online]. Available:
  \url{http://dx.doi.org/10.1109/ICDM.2007.90}
\BIBentrySTDinterwordspacing

\bibitem{Koren08}
\BIBentryALTinterwordspacing
Y.~Koren, ``Factorization meets the neighborhood: A multifaceted collaborative
  filtering model,'' in \emph{Proceedings of the 14th ACM SIGKDD International
  Conference on Knowledge Discovery and Data Mining}, ser. KDD '08.\hskip 1em
  plus 0.5em minus 0.4em\relax New York, NY, USA: ACM, 2008, pp. 426--434.
  [Online]. Available: \url{http://doi.acm.org/10.1145/1401890.1401944}
\BIBentrySTDinterwordspacing

\bibitem{WangDeVriesReinders06}
\BIBentryALTinterwordspacing
J.~Wang, A.~P. de~Vries, and M.~J.~T. Reinders, ``Unifying user-based and
  item-based collaborative filtering approaches by similarity fusion,'' in
  \emph{Proceedings of the 29th Annual International ACM SIGIR Conference on
  Research and Development in Information Retrieval}, ser. SIGIR '06.\hskip 1em
  plus 0.5em minus 0.4em\relax New York, NY, USA: ACM, 2006, pp. 501--508.
  [Online]. Available: \url{http://doi.acm.org/10.1145/1148170.1148257}
\BIBentrySTDinterwordspacing

\bibitem{BreslerChenShah14}
G.~Bresler, G.~H. Chen, and D.~Shah, ``A latent source model for online
  collaborative filtering,'' in \emph{Advances in Neural Information Processing
  Systems}, 2014, pp. 3347--3355.

\bibitem{BreslerShahVoloch15}
G.~Bresler, D.~Shah, and L.~F. Voloch, ``Collaborative filtering with low
  regret,'' in \emph{ACM Sigmetrics}, 2016.

\bibitem{AiroldiCostaChan13}
E.~M. Airoldi, T.~B. Costa, and S.~H. Chan, ``Stochastic blockmodel
  approximation of a graphon: Theory and consistent estimation,'' in
  \emph{Advances in Neural Information Processing Systems}, 2013, pp. 692--700.

\bibitem{ZhangLevinaZhu15}
Y.~Zhang, E.~Levina, and J.~Zhu, ``Estimating network edge probabilities by
  neighbourhood smoothing,'' \emph{Biometrika}, vol. 104, no.~4, pp. 771--783,
  2017.

\bibitem{kolda2009tensor}
T.~G. Kolda and B.~W. Bader, ``Tensor decompositions and applications,''
  \emph{SIAM review}, vol.~51, no.~3, pp. 455--500, 2009.

\bibitem{anandkumar2014tensor}
A.~Anandkumar, R.~Ge, D.~Hsu, S.~M. Kakade, and M.~Telgarsky, ``Tensor
  decompositions for learning latent variable models,'' \emph{The Journal of
  Machine Learning Research}, vol.~15, no.~1, pp. 2773--2832, 2014.

\bibitem{jain2014provable}
P.~Jain and S.~Oh, ``Provable tensor factorization with missing data,'' in
  \emph{Advances in Neural Information Processing Systems}, 2014, pp.
  1431--1439.

\bibitem{oh2014learning}
S.~Oh and D.~Shah, ``Learning mixed multinomial logit model from ordinal
  data,'' in \emph{Advances in Neural Information Processing Systems}, 2014,
  pp. 595--603.

\bibitem{liu2013tensor}
J.~Liu, P.~Musialski, P.~Wonka, and J.~Ye, ``Tensor completion for estimating
  missing values in visual data,'' \emph{IEEE transactions on pattern analysis
  and machine intelligence}, vol.~35, no.~1, pp. 208--220, 2013.

\bibitem{gandy2011tensor}
S.~Gandy, B.~Recht, and I.~Yamada, ``Tensor completion and low-n-rank tensor
  recovery via convex optimization,'' \emph{Inverse Problems}, vol.~27, no.~2,
  p. 025010, 2011.

\bibitem{tomioka2010estimation}
R.~Tomioka, K.~Hayashi, and H.~Kashima, ``Estimation of low-rank tensors via
  convex optimization,'' \emph{arXiv preprint arXiv:1010.0789}, 2010.

\bibitem{tomioka2011statistical}
R.~Tomioka, T.~Suzuki, K.~Hayashi, and H.~Kashima, ``Statistical performance of
  convex tensor decomposition,'' in \emph{Advances in Neural Information
  Processing Systems}, 2011, pp. 972--980.

\bibitem{bhojanapalli2015new}
S.~Bhojanapalli and S.~Sanghavi, ``A new sampling technique for tensors,''
  \emph{arXiv preprint arXiv:1502.05023}, 2015.

\bibitem{BarakMoitra16}
B.~Barak and A.~Moitra, ``Noisy tensor completion via the sum-of-squares
  hierarchy,'' in \emph{Conference on Learning Theory}, 2016.

\bibitem{PotechinSteurer17}
A.~Potechin and D.~Steurer, ``Exact tensor completion with sum-of-squares,'' in
  \emph{Conference on Learning Theory}, 2017, pp. 1619--1673.

\bibitem{xia2017statistically}
D.~Xia, M.~Yuan, and C.-H. Zhang, ``Statistically optimal and computationally
  efficient low rank tensor completion from noisy entries,'' \emph{arXiv
  preprint arXiv:1711.04934}, 2017.

\bibitem{xia2017polynomial}
D.~Xia and M.~Yuan, ``On polynomial time methods for exact low-rank tensor
  completion,'' \emph{Foundations of Computational Mathematics}, pp. 1--49,
  2017.

\bibitem{MontanariSun18}
A.~Montanari and N.~Sun, ``Spectral algorithms for tensor completion,''
  \emph{Communications on Pure and Applied Mathematics}, vol.~71, no.~11, pp.
  2381--2425, 2018.

\bibitem{arora2012learning}
S.~Arora, R.~Ge, and A.~Moitra, ``Learning topic models--going beyond svd,'' in
  \emph{Foundations of Computer Science (FOCS), 2012 IEEE 53rd Annual Symposium
  on}.\hskip 1em plus 0.5em minus 0.4em\relax IEEE, 2012, pp. 1--10.

\bibitem{arora2012computing}
S.~Arora, R.~Ge, R.~Kannan, and A.~Moitra, ``Computing a nonnegative matrix
  factorization--provably,'' in \emph{Proceedings of the 44th annual ACM
  symposium on Theory of computing}.\hskip 1em plus 0.5em minus 0.4em\relax
  ACM, 2012, pp. 145--162.

\bibitem{Aldous81}
D.~Aldous, ``Representations for partially eschangeable arrays of random
  variables,'' \emph{J. Multivariate Anal.}, vol.~11, pp. 581 -- 598, 1981.

\bibitem{Hoover82}
D.~Hoover, ``Row-column exchangeability and a generalized model for
  probability,'' in \emph{Exchangeability in Probability and Statistics (Rome,
  1981)}, 1981, pp. 281 -- 291.

\bibitem{Austin12}
T.~Austin, ``Exchangeable random arrays.'' \emph{Technical Report, Notes for
  IAS workshop.}, 2012.

\bibitem{OrbRoy2015}
P.~Orbanz and D.~M. Roy, ``Bayesian models of graphs, arrays and other
  exchangeable random structures,'' \emph{IEEE transactions on pattern analysis
  and machine intelligence}, vol.~37, no.~2, pp. 437 -- 461, 2015.

\bibitem{balakrishnan2011statistical}
S.~Balakrishnan, M.~Kolar, A.~Rinaldo, A.~Singh, and L.~Wasserman,
  ``Statistical and computational tradeoffs in biclustering,'' in \emph{NIPS
  2011 workshop on computational trade-offs in statistical learning}, vol.~4,
  2011.

\bibitem{ma2015computational}
Z.~Ma, Y.~Wu \emph{et~al.}, ``Computational barriers in minimax submatrix
  detection,'' \emph{The Annals of Statistics}, vol.~43, no.~3, pp. 1089--1116,
  2015.

\bibitem{abbe2016exact}
E.~Abbe, A.~S. Bandeira, and G.~Hall, ``Exact recovery in the stochastic block
  model,'' \emph{IEEE Transactions on Information Theory}, vol.~62, no.~1, pp.
  471--487, 2016.

\bibitem{shah2017stochastically}
N.~B. Shah, S.~Balakrishnan, A.~Guntuboyina, and M.~J. Wainwright,
  ``Stochastically transitive models for pairwise comparisons: Statistical and
  computational issues,'' \emph{IEEE Transactions on Information Theory},
  vol.~63, no.~2, pp. 934--959, 2017.

\bibitem{chatterjee2019estimation}
S.~Chatterjee and S.~Mukherjee, ``Estimation in tournaments and graphs under
  monotonicity constraints,'' \emph{IEEE Transactions on Information Theory},
  2019.

\bibitem{flammarion2019optimal}
N.~Flammarion, C.~Mao, P.~Rigollet \emph{et~al.}, ``Optimal rates of
  statistical seriation,'' \emph{Bernoulli}, vol.~25, no.~1, pp. 623--653,
  2019.

\bibitem{shah2018low}
N.~B. Shah, S.~Balakrishnan, and M.~J. Wainwright, ``Low permutation-rank
  matrices: Structural properties and noisy completion,'' in \emph{2018 IEEE
  International Symposium on Information Theory (ISIT)}.\hskip 1em plus 0.5em
  minus 0.4em\relax IEEE, 2018, pp. 366--370.

\bibitem{lee2016blind}
C.~E. Lee, Y.~Li, D.~Shah, and D.~Song, ``Blind regression: nonparametric
  regression for latent variable models via collaborative filtering,'' in
  \emph{Proceedings of the 30th International Conference on Neural Information
  Processing Systems}.\hskip 1em plus 0.5em minus 0.4em\relax Curran Associates
  Inc., 2016, pp. 2163--2173.

\bibitem{GaoLuZhou15}
C.~Gao, Y.~Lu, and H.~H. Zhou, ``Rate-optimal graphon estimation,'' \emph{The
  Annals of Statistics}, vol.~43, no.~6, pp. 2624--2652, 2015.

\bibitem{KloppTsybakovVerzelen15}
O.~Klopp, A.~B. Tsybakov, and N.~Verzelen, ``Oracle inequalities for network
  models and sparse graphon estimation,'' \emph{Annals of Statistics}, 2015.

\bibitem{indyk1}
P.~Indyk, ``Algorithmic applications of low-distortion geometric embeddings,''
  in \emph{focs}, vol.~1, 2001, pp. 10--33.

\bibitem{indyk2}
------, ``Nearest neighbors in high-dimensional spaces,'' 2004.

\bibitem{borgs2017thy}
C.~Borgs, J.~Chayes, C.~E. Lee, and D.~Shah, ``Thy friend is my friend:
  Iterative collaborative filtering for sparse matrix estimation,'' in
  \emph{Advances in Neural Information Processing Systems}, 2017, pp.
  4715--4726.

\bibitem{Kolda2008}
T.~G. Kolda and J.~Sun, ``Scalable tensor decompositions for multi-aspect data
  mining,'' in \emph{2008 Eighth IEEE International conference on data mining},
  2008, pp. 363 -- 372.

\bibitem{Sun2009}
J.~Sun, S.~Papadimitriou, C.~Y. Lin, N.~Cao, S.~Liu, and W.~Qian, ``Multivis:
  Content-based social network exploration through multi-way visual analysis,''
  in \emph{Proc. SIAM Intl. Conf. on Data Mining}, 2009, pp. 1064 -- 1075.

\bibitem{LMWY2013}
J.~Liu, P.~Musialski, P.~Wonka, and J.~Ye, ``Tensor completion for estimating
  missing values in visual data,'' \emph{IEEE Transactions on Pattern Analysis
  and Machine Intelligence}, vol.~35, no.~1, pp. 208 -- 220, 2013.

\bibitem{Zhang2014}
Z.~Zhang, G.~Ely, S.~Aeron, N.~Hao, and M.~Kilmer, ``Novel methods for
  multilinear data completion and denoising based on tensor-svd,'' in
  \emph{Proc. IEEE Conf. on CVPR}, 2014, pp. 3842 -- 3849.

\bibitem{inpainting-survey}
S.~Ravi, P.~Pasupathi, S.~Muthukumar, and N.~Krishnan, ``Image in-painting
  techniques-a survey and analysis,'' in \emph{Innovations in Information
  Technology (IIT), 2013 9th International Conference on}.\hskip 1em plus 0.5em
  minus 0.4em\relax IEEE, 2013, pp. 36--41.

\bibitem{Silva2008}
V.~De~Silva and L.~H. Lim, ``Tensor rank and the ill-posedness of the best
  low-rank approximation problem,'' \emph{SIAM Journal on Matrix Analysis and
  Applications}, vol.~30, no.~3, pp. 1084 -- 1127, 2008.

\bibitem{Gandy2011}
S.~Gandy, B.~Recht, and I.~Yamada, ``Tensor completion and low-n-rank tensor
  recovery via convex optimization,'' \emph{Inverse Problems}, vol.~27, no.~2,
  p. 025010, 2011.

\bibitem{Signoretto2011}
M.~Signoretto, R.~Van~de Plas, B.~De~Moor, and J.~A. Suykens, ``Tensor versus
  matrix completion: a comparison with application to spectral data,''
  \emph{IEEE Signal Processing Letters}, vol.~18, no.~7, pp. 403 -- 406, 2011.

\bibitem{Tomioka2011}
R.~Tomioka, T.~Suzuki, K.~Hayashi, and H.~Kashima, ``Statistical performance of
  convex tensor decomposition,'' in \emph{Advances in Neural Information
  Processing Systems}, 2011, pp. 972--980.

\bibitem{MHWG2014}
C.~Mu, B.~Huang, J.~Wright, and D.~Goldfarb, ``Square deal: Lower bounds and
  improved relaxations for tensor recovery,'' in \emph{ICML}, 2014, p. 2014.

\bibitem{mazumder2010tensor}
R.~Mazumder, T.~Hastie, and R.~Tibshirani, ``Spectral regularization algorithms
  for learning large incomplete matrices,'' \emph{The Journal of Machine
  Learning Research}, vol.~11, pp. 2287--2322, 2010.

\bibitem{liu2013}
J.~Liu, P.~Musialski, P.~Wonka, and J.~Ye, ``Tensor completion for estimating
  missing values in visual data,'' \emph{IEEE Trans. Pattern Analysis and
  Machine Intelligence}, vol.~35, no.~1, pp. 208--220, 2013.

\bibitem{fbcp}
Q.~Zhao, L.~Zhang, and A.~Cichocki, ``Bayesian cp factorization of incomplete
  tensors with automatic rank determination,'' \emph{IEEE Trans. Pattern
  Analysis and Machine Intelligence}, vol.~37, no.~9, pp. 1751--1763, 2015.

\bibitem{vershynin2018high}
R.~Vershynin, \emph{High-dimensional probability: An introduction with
  applications in data science}.\hskip 1em plus 0.5em minus 0.4em\relax
  Cambridge University Press, 2018, vol.~47.

\end{thebibliography}

\vskip -2\baselineskip plus -1fil

\begin{IEEEbiographynophoto}{Yihua Li}
Yihua Li graduated with B.S. and M.Eng degree in Electrical Engineering and Computer Science from MIT in 2015 and 2016, respectively. Her master thesis was supervised by Prof. Devavrat Shah and received David Adler Memorial EE MEng Thesis Award. Her interest lies in the theories and applications of statistical inference and deep learning, and she currently works as a software engineer at Google.
\end{IEEEbiographynophoto}

\vskip -2\baselineskip plus -1fil

\begin{IEEEbiographynophoto}{Devavrat Shah}
is a Professor with the department of Electrical Engineering and Computer Science at MIT. He is the director of Statistics and Data Science Center, Institute for Data, Systems and Society. He is a member of LIDS, CSAIL and ORC at MIT. He received his BTech in Computer Science from IIT Bombay in 1999 and PhD in Computer Science from Stanford University in 2004. 
His current research interest is in developing large-scale machine learning algorithms for unstructured data with particular interest in social data. He has made contributions to development of “gossip” protocols and “message-passing” algorithms for statistical inference which have been pillar of modern distributed data processing systems. 
Devavrat’s work has received broad recognition, including prize paper awards in Machine Learning, Operations Research and Computer Science, and career prizes including 2010 Erlang prize from the INFORMS Applied Probability Society, awarded bi-annually to a young researcher who has made outstanding contributions to applied probability. He is a distinguished young alumni of his alma mater IIT Bombay.
\end{IEEEbiographynophoto}

\vskip -2\baselineskip plus -1fil

\begin{IEEEbiographynophoto}{Dogyoon Song}
is currently pursuing his Ph.D. in the department of Electrical Engineering and Computer Science at MIT. He received his B.S. in Electrical Engineering, Mathematics, Physics and B.A. in Economics from Seoul National University in 2013. He is a recipient of the Samsung Scholarship and the Siebel Scholarship. His research interests are in theory and algorithms for optimization and statistical inference.
\end{IEEEbiographynophoto}

\vskip -2\baselineskip plus -1fil

\begin{IEEEbiographynophoto}{Christina Lee Yu} (formerly Christina E. Lee)
is an assistant professor at Cornell University in the Operations Research and Information Engineering Department. Prior to joining Cornell, she was a postdoc at Microsoft Research New England. She received her PhD and MS in Electrical Engineering and Computer Science from MIT in the Laboratory for Information and Decision Systems. She received her BS in Computer Science from California Institute of Technology. She is a recipient of the MIT Jacobs Presidential Fellowship, the NSF Graduate Research Fellowship, and the Claude E. Shannon Research Assistantship. Her research focuses on designing and analyzing scalable algorithms for processing social data based on principles from statistical inference. 
\end{IEEEbiographynophoto}


\vfill


\newpage

\appendices
\onecolumn
\section{Proof of Theorem \ref{thm:main_mse_user}}\label{sec:proof_main_thm}
\subsection{Outline of the Proof}\label{sec:proof_outline}

There are $m + n + 2mn$ independent sources of randomness in our model, $ \big\{ \featrow{u}\big\}_{u \in [m]}$, 
$\big\{ \featcol{i} \big\}_{i \in [n]}, \big\{ N(u,i) \big\}_{(u,i) \in [m] \times [n]}$,  $\big\{ M(u,i) \big\}_{(u,i) \in [m] \times [n]}$. 
We let $\Frandom$ denote the collection of these random variables. 

Due to the exchangeability of the model and the linearity of the expectation,
\[	\MSE(\hA) = \Exp{\big( \hA(u,i) - A(u,i) \big)^2 }	 = \Exp{\big( \hA(1,1) - A(1,1) \big)^2 }.	 \]
In Theorem \ref{thm:main_mse_user}, we provide an upper bound on $\Exp{\big( \hA(1,1) - A(1,1) \big)^2 ~\big|~ \featrow{1} }$, 
conditioned on the latent feature of the first row.

Recall from Algorithm \ref{alg:user_algorithm} that
\[	\hA(1,1) = 	 \frac{1}{| \baseest(1,1) |} \sum_{(v,1) \in \baseest(1,1)} Z(v,1) \] 
in our user-user fixed radius nearest neighbor algorithm when $\Nbase \geq 1$ and $\hA(1,1) = 0$ when $\Nbase = 0$. Therefore,
\begin{align*}
	\Exp{\big( \hA(1,1) - A(1,1) \big)^2 ~\big|~ \featrow{1}}
		&= \Exp{(A(1,1))^2 ~\Ind\big( { \Nbase = 0 }\big) ~\big|~ \featrow{1}} \\
			&\quad	+ \Exp{\big( \hA(1,1) - A(1,1) \big)^2 ~\Ind\big( { \Nbase \geq 1 }\big) ~\big|~ \featrow{1}}.
\end{align*}
By the assumption that the magnitude of the latent function $f$ is bounded by $\fbound$,
\begin{align*}
	\Exp{(A(1,1))^2 ~\Ind\big( { \Nbase = 0 }\big) ~\big|~ \featrow{1} }
		&\leq \fbound^2 \Prob{ \Nbase = 0 ~\big|~ \featrow{1}}
\end{align*}
and that
\begin{align}
	&\Exp{\big( \hA(1,1) - A(1,1) \big)^2 ~\Ind\big( { \Nbase \geq 1 }\big) ~\big|~ \featrow{1}}	\nonumber\\
		&\qquad=  \bbE\Bigg[ \bigg( \frac{1}{| \baseest(1,1) |} \sum_{(v,1) \in \baseest(1,1)} \big( A(v,1) - A(1,1) \big) \bigg)^2 
			~\Ind\big( { \Nbase \geq 1 }\big) ~\bigg|~ \featrow{1} \Bigg]		\nonumber\\
		&\qquad\quad + \bbE\Bigg[ \bigg( \frac{1}{| \baseest(1,1) |} \sum_{(v,1) \in \baseest(1,1)} N(v,1)\bigg)^2 
				~\Ind\big( { \Nbase \geq 1 }\big) ~\bigg|~ \featrow{1}\Bigg]		\nonumber\\
		&\qquad\stackrel{(a)}{\leq}
				\bbE_{\Frandom \setminus \featcol{1} } 
				\bigg[ \max_{(v,i) \in \baseest(1,1) } \Big\|  f(\featrow{v}, \cdot ) - f(\featrow{1}, \cdot) \Big\|_{L^2}^2 
					~\Ind\big( { \Nbase \geq 1 }\big)	~\Big|~ \featrow{1} \bigg]	\nonumber\\
			&\qquad\quad + C \sigma^2 \bbE \bigg[ \frac{1}{| \baseest(1,1) |} ~\Ind\big( { \Nbase \geq 1 }\big) ~\Big|~ \featrow{1} \bigg]	\label{eqn:mse_upper.0}
\end{align}
where $C > 0$ is an absolute constant (coming from $\| A \|_{\psi_2} \leftrightarrow \Var(A)$). 
The inequality in (a) follows from Lemmas \ref{lem:signal_MSE_upper} and \ref{lem:noise_MSE_upper}.

We need local properties of the probability measure on the latent space to upper bound the two terms in the last line; 
see Lemmas \ref{lem:signal_technical} and \ref{lem:noise_technical}.
We define the 
function $\localprob$ for $r \geq 0$ and $x \in \latsprow$
\begin{align*}
\localprob(x, r) = \bP_{\featrow{v} \sim \murow}\left(\big\| f(x, \cdot) - f(\featrow{v}, \cdot) \big\|_{L^2}^2 \leq r\right).
\end{align*} 
We prove in Lemma \ref{lem:prob_zero} that 
\begin{align*}
	\Prob{\Nbase = 0 ~\big|~ \featrow{1}}
		&\leq \exp\bigg( - \frac{(m-1)p}{8} \cdot \localprob\Big(\featrow{1}, \eta' - 2\sigma^2\Big) \bigg)\\
		&\quad	+ \exp \bigg[ -c \min \Big( \frac{t^2}{K^2}, \frac{t}{K} \Big) \frac{1}{2}(n-1)p^2 \bigg] 
			+ \exp\Big( - \frac{(n-1)p^2}{8} \Big).
\end{align*}

\subsection{Upper Bounding the Contribution of Signal on MSE}
\begin{lemma}\label{lem:signal_MSE_upper}
	\begin{align*}
		&\bbE\Bigg[ \bigg( \frac{1}{| \baseest(1,1) |} \sum_{(v,1) \in \baseest(1,1)} \big( A(v,1) - A(1,1) \big) \bigg)^2
		~\Ind\big( { \Nbase \geq 1 }\big) ~\bigg|~ \featrow{1} \Bigg]\\
			&\qquad= \bbE
				\bigg[ \max_{(v,i) \in \baseest(1,1) } \Big\|  f(\featrow{v}, \cdot ) - f(\featrow{1}, \cdot) \Big\|_{L^2}^2 
				~\Ind\big( { \Nbase \geq 1 }\big)	~\bigg|~ \featrow{1}\bigg].
	\end{align*}
\end{lemma}
\begin{proof}
Recall from Section \ref{sec:proof_outline} that we let $\Frandom$ denote the collection of $m + n + 2mn$ independent sources of randomness 
in our model, $ \big\{ \featrow{u}\big\}_{u \in [m]}$, $\big\{ \featcol{i} \big\}_{i \in [n]}, \big\{ N(u,i) \big\}_{(u,i) \in [m] \times [n]}$,  
$\big\{ M(u,i) \big\}_{(u,i) \in [m] \times [n]}$. 

By the tower property of expectation,
\begin{align*}
	&\bbE_{\Frandom | \featrow{1}}\Bigg[ \bigg( \frac{1}{\Nbase} \sum_{(v,1) \in \baseest(1,1)} \big( A(v,1) - A(1,1) \big) \bigg)^2	
	~\Ind\big( { \Nbase \geq 1 }\big) ~\bigg|~ \featrow{1} \Bigg]\\
		&= \bbE_{\Nbase}\Bigg[ ~ \bbE_{\Frandom \big| \featrow{1}, \Nbase }\bigg[ \bigg( \frac{1}{\Nbase} \sum_{(v,1) \in \baseest(1,1)} 
		\big( A(v,1) - A(1,1) \big) \bigg)^2 ~\Ind\big( { \Nbase \geq 1 }\big) ~\bigg|~ \featrow{1}, \Nbase \bigg] ~\Bigg].
\end{align*}
We investigate the conditional expectation, again using the tower property and the fact that $\Nbase$ is fully determined when we condition on $\Frandom \setminus \featcol{1}$. 
For the sake of readability, $\Ind\big( { \Nbase \geq 1 }\big)$ is omitted in the lines below.
\begin{align*}
	&\bbE_{\Frandom \big| \featrow{1}, \Nbase}\bigg[ \bigg( \frac{1}{| \baseest(1,1) |} \sum_{(v,1) \in \baseest(1,1)} \big( A(v,1) - A(1,1) \big) \bigg)^2 
		~\bigg|~ \featrow{1}, \Nbase \bigg] \\
		&	= \frac{1}{| \baseest(1,1) |^2} \bbE_{\Frandom \big| \featrow{1}, \Nbase }\bigg[ \bigg(  \sum_{(v,1) \in \baseest(1,1)} 
			\big( A(v,1) - A(1,1) \big) \bigg)^2 ~\bigg|~ \featrow{1}, \Nbase \bigg] \\
		&	\leq \frac{1}{| \baseest(1,1) |^2} \bbE_{\Frandom \big| \featrow{1}, \Nbase }\Bigg[ \Bigg(   \sum_{(v,1) \in \baseest(1,1) \atop (v',1) \in \baseest(1,1)} 
			\big| A(v,1) - A(1,1) \big|\big| A(v',1) - A(1,1) \big| \Bigg) ~\bigg|~ \featrow{1}, \Nbase \Bigg] \\
		&	\stackrel{(a)}{=} \frac{1}{| \baseest(1,1) |^2 } \bbE_{\Frandom \setminus \featcol{1} \big| \featrow{1}, \Nbase } \Bigg[~ \bbE_{\featcol{1} | \featrow{1}, \Nbase}
			\bigg[\\
			&\qquad \sum_{(v,1) \in \baseest(1,1) \atop (v',1) \in \baseest(1,1)}
			\big| A(v,1) - A(1,1) \big|\big| A(v',1) - A(1,1) \big| ~ \Big| ~ \featrow{1}, \Nbase \bigg] ~\bigg|~ \featrow{1}, \Nbase \Bigg]\\
		&	\stackrel{(b)}{=} \frac{1}{| \baseest(1,1) |^2 } \bbE_{\Frandom \setminus \featcol{1}  \big| \featrow{1}, \Nbase } \Bigg[ \\
			&\qquad\sum_{(v,1) \in \baseest(1,1) \atop (v',1) \in \baseest(1,1)}
				\bbE_{\featcol{1}  \big| \featrow{1} } \Big[ \big| A(v,1) - A(1,1) \big|\big| A(v',1) - A(1,1) \big| ~\Big|~ \featrow{1}\Big]
			~\bigg|~ \featrow{1}, \Nbase \Bigg]\\
		&	= \frac{1}{| \baseest(1,1) |^2 } \bbE_{\Frandom \setminus \featcol{1}  \big| \featrow{1}, \Nbase } \Bigg[\\
			&\qquad \sum_{(v,1) \in \baseest(1,1) \atop (v',1) \in \baseest(1,1)}
				\Big\| \big[ f(\featrow{v}, \cdot ) - f(\featrow{1}, \cdot) \big] \big[ f(\featrow{v'}, \cdot ) - f(\featrow{1}, \cdot) \big] \Big\|_{L^1}
			~\bigg|~ \featrow{1}, \Nbase \Bigg]\\
		&	\stackrel{(c)}{\leq} \frac{1}{| \baseest(1,1) |^2 } \bbE_{\Frandom \setminus \featcol{1} \big| \featrow{1}, \Nbase } \Bigg[ \\
			&\qquad\sum_{(v,1) \in \baseest(1,1) \atop (v',1) \in \baseest(1,1)}
				\Big\|  f(\featrow{v}, \cdot ) - f(\featrow{1}, \cdot) \Big\|_{L^2} \Big\| f(\featrow{v'}, \cdot ) - f(\featrow{1}, \cdot) \Big\|_{L^2}
			~\bigg|~ \featrow{1}, \Nbase \Bigg]\\
		&	\leq \frac{1}{| \baseest(1,1) |^2 } \bbE_{\Frandom \setminus \featcol{1} \big| \featrow{1}, \Nbase } 
			\Bigg[ \Nbase^2  \max_{(v,i) \in \baseest(1,1) } \Big\|  f(\featrow{v}, \cdot ) - f(\featrow{1}, \cdot) \Big\|_{L^2}^2
			~\bigg|~ \featrow{1}, \Nbase \Bigg]\\
		&	=  \bbE_{\Frandom \setminus \featcol{1} \big| \featrow{1}, \Nbase } 
			\Bigg[ \max_{(v,1) \in \baseest(1,1) } \Big\|  f(\featrow{v}, \cdot ) - f(\featrow{1}, \cdot) \Big\|_{L^2}^2
			~\bigg|~ \featrow{1},  \Nbase \Bigg].
\end{align*}
Equation (a) follows from the conditional independence; (b) follows from Fubini's theorem and the independence between 
$\fcol{1}$ and $\Nbase$; and (c) follows from Cauchy-Schwarz inequality.

Consequently, it follows that 
\begin{align*}
	&\bbE\Bigg[ \bigg( \frac{1}{\Nbase} \sum_{(v,1) \in \baseest(1,1)} \big( A(v,1) - A(1,1) \big) \bigg)^2 \Ind\big( { \Nbase \geq 1 }\big) ~\bigg|~ \featrow{1}\Bigg]\\
		&\qquad= \bbE
			\Bigg[ \max_{(v,i) \in \baseest(1,1) } \Big\|  f(\featrow{v}, \cdot ) - f(\featrow{1}, \cdot) \Big\|_{L^2}^2 \Ind\big( { \Nbase \geq 1 }\big) ~\bigg|~ \featrow{1} \Bigg].
\end{align*}
\end{proof}

\begin{lemma}\label{lem:signal_technical}
Let $\eta \geq 2\sigma^2$. The following inequality holds for user-user fixed radius neighbor algorithm:
	\begin{align*}
		&\bbE
			\bigg[ \max_{(v,1) \in \baseest(1,1) } \Big\|  f(\featrow{v}, \cdot ) - f(\featrow{1}, \cdot) \Big\|_{L^2}^2
			\Ind\big( { \Nbase \geq 1 }\big) ~\bigg|~ \featrow{1}	\bigg]\\
			&\qquad\leq \big( \eta - 2\sigma^2 \big)
			+ \frac{K \sqrt{\pi}}{\sqrt{2 c (n-1)p^2}} + \frac{2K e^{-\frac{c}{2} (n-1)p^2}}{c (n-1)p^2}
			+ \bigg( \frac{K \sqrt{\pi}}{2\sqrt{c }} + \frac{K e^{-c}}{c} \bigg) \exp\Big( - \frac{(n-1)p^2}{8} \Big)\\
			&\qquad\leq \big( \eta - 2\sigma^2 \big)
			+ CK \bigg[  \frac{1}{\sqrt{ (n-1) p^2}} + \exp \big(-c' (n-1)p^2\big) \bigg].
	\end{align*}
In the above expression, $K \triangleq \Big(  \frac{\fbound}{\sqrt{\ln 2}} + 2 \sigma \Big)^2$ and $C, c, c' > 0$ are absolute constants. 
\end{lemma}
\begin{proof}
\textbf{Step 1:}
Choose $v \in [m] \setminus \{1\}$ such that $(v,1) \in \baseest(1,1)$. Observe that 
\[
	\dissimrow(1,v) = \frac{1}{|\basesim(1) \cap \basesim(v)|} \sum_{j \in \basesim(1) \cap \basesim(v)} \big( Z(1,j) - Z(v,j) \big)^2.
\] 
Note that $Z(1,j) - Z(v,j) = \big( A(1,j) -A(v,j) \big) + \big( N(1,j) - N(v,j) \big)$. By the model assumptions, we have
\begin{align*}
	&|A(1,j) - A(v,j)| \leq 2\fbound \quad\forall v, j	\qquad \Rightarrow \qquad \big\| A(1,j) - A(v,j) \big\|_{\psi_2} 
		\leq \frac{2\fbound}{\sqrt{\ln 2}},\\
	&\big\| N(1,j) - N(v,j) \big\|_{\psi_2} \leq \big\| N(1,j) \big\|_{\psi_2} + \big\| N(v,j) \big\|_{\psi_2} = 2\sigma.
\end{align*}
Therefore (cf. \cite[Lemma 2.7.6]{vershynin2018high}),
\begin{align*}
	\big\| Z(1,j) - Z(v,j) \big\|_{\psi_2} &\leq \big\| A(1,j) - A(v,j) \big\|_{\psi_2} + \big\| N(1,j) - N(v,j) \big\|_{\psi_2} 
		\leq  \frac{2\fbound}{\sqrt{\ln 2}} + 2 \sigma,\\
	\big\| \big(Z(1,j) - Z(v,j) \big)^2 \big\|_{\psi_1} &= \big\| Z(1,j) - Z(v,j) \big\|_{\psi_2}^2
		\leq \Big(  \frac{2\fbound}{\sqrt{\ln 2}} + 2 \sigma \Big)^2.
\end{align*}

For notational conciseness, we let $f_{1v}^2$ denote the conditional expectation of $\dissimrow(1,v)$:
\begin{equation}\label{eqn:exp_dissim}
	f_{1v}^2 \triangleq \Exp{\dissimrow(1,v) | \featrow{1},\featrow{v}} 
		= \big\|  f(\featrow{v}, \cdot ) - f(\featrow{1}, \cdot) \big\|_{L^2}^2 + 2\sigma^2.
\end{equation}
Now, suppose that $|\basesim(1) \cap \basesim(v)| = N_v$. We observe that 
$\Exp{\dissimrow(1,v) | \featrow{1},\featrow{v}, |\basesim(1) \cap \basesim(v)| = N_v} = f_{1v}^2$ regardless of $N_v$. 
Then it follows from Bernstein's inequality (cf. \cite[Theorem 2.8.2]{vershynin2018high}) that for any $t \geq 0$,
\begin{equation}\label{eqn:conc_dissim}
	\Prob{ \dissimrow(1,v) - f_{1v}^2 \leq - t ~ \Big|~ \featrow{1}, \featrow{v}, |\basesim(1) \cap \basesim(v)| = N_v} 
		\leq \exp \bigg[ -c \min \Big( \frac{t^2}{K^2}, \frac{t}{K} \Big) N_v \bigg].
\end{equation}
Note that we provided an upper bound on the conditional probability in \eqref{eqn:conc_dissim} that holds for 
any realization of $\featrow{1}$ and $\featrow{v}$, and thus we could remove the conditioning on $\featrow{v}$.
By plugging in the value of $f_{1v}^2$, this leads to the following inequality: for any $t \geq 0$,
\begin{align*}
	&\Prob{  \big\|  f(\featrow{v}, \cdot ) - f(\featrow{1}, \cdot) \big\|_{L^2}^2 \geq \dissimrow(1,v) - 2\sigma^2 + t  
		~ \Big|~ \featrow{1}, |\basesim(1) \cap \basesim(v)| = N_v}\\
		&\leq  \exp \bigg[ -c \min \Big( \frac{t^2}{K^2}, \frac{t}{K} \Big) N_v \bigg].	
\end{align*}


\textbf{Step 2:} Let $v^* \in [m] \setminus \{1\}$ denote the maximizer such that $(v^*, 1) \in \baseest(1,1)$ and 
\[	\Big\|  f(\featrow{v^*}, \cdot ) - f(\featrow{1}, \cdot) \Big\|_{L^2}^2
		= \max_{(v,1) \in \baseest(1,1) } \Big\|  f(\featrow{v}, \cdot ) - f(\featrow{1}, \cdot) \Big\|_{L^2}^2.	\]
According to the description of our user-user fixed radius nearest neighbor algorithm, $\dissimrow(1,v^*) \leq \eta$.
\begin{align}
	&\bbE
			\bigg[ \max_{(v,i) \in \baseest(1,1) } \big\|  f(\featrow{v}, \cdot ) - f(\featrow{1}, \cdot) \big\|_{L^2}^2 
			~\Big|~ \featrow{1}	\bigg]	\nonumber\\
		&\qquad= \bbE
			\bigg[ \big\|  f(\featrow{v^*}, \cdot ) - f(\featrow{1}, \cdot) \big\|_{L^2}^2 ~\Big|~ \featrow{1}	\bigg]	\nonumber\\
		&\qquad= \bbE_{N_v }  \Bigg[ \int_0^{\infty}
			\bbP\bigg( \big\|  f(\featrow{v^*}, \cdot ) - f(\featrow{1}, \cdot) \big\|_{L^2}^2	 \geq s 
			~ \Big| ~ \featrow{1},  |\basesim(1) \cap \basesim(v)| = N_v \bigg) ~ds \Bigg].	\label{eqn:upper_signal_term.1}
\end{align}

Now we observe that 
\begin{align}
	& \int_0^{\infty} \bbP\bigg( \big\|  f(\featrow{v^*}, \cdot ) - f(\featrow{1}, \cdot) \big\|_{L^2}^2	 \geq s 
			~ \Big| ~ \featrow{1}, |\basesim(1) \cap \basesim(v)| = N_v \bigg) ~ds	\nonumber\\
		&\qquad\leq \int_0^{\eta - 2\sigma^2} ds
			+ \int_{0}^{K} \exp \bigg[ -c N_v \frac{t^2}{K^2} \bigg] ds
			+ \int_{K}^{\infty} \exp \bigg[ -c N_v \frac{t}{K} \bigg] ds	\nonumber\\
		&\qquad\leq \big( \eta - 2\sigma^2 \big) 
			+ \int_{0}^{\infty} \exp \bigg[ -c N_v \frac{t^2}{K^2} \bigg] ds
			+ \int_{K}^{\infty} \exp \bigg[ -c N_v \frac{t}{K} \bigg] ds	\nonumber\\
		&\qquad= \big( \eta - 2\sigma^2 \big)
			+ \frac{K \sqrt{\pi}}{2\sqrt{c N_v}} + \frac{K e^{-c N_v}}{cN_v}.	\label{eqn:upper_signal_term.2}
\end{align}
Recall that for $a > 0$ and $b \in \Reals$,
\[	
	\int_0^{\infty} e^{-a x^2} dx = \frac{1}{2}\sqrt{\frac{\pi}{a}}
	\qquad\text{and}\qquad
	\int_b^{\infty} e^{-a x} dx = \frac{1}{a}e^{-ab}.	
\]

\textbf{Step 3:} Observe that $ |\basesim(1) \cap \basesim(v)| \sim \textrm{Binomial}(n-1, p^2)$. By the binomial Chernoff theorem,
\begin{equation}\label{eqn:overlap_lower}
	\Prob{ |\basesim(1) \cap \basesim(v)| \leq \frac{1}{2} (n-1)p^2} \leq \exp\Big( - \frac{(n-1)p^2}{8} \Big).
\end{equation}
Therefore, with the shorthand notation $N_v = |\basesim(1) \cap \basesim(v)|$, we can see that
\begin{align*}
	&\bbE \Bigg[ \big( \eta - 2\sigma^2 \big)	+ \frac{K \sqrt{\pi}}{2\sqrt{c N_v}} + \frac{K e^{-c N_v}}{cN_v} \Bigg]\\
		&\qquad= \big( \eta - 2\sigma^2 \big)
			+ \bbE \bigg[ \frac{K \sqrt{\pi}}{2\sqrt{c N_v}} + \frac{K e^{-c N_v}}{cN_v} ~\Big|~N_v > \frac{1}{2} (n-1)p^2 \bigg]
			\Prob{N_v > \frac{1}{2} (n-1)p^2}\\
			&\qquad\quad
			+ \bbE \bigg[ \frac{K \sqrt{\pi}}{2\sqrt{c N_v}} + \frac{K e^{-c N_v}}{cN_v} ~\Big|~N_v \leq \frac{1}{2} (n-1)p^2 \bigg]
			\Prob{N_v \leq \frac{1}{2} (n-1)p^2}\\
		&\qquad\stackrel{(a)}{\leq} \big( \eta - 2\sigma^2 \big)
			+ \bbE \bigg[ \frac{K \sqrt{\pi}}{2\sqrt{c N_v}} + \frac{K e^{-c N_v}}{cN_v} ~\Big|~N_v > \frac{1}{2} (n-1)p^2 \bigg]
			+ \bigg( \frac{K \sqrt{\pi}}{2\sqrt{c }} + \frac{K e^{-c}}{c} \bigg)
			\Prob{N_v \leq \frac{1}{2} (n-1)p^2}\\
		&\qquad\leq \big( \eta - 2\sigma^2 \big)
			+ \frac{K \sqrt{\pi}}{\sqrt{2 c (n-1)p^2}} + \frac{2K e^{-\frac{c}{2} (n-1)p^2}}{c (n-1)p^2}
			+ \bigg( \frac{K \sqrt{\pi}}{2\sqrt{c }} + \frac{K e^{-c}}{c} \bigg) \exp\Big( - \frac{(n-1)p^2}{8} \Big).
\end{align*}
The inequality in (a) follows from that $N_v \geq 1$ and $\frac{K e^{-c N_v}}{cN_v}$ is maximized when $N_v = 1$.
\end{proof}


\subsection{Upper Bounding the Contribution of Noise on MSE}
\begin{lemma}\label{lem:noise_MSE_upper}
	\begin{align*}
		 \bbE\Bigg[ \bigg( \frac{1}{| \baseest(1,1) |} \sum_{(v,1) \in \baseest(1,1)} N(v,1)\bigg)^2~\Ind\big( { \Nbase \geq 1 }\big) 
		 	~\bigg|~ \featrow{1}\Bigg]
		 	&\leq C\sigma^2 \bbE \bigg[ \frac{1}{| \baseest(1,1) |} ~\Ind\big( { \Nbase \geq 1 }\big) ~\Big|~ \featrow{1} \bigg].
	\end{align*}
\end{lemma}
\begin{proof}
By the tower property of expectation,
\begin{align*}
	& \bbE\Bigg[ \bigg( \frac{1}{| \baseest(1,1) |} \sum_{(v,1) \in \baseest(1,1)} N(v,1)\bigg)^2~\Ind\big( { \Nbase \geq 1 }\big) 
		~\bigg|~ \featrow{1} \Bigg]\\
	 	&\qquad = \bbE\Bigg[ ~ \bbE\bigg[ \bigg( \frac{1}{| \baseest(1,1) |} \sum_{(v,1) \in \baseest(1,1)} N(v,1)\bigg)^2
			~\Ind\big( { \Nbase \geq 1 }\big) ~ \bigg|~ \featrow{1}, \Nbase \bigg] ~\bigg|~ \featrow{1}~\Bigg]\\
		&\qquad \leq \bbE\Bigg[ ~ \bbE\bigg[  \frac{1}{| \baseest(1,1) |^2} \sum_{(v,1) \in \baseest(1,1)} | N(v,1) |^2
			~\Ind\big( { \Nbase \geq 1 }\big) ~ \bigg|~ \featrow{1}, \Nbase \bigg] ~\bigg|~ \featrow{1}~\Bigg]\\
		&\qquad \leq \bbE\Bigg[ ~  \frac{1}{| \baseest(1,1) |^2}~ \sum_{(v,1) \in \baseest(1,1)}  \bbE\Big[ | N(v,1) |^2 
			~\big|~  \featrow{1}, \Nbase \Big] ~ \Ind\big(\Nbase \geq 1\big) ~\bigg|~ \featrow{1}~\Bigg]	\\
		&\qquad \stackrel{(a)}{\leq} C\sigma^2 \bbE \bigg[ \frac{1}{| \baseest(1,1) |} ~ \Ind\big(\Nbase \geq 1\big) 
			~\bigg|~ \featrow{1} \bigg].
\end{align*}
(a) follows from $ \bbE\Big[ | N(v,1) |^2 ~\big|~  \featrow{1}, \Nbase \Big] =  \bbE\big[ | N(v,1) |^2 \big] \leq C \sigma^2$, 
as a result of $\| N(v,1) \|_{\psi_2} \leq \sigma$.
\end{proof}

\begin{lemma}\label{lem:noise_technical}
	Let $\eta' \geq 2\sigma^2$ and $\eta \geq \eta' + K \max\Big( \sqrt{\frac{4 \log (m-1)}{c(n-1)p^2}}, \frac{4 \log (m-1)}{c(n-1)p^2} \Big)$. 
	The following inequality holds for user-user fixed radius nearest neighbor algorithm:
	\begin{align*}
		\bbE \bigg[ \frac{1}{| \baseest(1,1) |} ~ \Ind\big(\Nbase \geq 1\big)  ~ \Big| ~  \featrow{1} \bigg]
			&\leq 2 \bigg[  (m-1)p \cdot \localprob\Big(\featrow{1},  \eta' - 2\sigma^2\Big)  \bigg]^{-1}	+ \frac{1}{m-1} \\
				&\qquad+ \exp\bigg( - \frac{(m-1)p}{8} \cdot \localprob\Big(\featrow{1},  \eta' - 2\sigma^2\Big) \bigg)
				+ (m-1)\exp \Big( - \frac{(n-1)p^2}{8} \Big).
%
%
%
	\end{align*}
where 
$\localprob(x, r) = \bP_{\featrow{v} \sim \murow}\left(\big\| f(x, \cdot) - f(\featrow{v}, \cdot) \big\|_{L^2}^2 \leq r\right)$
for any $x \in \latsprow$ and $r > 0$.
\end{lemma}


\begin{proof}
\textbf{Step 1:}
%
%
%
%
Our interest is in bounding conditional expectation of $\frac{1}{| \baseest(1,1) |} ~ \Ind\big(\Nbase \geq 1\big)$ given 
$\featrow{1} \in \latsprow$. To that end, let us fix $\featrow{1}$.  
Let $\eta' = \eta - K \max\Big( \sqrt{\frac{4 \log (m-1)}{c(n-1)p^2}}, \frac{4 \log (m-1)}{c(n-1)p^2} \Big)$.
For each $v \in [m] \setminus \{ 1 \}$, the following two conditions
\begin{equation}\label{eqn:cond_good_row}
	M(v,1) = 1	\qquad\text{and}\qquad   \big\| f(\featrow{1}, \cdot) - f(\featrow{v}, \cdot) \big\|_{L^2}^2 \leq \eta' - 2\sigma^2
\end{equation}
are satisfied with success probability\footnote{Note that the first condition solely depends on $M(v,1)$, 
while the second condition depends only on $\featrow{1}, \featrow{v}$, hence, they are independent events.} 
$p \cdot \localprob\big(\featrow{1}, \eta' - 2\sigma^2 \big)$.

Let $\goodrows \subset [m] \setminus \{1 \}$ denote the set of row indices 
that satisfy the two conditions described in \eqref{eqn:cond_good_row}. By the binomial Chernoff bound, 
\begin{equation}\label{eqn:upper_noise_term.1}
	\Prob{ |\goodrows| \leq \frac{1}{2} (m-1)p  \cdot \localprob\Big(\featrow{1}, \eta' - 2\sigma^2\Big) ~ \bigg|~ \featrow{1} }
		\leq \exp\bigg( - \frac{(m-1)p}{8} \cdot \localprob\Big(\featrow{1}, \eta' - 2\sigma^2\Big) \bigg).
\end{equation}

\textbf{Step 2:}
Next, we want to show that $\{(v,1): v \in \goodrows\} \subset \baseest(1,1)$ with high probability.
For that purpose, we first require $|\basesim(1) \cap \basesim(v) |$ to be sufficiently large for all $v \in \goodrows$. Observe that 
for each $v \in \goodrows$, $|\basesim(1) \cap \basesim(v) |$ is distributed following $\textrm{Binomial}(n-1, p^2)$. Again by the binomial 
Chernoff bound (as in \eqref{eqn:overlap_lower}), we have 
\[	\Prob{|\basesim(1) \cap \basesim(v)| \leq \frac{(n-1)p^2}{2}} \leq \exp \Big( - \frac{(n-1)p^2}{8} \Big).	\]
Since $|\goodrows| \leq m-1$, it follows from the union bound that 
\begin{equation}\label{eqn:upper_noise_term.2}
	\Prob{ \min_{v \in \goodrows }|\basesim(1) \cap \basesim(v)| \leq \frac{(n-1)p^2}{2}} \leq (m-1)\exp \Big( - \frac{(n-1)p^2}{8} \Big).
\end{equation}

By construction, for all $v \in \goodrows$,
\begin{align*}
	\big\| f(\featrow{1}, \cdot) - f(\featrow{v}, \cdot) \big\|_{L^2}^2 \leq \eta' - 2\sigma^2.
\end{align*}
Note that we can obtain the following concentration inequality by similar arguments\footnote{This provides a probabilistic 
tail bound on the opposite side of that in \eqref{eqn:conc_dissim}. The proof remains the same.} as in \eqref{eqn:conc_dissim}:
\begin{equation}\label{eqn:conc_dissim2}
	\Prob{ \dissimrow(1,v) - f_{1v}^2 \geq t ~ \Big|~ \featrow{1}, \featrow{v}, |\basesim(1) \cap \basesim(v)| = N_v} 
		\leq \exp \bigg[ -c \min \Big( \frac{t^2}{K^2}, \frac{t}{K} \Big) N_v \bigg].
\end{equation}
With another application of the union bound, this implies that 
\begin{align*}
	&\Prob{ \max_{v \in \goodrows }\dissimrow(1,v) \geq \eta' + t ~\bigg| ~ \min_{v \in \goodrows }|\basesim(1) \cap \basesim(v)| > \frac{(n-1)p^2}{2}}\\ 
		&\qquad\leq (m-1)  \exp \bigg[ -c \min \Big( \frac{t^2}{K^2}, \frac{t}{K} \Big) \frac{(n-1)p^2}{2} \bigg].
\end{align*}
Therefore\footnote{Recall that $\eta' = \eta - t$ with $t = K \max\Big( \sqrt{\frac{4 \log (m-1)}{c(n-1)p^2}}, \frac{4 \log (m-1)}{c(n-1)p^2} \Big)$.},
\begin{equation}\label{eqn:upper_noise_term.3}
	\Prob{ \max_{v \in \goodrows }\dissimrow(1,v) \geq \eta ~\bigg| ~ \min_{v \in \goodrows }|\basesim(1) \cap \basesim(v)| > \frac{(n-1)p^2}{2}}
		\leq \frac{1}{m-1}.
\end{equation}
Observe that 
\begin{equation}\label{eqn:upper_noise_term.4}
	\text{if}\quad
	\max_{v \in \goodrows }\dissimrow(1,v) \leq \eta
	\qquad\text{then}\quad
	\{(v,1): v \in \goodrows\} \subset \baseest(1,1),
\end{equation}
which implies that $|\baseest(1,1)| \geq |\goodrows|$.

\textbf{Step 3:}
Let 
\begin{align*}
	\Eva	&\triangleq \bigg\{  |\goodrows| > \frac{1}{2} (m-1)p  \cdot \localprob\Big(\featrow{1}, \eta' - 2\sigma^2\Big) \bigg\}\\
	\Evb	&\triangleq \bigg\{  \max_{v \in \goodrows }\dissimrow(1,v) \leq \eta  \bigg\},\\
	\Evc	&\triangleq \bigg\{  \min_{v \in \goodrows }|\basesim(1) \cap \basesim(v)| > \frac{(n-1)p^2}{2} \bigg\}.
\end{align*}
By the law of total probability,
\begin{align*}
	~& \bbE \bigg[ \frac{1}{| \baseest(1,1) |} ~ \Ind\big(\Nbase \geq 1\big) ~\Big|~ \featrow{1} \bigg]  \\
	~	& \qquad \stackrel{(a)}{\leq} \bbE \bigg[ \frac{1}{| \baseest(1,1) |} ~ \Ind\big(\Nbase \geq 1\big) ~\Big|~\featrow{1},  \Eva \cap \Evb  \bigg] \Prob{\Eva \cap \Evb} 
			+ \Prob{\Eva^c \cup \Evb^c }\\
	~	&\qquad \leq \bbE \bigg[ \frac{1}{| \baseest(1,1) |} ~ \Ind\big(\Nbase \geq 1\big) ~\Big|~ \featrow{1}, \Eva \cap \Evb  \bigg] 
			+ \Prob{\Eva^c} + \Prob{\Evb^c|\Evc} + \Prob{\Evc^c}\\
	~	&\qquad  \stackrel{(b)}{\leq}
			2 \bigg[  (m-1)p  \cdot \localprob\Big(\featrow{1}, \eta' - 2\sigma^2\Big) \bigg]^{-1}	+ \frac{1}{m-1} \\
				&\qquad+ \exp\bigg( - \frac{(m-1)p}{8} \cdot \localprob\Big(\featrow{1}, \eta' - 2\sigma^2\Big) \bigg)
				+ (m-1)\exp \Big( - \frac{(n-1)p^2}{8} \Big)
\end{align*}
where (a) follows from $ \frac{1}{| \baseest(1,1) |} ~ \Ind\big(\Nbase \geq 1\big) \leq 1$; and (b) follows from that \eqref{eqn:upper_noise_term.4} and
\begin{align*}
	\Prob{\Eva^c}	&\leq \exp\bigg( - \frac{(m-1)p}{8} \cdot \localprob\Big(\featrow{1}, \eta' - 2\sigma^2\Big) \bigg),
			&\because \eqref{eqn:upper_noise_term.1}\\
	\Prob{\Evb^c|\Evc} &\leq	\frac{1}{m-1},		
			&\because \eqref{eqn:upper_noise_term.3}\\
	\Prob{\Evc^c}	&\leq (m-1)\exp \Big( - \frac{(n-1)p^2}{8} \Big).		
			&\because\eqref{eqn:upper_noise_term.2}
\end{align*}

\end{proof}

\subsection{Upper Bounding the Probability of $\Nbase = 0$}

\begin{lemma}\label{lem:prob_zero}
Let $\eta' \geq 2\sigma^2$ and $\eta \geq \eta' + K \max\Big( \sqrt{\frac{4 \log (m-1)}{c(n-1)p^2}}, \frac{4 \log (m-1)}{c(n-1)p^2} \Big)$. 
The following inequality holds for user-user fixed radius nearest neighbor algorithm:
\begin{align*}
	\Prob{\Nbase = 0 ~\big|~ \featrow{1}}
		&\leq \exp\bigg( - \frac{(m-1)p}{8} \cdot \localprob\Big(\featrow{1}, \eta' - 2\sigma^2\Big) \bigg)
			+ \frac{1}{(m-1)^2} + \exp\Big( - \frac{(n-1)p^2}{8} \Big).
\end{align*}
\end{lemma}
\begin{proof}
Recall the definition of $\goodrows$ from \eqref{eqn:cond_good_row}:
\[	\goodrows = \left\{ v \in [m] \setminus \{1 \} ~\text{ such that }~
		M(v,1) = 1	~\text{ and }~  \big\| f(\featrow{1}, \cdot) - f(\featrow{v}, \cdot) \big\|_{L^2}^2 \leq \eta' - 2\sigma^2\right\}.	\]
Note that we obtain in \eqref{eqn:upper_noise_term.1} that
\[
	\Prob{ |\goodrows| \leq \frac{1}{2} (m-1)p  \cdot \localprob\Big(\featrow{1}, \eta' - 2\sigma^2\Big) ~\bigg|~ \featrow{1} }
		\leq \exp\bigg( - \frac{(m-1)p}{8} \cdot \localprob\Big(\featrow{1}, \eta' - 2\sigma^2\Big) \bigg)
\]
by applying the binomial Chernoff bound.

Now we observe that if there exists at least one $v \in \goodrows$, for which 
\[	\dissimrow(1,v) - \big\|  f(\featrow{v}, \cdot ) - f(\featrow{1}, \cdot) \big\|_{L^2}^2 - 2\sigma^2 \leq \eta - \eta', 	\]
then $\Nbase \geq 1$ due to the Lipschitzness assumption on $f$. Therefore, we have 
\begin{align}
	\Prob{\Nbase = 0}
		&\leq \Prob{ |\goodrows| = 0 } \nonumber\\
		&\quad + \Prob{ \dissimrow(1,v) - \big\|  f(\featrow{v}, \cdot ) - f(\featrow{1}, \cdot) \big\|_{L^2}^2 - 2\sigma^2 \leq \eta - \eta', ~~ \forall v \in \goodrows} \Prob{ |\goodrows| \geq 1 }\nonumber\\
		&\leq \Prob{ |\goodrows| = 0 } + \Prob{ \dissimrow(1,v_0) - \big\|  f(\featrow{v_0}, \cdot ) - f(\featrow{1}, \cdot) \big\|_{L^2}^2 - 2\sigma^2 \leq \eta - \eta' } \label{eq:prob_0_term.1}
\end{align}
where $v_0 = \min \{ v \in \goodrows \}$.

We bound the two terms in \eqref{eq:prob_0_term.1} separately. First, as long as $\localprob\Big(\featrow{1}, \eta' - 2\sigma^2\Big) > 0$, it follows from \eqref{eqn:upper_noise_term.1} that
\begin{align*}
	\Prob{ |\goodrows| = 0 ~\big|~ \featrow{1} }
		&\leq \Prob{ |\goodrows| \leq \frac{1}{2} (m-1)p  \cdot \localprob\Big(\featrow{1}, \eta' - 2\sigma^2\Big) ~\bigg|~ \featrow{1}}\\
		&\leq \exp\bigg( - \frac{(m-1)p}{8} \cdot \localprob\Big(\featrow{1}, \eta' - 2\sigma^2\Big) \bigg).
\end{align*}
Second, by the usual trick of total probability (recall from \eqref{eqn:exp_dissim} that $f_{1v_0}^2$ denotes 
$\Exp{\dissimrow(1,v_0) | \featrow{1}, \featrow{v_0}}$),
\begin{align*}
	&\Prob{ \dissimrow(1,v_0) - f_{1v_0}^2 \leq \eta - \eta' ~\big|~ \featrow{1}}\\
		&\qquad\leq \Prob{ \dissimrow(1,v_0) - f_{1v_0}^2 \leq \eta - \eta' ~\Big|~ \featrow{1}, |\basesim(1) \cap \basesim(v_0)| \geq \frac{1}{2}(n-1)p^2 }\\
		&\qquad\quad
			+ \Prob{  |\basesim(1) \cap \basesim(v_0)| < \frac{1}{2}(n-1)p^2 ~\bigg|~ \featrow{1}}\\
	&\qquad\stackrel{(a)}{\leq} \exp \bigg[ -c \min \Big( \frac{(\eta - \eta')^2}{K^2}, \frac{\eta - \eta'}{K} \Big) \frac{1}{2}(n-1)p^2 \bigg] 
		+ \Prob{  |\basesim(1) \cap \basesim(v_0)| < \frac{1}{2}(n-1)p^2 }\\
	&\qquad\stackrel{(b)}{\leq} \exp \big[ -2 \log(m-1) \big] 
		+ \exp\Big( - \frac{(n-1)p^2}{8} \Big)	
\end{align*}
where (a) follows from \eqref{eqn:conc_dissim} and the independence between $|\basesim(1) \cap \basesim(v_0)|$ and $\featrow{1}$; 
and (b) follows from \eqref{eqn:overlap_lower} and the assumption that $\eta - \eta' \geq K \max\Big( \sqrt{\frac{4 \log (m-1)}{c(n-1)p^2}}, \frac{4 \log (m-1)}{c(n-1)p^2} \Big)$.

All in all, 
\begin{align*}
	\Prob{\Nbase = 0 ~\big|~ \featrow{1}}
		&\leq \exp\bigg( - \frac{(m-1)p}{8} \cdot \localprob\Big(\featrow{1}, \eta' - 2\sigma^2\Big) \bigg)
			+ \frac{1}{(m-1)^2} + \exp\Big( - \frac{(n-1)p^2}{8} \Big).
\end{align*}
\end{proof}

\newpage

    

\end{document}